%% file: 2003-11.tex
\documentclass{gtart}
\input gtoutput

\lognumber{245}
\volumenumber{7}\papernumber{11}
\volumeyear{2003}
\pagenumbers{329}{397}
\received{6 April 2002}
\revised{1 June 2003}
\accepted{23 May 2003}
\published{19 June 2003}
\proposed{Cameron Gordon}
\seconded{Walter Neumann, Joan Birman}

\usepackage{amssymb}
\usepackage{latexsym}
\usepackage{enumerate}
\usepackage{epsfig}
\usepackage{subfigure}
\usepackage{amsmath}

\def\orbit{G \centerdot}
\def\a{\orbit A}
\def\b{\orbit B}
\def\g{\gamma}
\def\D{\Delta}
\def\phibar{\overline{\phi}}
\def\R{\mathbb{R}}
\def\s{\sigma}
\def\S{\Sigma}
\def\Z{\mathbb{Z}}
\def\d{\partial}
\def\cross{\times}

\def\fol{\mathcal{F}}

\def\half{\frac{1}{2}}
\def\complement{-}

\newtheorem{theorem}{Theorem}[section]
\newtheorem{lemma}[theorem]{Lemma}
\newtheorem{corollary}[theorem]{Corollary}
\newtheorem{proposition}[theorem]{Proposition}
\newtheorem*{lemma1}{Lemma \ref{lemma1}}
\newtheorem*{lemma2}{Lemma \ref{lemma2}}
\newtheorem*{lemma3}{Lemma \ref{lemma3}}

\theoremstyle{definition}
\newtheorem{definition}[theorem]{Definition}
\newtheorem{notation}[theorem]{Notation}

\newtheorem{remark}[theorem]{Remark}

\newtheorem{example}[theorem]{Example}

\newcounter{case}

\newenvironment{case}[1]{\stepcounter{case} \addvspace{.5\baselineskip} \noindent\textbf{Case \thecase}\qua \textsl{#1}}{\hfill\fbox{Case \thecase}}

\newtheorem{subcase}{Case}[case]

\newcounter{step}

\newenvironment{step}[1]{\stepcounter{step} \addvspace{.5\baselineskip} \noindent\textbf{Step \thestep}\qua \textsl{#1}}{\hfill \framebox[1.1\width]{\raisebox{0pt}[0.8\height][0.1\depth]{Step \thestep}}}

\begin{document}

\title{Period three actions on the three-sphere}

\author{Joseph Maher\\J Hyam Rubinstein}

\address{Mathematics 253-37, California Institute of 
Technology\\Pasadena, CA 91125, USA}

\secondaddress{Department of Mathematics and Statistics\\
University of Melbourne,
Parkville VIC 3010, Australia}

\asciiaddress{Mathematics 253-37, California Institute of 
Technology\\Pasadena, CA 91125, USA\\and\\Department of 
Mathematics and Statistics\\
University of Melbourne,
Parkville VIC 3010, Australia}

\email{maher@its.caltech.edu, rubin@ms.unimelb.edu.au}

\asciiabstract{We show that a free period three action on the
three-sphere is standard, ie, the quotient is homeomorphic to a lens
space. We use a minimax argument involving sweepouts.}

\begin{abstract}
We show that a free period three action on the
three-sphere is standard, ie, the quotient is homeomorphic to a lens
space. We use a minimax argument involving sweepouts.
\end{abstract}

\title{Period three actions on the three-sphere}

\keywords{3-manifold, 3-sphere, group action, spherical 3-manifold, lens space}

\primaryclass{57M60}

\secondaryclass{57M50}

\maketitle

\section{Introduction}

\subsection{Background}

This paper presents a proof that a free period three action on the three-sphere is standard. In this section we give an outline of the proof, but we begin by giving a brief summary of some related results.

Thurston's Geometrization Conjecture implies that all $3$-manifolds with finite fundamental groups are homeomorphic to quotients of $S^3$ by a finite group of isometries. This conjecture splits into the following two conjectures: the Poincar\'e Conjecture, which states that the universal cover of every closed $3$-manifold with finite fundamental group is $S^3$, and the Spherical Spaceform Conjecture, which states that every $3$-manifold whose universal cover is $S^3$ is homeomorphic to a quotient of $S^3$ by isometries.

The $3$-manifolds which are quotients of $S^3$ by isometries are known as spherical or elliptic $3$-manifolds. They have been classified by Hopf \cite{hopf} and Seifert and Threlfall \cite{st}, by identifying the finite subgroups of $SO(4)$ which act freely on $S^3$. Wolf \cite{wolf} has extended this to higher dimensions. There are more recent accounts of the $3$-dimensional case in Scott \cite{scott} and Thurston \cite{thurston}.

Milnor \cite{milnor} and Lee \cite{lee} investigated which finite groups might be able to act freely on $S^3$. They recovered the list of finite subgroups of $SO(4)$ which act freely, but also another family of groups. These extra groups are an infinite family of finite subgroups of $SO(4)$, which do not act freely on $S^3$ as subgroups of $SO(4)$, but may possibly have some non-standard free action on $S^3$.  See table 1.

\begin{table}[ht!]
\begin{center}
\begin{tabular}{p{4in}lll}
Subgroups of $SO(4)$ that act freely & Order \\
Cyclic, $\Z_{n}$ & $n$ \\
Quaternionic, $Q_{4n} \cross \Z_m$, where $n$ and $m$ are coprime & $4nm$ \\
Binary tetrahedral, $T \cross \Z_m$, where $m$ is coprime to 24 & $24m$ \\
Binary octahedral, $O \cross \Z_m$, where $m$ is coprime to 48 & $48m$ \\
Binary icosahedral, $I \cross \Z_m$, where $m$ is coprime to 120 & $120m$ \\
Index two subgroups of $Q_{4n} \cross \Z_{2m}$  &  $4nm$ \\
Index three subgroups of $T \cross \Z_{3m}$ & $24m$ \\
&&\\
Subgroups of $SO(4)$ which could only have non-standard actions &&\\
$Q(8n, k, l)$ & $8nkl$ \\

\end{tabular}
\end{center}
\caption{Finite groups that may act freely on $S^{3}$}
\end{table}

The group $Q(8n, k, l)$ is the semi-direct product of $\Z_{kl}$ with $Q_{8n}$, and has presentation $\{ a,b,c \mid a^2 = (ab)^2 = b^{2n} = 1, c^{kl} = 1, aca^{-1} = c^r, bcb^{-1} = c^{-1} \}$, where $r \equiv -1$ mod $k$, and $r \equiv 1$ mod $l$.  See Milnor \cite{milnor} and Lee \cite{lee} for more details.

By work of Thomas \cite{thomas86}, if the group is solvable, it suffices to show that the actions of its cyclic subgroups are standard. The only groups in table 1 which are not solvable are those containing the binary icosahedral group. In particular, showing that the cyclic actions are standard would eliminate the exceptional groups $Q(8n,k,l)$.

The first group action shown to be standard was $\Z_2$, by Livesay, \cite{livesay}. His argument uses the standard sweepout by $2$-spheres on $S^3$. An invariant unknotted curve is found by considering the intersections of the $2$-spheres by their images under the involution. This was generalised to higher powers of two by Rice \cite{rice}, Ritter \cite{ritter} and Myers \cite{myers}. The other results due to Rubinstein \cite{rubinstein79a, rubinstein79b} use the fact that certain spherical manifolds contain embedded Klein bottles, and involve considering the intersections of the images of the Klein bottle under the group action.

Table 2 summarises which group actions are known to be standard.

\begin{table}[ht!]
\begin{center}
\begin{tabular}{ll}
Group & Reference \\
$\Z_{2}$ & Livesay \cite{livesay} \\
$\Z_{4}$ & Rice \cite{rice} \\
$\Z_{8}$ & Ritter \cite{ritter} \\
$Q_{2^{k}}$ & Evans, Maxwell \cite{em} \\
$\Z_{6}, \Z_{12}, Q_{2^{k}}$ & Rubinstein \cite{rubinstein79a} \\
non-cyclic, non-quaternionic groups of order $2^{a}3^{b}$ & Rubinstein \cite{rubinstein79b} \\
$\Z_{2^{k}}, \Z_{2^{k}.3}, Q_{2^{k}.3}, k \geqslant 2$ & Myers \cite{myers} \\
\end{tabular}
\end{center}
\caption{Group actions known to be standard}
\end{table}

\subsection{Acknowledgements}

The first named author would like to thank his PhD advisor Daryl Cooper for his exceptional fortitude and patience during the writing of this paper. He would also like to thank the National Science Foundation for providing partial support during the writing of this paper. Both authors would like to thank the referee for many helpful comments.

\subsection{Outline}

An action of $\Z_3$ on $S^3$ is given by a diffeomorphism $g\co S^3 \to S^3$, which is free, and which is period three. This means that $g$ generates a group $G \cong \Z_3$, and as $g$ is free, the quotient $S^3/G$ is a manifold. If $g$ is linear, namely an element of $SO(4)$, then the quotient is a lens space. We say the action of $g$ is standard if $g$ is conjugate by a diffeomorphism to an element of $SO(4)$. This is equivalent to the quotient being diffeomorphic to a lens space.

We can show that the action is standard by finding an invariant unknotted circle in $S^3$. An unknotted invariant circle has an invariant neighbourhood which is a solid torus. The complement of this neighbourhood is also an invariant solid torus. The quotient of a solid torus by a finite group acting freely is again a solid torus, so the quotient manifold is the union of two solid tori, a lens space. 

We will find an invariant unknotted curve by studying sweepouts of $S^3$. A sweepout of $S^3$ is a family of surfaces which ``fill up'' the manifold. A simple example is the foliation of $S^3$ by $2$-spheres with two singular leaves which are points. Think of the leaves as parameterised by time, starting with one singular leaf at $t = 0$ and ending with the other one at $t = 1$. At a non-singular time $t$, the sweepout consists of a single $2$-sphere $S_t$. We can think of the union of the leaves as a $3$-manifold, in this case $S^3$, with a height function for which each level set is a $2$-sphere. The map from the leaf space to $S^3$ is degree one, and is an embedding on each level set of the height function.

For our purposes, we require a more general definition, in which the sweepout surfaces at non-singular times may be finitely many spheres. We say a generalised sweepout is a $3$-manifold $M$, with a height function $h$ on it, so that the level sets at regular values are unions of $2$-spheres, and a degree one map $g\co  M \to S^3$, which is an embedding on each level set of the height function. We shall think of the height function on $M$ as time. The $3$-manifold $M$ is in fact either $S^3$, or a connect sum of $S^2 \cross S^1$'s.

We can can look at the three images of the sweepout spheres under the group $\Z_3$. Generically, they will intersect in double curves and triple points. In order to distinguish the three images of the spheres under $\Z_3$, we will colour them red, blue and green.

By general position, we can arrange that the sweepout spheres intersect transversely for all but finitely many times, and that the non-transverse intersections all come from a finite list of possibilities, corresponding to critical points of the height function on either $M$ itself, or on the double or triple point sets of $M$. 
We call these non-transverse intersections {\sl moves}, and we can describe a sweepout by drawing the configurations in $S^3$ in between the critical times, with each neighbouring pair of pictures differing by a non-transverse intersection. 
Critical points of the height function on $M$ change the number of $2$-spheres by one. 
So a $2$-sphere can either appear or disappear, or two spheres can either split apart or join together. We will call the appearance of a sphere an {\sl appear} move, and the disappearance of a sphere a {\sl vanish} move. A critical point which splits two spheres apart will be called a {\sl cut} move, and a critical point which joins two spheres together will be called a {\sl paste} move.
Critical points of the height function on the double set change the number of double curves, by either creating or destroying a double curve, which we will call {\sl birth} or {\sl death} moves respectively, or by saddling curves together, which we will call saddle moves. Critical points of the height function on the triple set change the number of triple points, which we will call {\sl triple point} moves. In fact the number of triple points is always a multiple of six, as every triple point has three images under $G$, and there must also be an even number of triple points.

We look for a sweepout that is ``simple'', by defining a complexity for sweepouts, and showing how to change the sweepout to reduce complexity. We say that the {\sl complexity} of the sweepout at a generic time $t$ is the ordered pair $(n, d)$, where $n$ is the number of triple points, and $d$ is the number of simple closed double curves without triple points. We order the pairs $(n,d)$ lexicographically. We say that the complexity of a sweepout is the maximum complexity that occurs over all generic times. We say that a sweepout is a minimax sweepout if it has minimal complexity. The complexity only depends on the graph of double curves, so moves which do not change the graph of double curves do not change complexity.

At a generic time, the double curves form an equivariant graph in $S^3$, which may contain an invariant unknotted curve. The basic strategy is to show that a minimax sweepout must have an unknotted invariant curve in its graph of double curves at some time during the sweepout. This follows if we can show that an arbitrary local maximum can either be reduced in height by changing the sweepout, or else contains an invariant unknotted curve in its graph of double curves. A minimax sweepout contains local maxima that cannot be removed, so it contains an unknotted invariant curve.

We now give a more detailed outline of the argument.

\subsection{Modifications}

In order to apply this strategy we need some way to simplify a sweepout with many triple points. One of the basic operations we will use is to change the sweepout by cutting out a subset of the sweepout, and replacing it with a different subset. We now give an informal description of this procedure, which we will call a {\sl modification}. A precise description of this is given in Section \ref{subsection:modifications}.

Imagine choosing an equivariant $3$-dimensional subset $N$ of $S^3$, and watching the images of the sweepout surfaces inside it for some time interval $I$. Assume that the boundary of $N \cross I$ is disjoint from any of the moves of the sweepout. We can think of the image of the sweepout in $N \cross I$ as a sequence of pictures of spheres and planar surfaces in $N$, with each picture in the sequence differing from its neighbours by a move. If we draw a sequence of pictures in $N$ which is different, but which agrees with the original one on the boundary of $N \cross I$, then we can try to construct a sweepout by replacing the original map into $N \cross I$ by a new one determined by the new set of pictures. We need to check that the new map really is a sweepout, by checking that the level sets are still spheres, and that the projection to $S^3$ is still degree one, but if it is a sweepout, then we say that we have produced a new sweepout from the original one by a modification.

The two main modifications we will use will involve removing double curves, and removing bigons. A {\sl bigon} is a subdisc of a sweepout surface whose boundary consists of a pair of double arcs which share common endpoints, and whose interiors are disjoint. We say a subsurface of a sweepout sphere is {\sl double curve free} if its interior is disjoint from all of the double curves.

\begin{description}

\item[Removing double curves] 

Suppose we have a double curve that bounds a double curve free disc $\D$ in a sweepout sphere, for some time interval $I$. 
Think of the sweepout surface containing the disc as ``horizontal'', and the other sweepout surface intersecting the double curve as ``vertical''.
Then we can change the sweepout to remove this double curve for a subinterval of $I$, by ``pinching off'' the vertical annulus at the beginning of the time interval, and then replacing it at the end of the time interval, so that the sweepout returns to its initial configuration.
We want to do this generically, using moves, so first use a cut move to pinch the vertical annulus into two discs, one of which is disjoint from the horizontal surface, and the other intersects it in the double curve.
Now use a double curve death move to remove the double curve. We will call this special pair of moves which removes a double curve a compound double curve death, and its time reverse a compound double curve birth.
This removes the double curve, and also splits one of the spheres into two.  
We can also do this equivariantly, as the disc $\D$ is disjoint from its images, so in fact we will remove the orbit of a double curve, and split three spheres apart. This is illustrated in Figure \ref{figure:removing a double curve}.

\item[Removing bigons] 

Suppose we have a bigon $A$, which is double curve free for a time interval $I$. We can reduce the number of triple points in the sweepout for a subinterval of $I$ by ``undoing'' the bigon. Think of the bigon as being ``horizontal'', and the sweepout surfaces which intersect it in the double arcs of the boundary of $A$ as being ``vertical''. We can push the two vertical surfaces sideways through each other to remove the bigon, and reduce the number of triple points. As before, we wish to do this generically, and this in fact requires a saddle move followed by a triple point death move. We will call the pair of moves that remove the bigon a compound triple point death, and we will call their time reverse a compound triple point birth. We can also also do this equivariantly, as a double curve free bigon is disjoint from its images. This is illustrated in Figure \ref{figure:removing a bigon}.

\end{description}

Both of these modifications are in fact guaranteed to produce new sweepouts. Removing bigons only changes the level sets up to isotopy, while removing a double curve splits one $2$-sphere into two spheres, and then joins them back together again. The degree of the map into $S^3$ is also preserved, because as the modification neighbourhood $N_t$ is not all of $S^3$, (perhaps after adjusting the sweepout by an isotopy) we can find a point in $S^3$ which does not lie in any of the $N_t$, so the pre-image of this point does not change.

We can start using these modifications to try to simplify the sweepout. For example, suppose there is a local maximum in the graph of complexity against time, for which there is a double curve free bigon which is disjoint from all the moves occurring during the local maximum. Then we can remove the bigon for the duration of the local maximum, which reduces the height of the local maximum.

\subsection{The main argument}

Suppose we have a time interval which contains a local maximum of complexity. This will contain a move which increases complexity, immediately followed by a move which decreases complexity. If these two moves have disjoint move neighbourhoods, then we can just swap the order in which they occur,  reducing the height of the local maximum. In general, if we change the sweepout in any way so that the height of a local maximum is reduced, then we say that the local maximum has been {\sl undermined}.

The main argument, in Section \ref{section:undermining}, shows how to undermine an arbitrary local maximum, assuming the following lemmas:

\begin{lemma1}
Every sweepout contains triple points.
\end{lemma1}

We say that a bigon is {\sl vertex free} if it contains no triple points in its interior, and precisely two triple points in its boundary.

\begin{lemma2}
Suppose a non-triple point move occurs while the configuration contains triple points. Then there is a vertex free bigon orbit disjoint from the move neighbourhood of the move.
\end{lemma2}

We say that a local maximum which consists of a compound triple point birth, followed by a compound triple point death, is a {\sl special case local maximum}. There are only finitely many different ways in which the two compound moves can intersect.

\begin{lemma3}
A special case local maximum can either be undermined, or else contains an invariant unknotted circle.
\end{lemma3}

The first step of the argument in Section \ref{section:undermining}, is to change the sweepout, so that the highest local maxima consist of compound moves only. Lemma \ref{lemma2} implies that for every non-triple point move there is a vertex free bigon which is disjoint from the move. If this contains simple closed double curves, then there is an innermost one which we can remove using a double curve removal modification. If there are no double curves inside the bigon, then we can remove the bigon using a bigon removal modification. Both of these modifications reduce the height at which the non-triple point move occurs, and insert extra compound moves into the sweepout. We can use Lemma \ref{lemma2} repeatedly until all non-triple point moves lie far below the height of the local maxima. We then show that we can replace a non-compound triple point move with a compound triple point move, possibly adding a non-triple point move at a lower height than the local maxima. In this way we may assume that all moves above  a certain height are compound moves, so it suffices to show how to undermine local maxima in which all the moves are compound moves.

The next step is to show that we can undermine all local maxima consisting of compound moves only, which are not special case local maxima, ie, they contain at least one compound double curve move. If both compound moves are double curve moves which are disjoint, then we can swap the order of the moves to undermine the local maximum. If they are not disjoint, then we show that we can find a disjoint double curve or bigon to remove for the duration of the local maximum. If one of the compound moves is a double curve move and the other is a triple point move, then we can swap the order, as the bigon for the triple point move must be disjoint from the disc for the double curve move.

Lemma \ref{lemma3} now implies that either we can remove all the special case local maxima, or else there is an invariant unknotted curve. Lemma \ref{lemma1} ensures that we can not remove all the triple points from the sweepout, so we cannot undermine all the special case local maxima, so there must be at least one that contains an unknotted invariant curve.

In order to complete the proof it remains to prove the three lemmas, and we now give a brief summary of the arguments that we will use.

\subsection{Every sweepout contains triple points}

We can choose a continuously varying ``inside'' for the sweepout spheres $S_t$. As $g$ is degree one, we can choose the inside so that it starts out small and ends up large. The three images of $S_t$ under $G$ are coloured, and we colour the inside of each image of $S_t$ with the same colour as the surface. So each component of $S^3 - \orbit S_t$ may be coloured by some combination of the colours red, blue and green. If a region is outside of all of the spheres, then we say that it is a clear region.

Before any of the sweepout spheres have appeared, $S^3$ is a clear connected invariant region. After the sweepout, $S^3$ is coloured with all three colours, and there are no clear regions at all. So some move must break up the clear connected invariant region. We show that if there are no triple points, none of the non-triple point moves can do this, using a simple case by case argument.

\subsection{Disjoint bigons}

We show that if there are bigons, then there is always a bigon disjoint from a non-triple point move. We use an elementary combinatorial argument to show this.

\subsection{Special cases}

A special case local maximum consists of a compound triple point birth, followed  by a compound triple point death. The compound birth move creates a bigon, $A$ say, and the compound death destroys a bigon, $B$ say. Between the compound moves, both bigons are present at the same time. If the two bigon orbits are disjoint, then we can undermine the local maximum by just swapping the order in which the compound moves occur. If the two bigon orbits are the same, then the configuration before the local maximum is isotopic to the configuration after the local maximum, so we can just replace the moves creating the local maximum with an isotopy, thus undermining the local maximum. If they are not disjoint, and not the same, then there are only finitely many different ways in which they can intersect, and we deal with each possibility in turn.

For example, if the bigon orbits share a single vertex orbit in common, then we explicitly construct a new sweepout with fewer triple points. Otherwise, the two bigon orbits share all six vertices in common. If the union of the bigon orbits is connected, then we show that there is an invariant unknotted curve. In the remaining cases, we prove a useful undermining lemma, which shows that if we have a $3$-ball disjoint from its images under $G$, which intersects the sweepout in a ``simple'' way, then we can replace the sweepout inside the $3$-ball with one which has no triple points. We then show how to find such a $3$-ball which contains the compound moves, for each of the remaining special cases.

We now give a brief description of the undermining lemma, and explain the ``simple'' condition on the intersection of the sweepout surface with the $3$-ball. Assume we have chosen a $3$-ball $B$ which is disjoint from its images under $G$, which contains all the compound moves, and contains no triple points before or after the local maximum. The intersection of the sweepout surfaces with the boundary of $B$ is constant up to isotopy, and consists of a pattern of intersecting circles, coloured red, blue and green, which we shall call the boundary pattern. If there is a simple closed curve of intersection which bounds a double curve free disc in $\d B$ then we can use this disc as a cut disc for a cut move, to reduce the number of curves of intersection between the sweepout surfaces and $\d B$. If the curves of intersection create a double curve free bigon in the boundary, then we can use the bigon as a saddle disc for a saddle move which reduces the number of intersections of the double curves with $\d B$. If we can remove all the intersections of the sweepout surfaces with $\d B$ in this way, then we say that the boundary pattern is {\sl saddle reducible}. This is what we mean by ``simple'' intersections.

We can construct a partial sweepout which agrees with the original one on $\d B$. Start with the configuration in $B$ before the local maximum, and now do cut moves and saddle moves in some order to remove all the intersections of the sweepout surface with a $2$-sphere parallel to $\d B$, but just inside $B$. The remaining components of the sweepout surfaces inside $B$ can now be removed using cut and death moves, to remove the double curves, and vanish moves to remove $2$-sphere components. We can construct a similar partial sweepout starting with the configuration after the local maximum, and working back in time, using the same sequence of cut and saddle moves near the boundary. These two partial sweepouts can then be patched together in the middle by an isotopy to create the desired replacement partial sweepout.

\section{Preliminary definitions} \label{section:definitions}

\subsection{Generalised sweepouts}

\begin{definition}{\sl Free action} 
\index{free action} \index{action!free}

We say that a {\sl free action} of the group $G \cong \Z_3$ on $S^3$ is generated by the diffeomorphism $g\co S^3 \to S^3$ if $g$ has no fixed points, and $g^3$ is the identity. Therefore $f$ generates a cyclic group of order three, $G = <\!g\!> \cong \Z_3$.
\end{definition}

This means that $g$ is orientation preserving.

\begin{remark}
Throughout this paper, we assume that we have been given some fixed $g$, which we do not change.
\end{remark}

\begin{notation}

If $N$ is a subset of $S^3$, then we write $\orbit N$ for the orbit of $N$ under $G$. As $g$ is free, $S^3/G$ is a manifold, we will call this manifold $L$.
\end{notation}

\begin{definition}{\sl Standard action} 
\index{standard action} \index{action!standard}

We say that the action is {\sl standard} if $g$ is conjugate by a diffeomorphism to a linear map, ie, an element of $SO(4)$. 
\end{definition}

This is equivalent to $L$ being diffeomorphic to the quotient of $S^{3}$ by a linear group, in this case a lens space. The following theorem is well known:

\begin{theorem}{}
If there exists a smooth unknotted invariant circle in $S^3$, then the action of $G$ is standard.
\end{theorem}

\begin{proof}
Given an unknotted invariant circle, we can find an invariant tubular neighbourhood, which is a solid torus. As the curve is unknotted, the complement of the tubular neighbourhood is also a solid torus. So if we can show that a free action on a solid torus gives a solid torus, then this gives a genus one Heegaard splitting of the quotient, which is therefore a Lens space. We now show that a free $\Z_3$ action on a solid torus is standard.

The quotient manifold has torus boundary, as the boundary of the solid torus is a torus, and the quotient of a torus by a free $\Z_3$ action is also a torus. The quotient manifold is also irreducible, as it is covered by a solid torus, which is irreducible. A meridional disc in a solid torus is a compressing disc for the boundary, and this projects down to an immersed disc in the quotient, so the quotient manifold also has compressible boundary. By the loop theorem, there is an embedded compressing disc in the quotient. If we cut the quotient along this disc, we get a manifold with $2$-sphere boundary, which bounds a $3$-ball, as the quotient is irreducible. Therefore the quotient manifold is a $3$-ball with a single one-handle attached, ie, a solid torus.  
\end{proof}

\begin{definition}{\sl A generalised sweepout} 
\index{sweepout} \index{generalised sweepout} \index{sweepout!generalised} \index{$M$} \index{$f$} \index{$h$} \index{$\phi$}

A {\sl generalised sweepout} is a triple $(M, f, h)$, where

\begin{itemize}

\item $M$ is a closed, orientable $3$-manifold. 

\item The smooth map $h\co M \to \R$ is a Morse function, such that for all but finitely many $t \in \R$, the inverse image, $h^{-1}(t)$ is a collection of $2$-spheres.

\item The smooth map $f\co M \to S^{3}$  is degree one. 

\item The map $f|_{h^{-1}(t)}$ is an embedding on the level set $h^{-1}(t)$ for every $t$.

\end{itemize}

We will often write a sweepout as $(M, \phi)$, where $\phi$ denotes the map $(f \cross h)\co M \to S^3 \cross \R$. We will think of $t \in \R$ as the time coordinate.
\end{definition}

\begin{remark} \index{$\pi_{S^3}$} \index{$\pi_\R$}
The map $\phi\co M \to S^3 \cross \R$ is a smooth embedding. 
We will write $\pi_{S^3}$ and $\pi_\R$ for the projection maps from the product $S^3 \cross \R$ to its factors. 
From now on, whenever we say ``sweepout'', we mean ``generalised sweepout''.
\end{remark}

\begin{notation}
\index{$M_t$} \index{$S_t$}

We will write $M_t$ for $h^{-1}(t)$. The $M_{t}$ form the leaves of a singular foliation $\fol$ of $M$.
We will write $S_t$ for $f(M_t)$, and call these the sweepout surfaces at time $t$.
\end{notation}

\begin{remark}
The manifold $M$ need not be connected. In fact each component of $M$ will be either $S^3$ or a connect sum of $S^2 \cross S^1$'s. We prove this later on as Proposition \ref{proposition:connect_sum}.
\end{remark}

\begin{example}
Let $M = S^{3}$ be the set of points a unit distance from the origin in $\R^{4}$, with the Morse function $h$ given by the $x$-coordinate. Let $f \co M \to S^3$ be the identity map. Then $(S^3, f \cross h)$ is a sweepout.
\end{example}

This example gives a foliation of $S^{3}$ by $2$-spheres, with two singular leaves consisting of points. We will think of the leaves of this foliation as being parameterised by a time interval. Time starts at $t=-1$ at one singular leaf, and ends at $t=1$ at the other. So in the standard sweepout $S_{-1}$ is a single point. As $t$ increases, $S_t$ becomes a sphere which moves away from the initial point, increasing in size, till at $t=0$, $S_t$ is an equatorial sphere. It then decreases in size, until at $t=1$, $S_1$ is again a single point, in fact the antipodal point to $S_{-1}$.

\begin{definition}{\sl Morse function on a manifold with boundary}

Let $M$ be a manifold with boundary, and let $h\co M \to \R$ be a function, which extends to a Morse function on $M$ union an open collar neighbourhood of $\d M$. There should be no singularities on $\d M$. Then $h$ is a Morse function on $M$.
\end{definition}

\begin{definition}{\sl Compatible product structures} \index{compatible product structure}\index{product structure!compatible} \index{$N_I$}

Let $K$ be a compact equivariant submanifold of $S^3 \cross \R$, for which there is a level preserving diffeomorphism $\psi\co  S^3 \cross \R \to S^3 \cross \R$, so that $\psi(K)$ is a product $N \cross I \subset S^3 \cross \R$. We will always choose $N$ to be $3$-dimensional. Then we say that $K$ has a compatible product structure. We can think of $K$ as a family of equivariant submanifolds of $S^3$ that varies continuously with time.  We will write $K$ as $N_I$, and we will write $N_t$ to refer to $K \cap  ( S^3 \cross  t )$.
\end{definition}

\begin{figure}[ht!]
\begin{center}
\epsfig{file=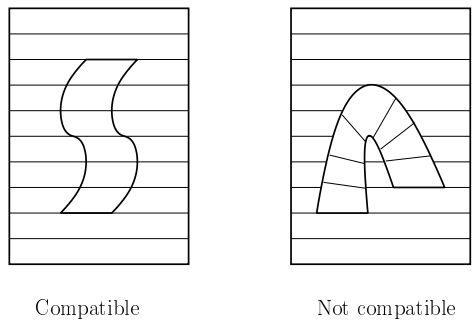, height=120pt}
\end{center}
\caption{Compatible and incompatible product structures}
\end{figure}

The submanifold $K$ need not be connected. Note that $K$ need not be the product of a subset of $S^3$ with an interval in $\R$, it just has to be isotopic to one under a level-preserving isotopy of $S^3 \cross \R$.

We will often be interested in how a compatible product submanifold intersects the sweepout surfaces.

\begin{definition}{\sl Constant boundary pattern}
\index{boundary pattern!constant} \index{constant boundary pattern}

Let $(M, \phi)$ be a sweepout, and let $N_I$ be a compatible product submanifold of $S^3 \cross \R$. If the intersection of the sweepout surfaces with $\d N_t$ only changes by an isotopy during the time interval $I$, then we say that $N_I$ has {\sl constant boundary pattern}, with respect to the sweepout $(M, \phi)$.
\end{definition}

\begin{definition}{\sl A partial sweepout}
\index{partial sweepout} \index{sweepout!partial}

Let $N_I$ be a compatible product submanifold of $S^3 \cross \R$, where $N$ is three dimensional. Let $P$ be a $3$-manifold with boundary, which need not be connected, and let $(f \cross h) \co P \to N_I$ be a proper embedding so that $h$ is a Morse function and $h^{-1}(t)$ is a union of spheres and planar surfaces for all but finitely many $t$. If the surface $h^{-1}(t)$ has singularities, they should be disjoint from $\d P$. 

Then $(P, f \cross h)$ is a {\sl partial sweepout} with {\sl image neighbourhood} $N_I$.
\end{definition}

\begin{remark}
An image neighbourhood has constant boundary pattern, as there are no singularities on the boundary.
\end{remark}

\begin{definition}{\sl A partial sweepout which is a product}
\index{product partial sweepout}\index{partial sweepout!product}

Let $(P, \phi)$ be a partial sweepout, in which $P$ has a product structure, so that the image of the product structure on $\phi(P)$ is compatible with the product structure on $S^3 \cross \R$. Then we say that $(P, \phi)$ is a {\sl product partial sweepout.}
\end{definition}

\subsection{Diagram Conventions}

\begin{notation}
The sweepout surface $S_t$ is a union of $2$-spheres, and we will
label them red. The spheres in $S_t$ have two sets of images under
$G$, we will label the spheres in $gS_t$ green, and the spheres in
$g^2S_t$ blue. If two spheres intersect in a double curve, we will
label the double curve by the complementary colour, ie, the colour of
the sphere not involved in the intersection. For example, if the red
sphere intersects the green sphere in a double curve, we will label
that double curve blue.  To distinguish colours when printed in black
and white (or viewed on a b/w monitor) we have used solid lines for
green, long dashes for blue and short dashes for red.

\end{notation}

\begin{definition}{\sl The outside of a sphere}
\index{sphere!outside}

At a non-singular time $t$ we can choose a normal vector for one of the disjoint red spheres. As each $2$-sphere is separating in $S^3$, we can choose compatible normal vectors for the other spheres in $S_t$, so that in each complementary region the normal vectors all point either in or out. As $t$ varies, we can choose a continuously varying family of such normal vectors defined on $S^3 -\{  \text{singular points} \}$. 
We define the {\sl outside} of the family of spheres to be the side that the normal vector points to, and we choose the normal vector so that the outside starts out large and ends up small.
\end{definition}

\begin{notation}
The inside of the sphere will be labelled with the same colour as the sphere. Double curves occur in threes, one of each colour, and generically intersect spheres of the same colour transversely.  
\end{notation}

\begin{definition}{\sl A region}
\index{region}

A {\sl region} (at time $t$) is the closure of a connected component of the complement of the orbit of $S_t$ in $S^3$.
\end{definition}

A region may be coloured by any combination of the colours red, green and blue, so there are eight different ways in which a region may be coloured.

\begin{definition}{\sl A clear region}
\index{clear region} \index{region!clear}

If a region is not coloured by any colour, we say that it is a {\sl clear region}.
\end{definition}

In order to draw diagrams of the intersections of the spheres, it is helpful to think of $S^{3}$ as $\R^3 \cup \{ \infty \}$. There is an isotopy of $S^{3}$ that takes one of the red spheres to the $xy$-plane, so that the normal vector points upwards in the direction of the positive $z$-axis. 

\begin{definition}{\sl The blue-green diagram and graph}
\index{blue-green diagram} \index{diagram!blue-green}
\index{blue-green graph} \index{graph!blue-green}

A {\sl blue-green diagram} consists of finitely many red spheres, which may contain finitely many blue and green circles. The union of the green circles is embedded, as is the union of blue circles. The green circles may intersect the blue circles transversely to form a four-valent graph, which we call the {\sl blue-green graph}. 
\end{definition}

\begin{definition}{\sl The configuration}
\index{configuration}

We will refer to the image of the three families of spheres in $S^3$, or some subset of $S^3$, as a {\sl configuration}. 
\end{definition}

Different configurations may give rise to the same blue-green diagram. 
The blue-green graph may have many connected components in each red sphere.
The vertices of the blue-green graph correspond to triple points, and are four-valent in the blue-green graph. They are six-valent in the graph of double curves in the $3$-dimensional configuration.

\begin{definition}{\sl A red bigon}
\index{bigon}

A {\sl red bigon} is a closed disc in a red sweepout sphere $S_t$, whose boundary consists of a pair of simple arcs, one green and one blue, which contain no triple points in their interiors. The endpoints of the arcs are a pair of triple points.

A green bigon is the image of a red bigon under $g$, and a blue bigon is the image of a red bigon under $g^2$. 
\end{definition}

\begin{definition}{\sl A bigon orbit}
\index{bigon-orbit}

A {\sl bigon orbit} $\a$ is the orbit of a red bigon $A$.
\end{definition}

A bigon orbit is equivariant, and consists of a red bigon, a green bigon and a blue bigon. 
The interior of a red bigon may intersect its images if it contains double curves in its interior. However the boundary of a red bigon is disjoint from its images under $G$, see Corollary \ref{corollary:bigon orbit disjoint boundaries}.
There may be triple points in the interior of the bigon orbit, if the red bigon contains more than one component of the blue-green graph.

\begin{definition}{\sl Vertex free}
\index{vertex free}

We say a subsurface of the sweepout spheres is {\sl vertex free} if it contains no triple points in its interior, and there are no double arcs that cross the boundary transversely.
\end{definition}

A vertex free subsurface of the sweepout spheres may still contain double curves without triple points in its interior.

\begin{definition}{\sl Double curve free}
\index{double curve free}

We say a subsurface of the sweepout spheres is {\sl double curve free} if its interior is disjoint from the double curves.
\end{definition}

\begin{remark}
Double curve free implies vertex free.
\end{remark}

\subsection{General position}

\begin{definition}{\sl Singular sets}
\index{double set} \index{triple set}

Let $(M, f \cross h)$ be a generalised sweepout, and let $g\co S^3 \to S^3$ be a diffeomorphism of period three. Let $\bar f = f/G$, and let $\phibar = \bar f \cross h\co M \to L \cross \R = (S^3/G) \cross \R$. The {\sl double set}, $\S_2$, of $M$ is $\{ x \in M \mid \phibar^{-1}(\phibar(x))$ contains at least two points $\}$. The {\sl triple set}, $\S_3$, of $M$ is $\{ x \in M \mid \phibar^{-1}(\phibar(x))$ contains at least three points$\}$.
\end{definition}

As $\phibar$ is a map from a $3$-manifold to a $4$-manifold, $\S_2$ will be a $2$-dimensional manifold and $\S_3$ will be a $1$-dimensional manifold,  if $\phibar$ is self-transverse.

\begin{definition}{\sl General position for a generalised sweepout with respect to~$g$}
\index{sweepout!general position}

Let $(M, f \cross h)$ be a generalised sweepout, and let $g\co S^3 \to S^3$ be a period three diffeomorphism.

We say that the sweepout $(M, f \cross h)$ is in {\sl general position with respect to $g$}, if $h$ is a Morse function when restricted to $\phibar(\S_2)$ and $\phibar(\S_3)$. Furthermore, we require the critical times of the Morse functions on $M$ and the singular sets to all be distinct. 
\end{definition}

\begin{remark}
We require $h / G$ to be a Morse function on the \emph{image} of the double set, as critical heights on the image will always have pairs of singularities in the pre-image in $M$.
\end{remark}

\begin{theorem}
Let $(M, f \cross h)$ be a generalised sweepout, and let $g$ be a period three diffeomorphism of $S^3$. Then we can alter the map $f \cross h$ by a small homotopy so that $(M, f \cross h)$ is a sweepout in general position with respect to $g$.
\end{theorem}

\begin{proof}
The map $\phibar$ is a smooth immersion, so we can change $\phibar$ to $\phibar'$ by a small homotopy so that $\phibar'$ is a self-transverse immersion. The singular sets $\S_2$ and $\S_3$ will now be nested submanifolds of $M$. A small homotopy of $\phibar = \bar f \cross h$ is the same as small homotopies of $\bar f$ and $h$, so we may assume that $h'$ is a Morse function, as Morse functions are dense in $C^\infty(M, \R)$.

We can now perturb the product structure on $L \cross \R$ to make coordinate projection onto $\R$ a Morse function, not just on $M$, but on the singular sets as well. Note that perturbing the product structure changes the height function on $M$, without changing the singular sets. One way to construct such a perturbation is to choose an embedding of $L$ in $\R^k$, for some $k$, and then use this to define an embedding $e\co L \cross \R \to \R^k \cross \R$. Let $\pi_i$ be projection on to the $i$-th coordinate of $\R^{k+1}$, and consider the functions $e_a$ on $M$ defined by $e_a = \pi_{k+1} \circ e \circ \phibar' + a_1 \pi_1 + \cdots + a_{k+1} \pi_{k+1}$, for $a \in \R^{k+1}$. The set of $a \in \R^{k+1}$ for which $e_a$ is a Morse function on the images of $M$ and the singular sets is an open dense set of $\R^{k+1}$, as they are immersed manifolds in $\R^{k+1}$. Therefore we can choose $a$ to be close to zero, and furthermore $\bar f' \cross h'$ is homotopic to $\phibar'' =  \bar f' \cross (h' + e_a)$ using the homotopy $\bar f' \cross (h' + e_{ta})$ for $t \in [0,1]$.

We have changed $\phibar$ to $\phibar''$ by a small homotopy, so that $\phibar''$ is transverse, and $h''$ is a Morse function on $M$ and the singular sets. As $f''$ is homotopic to $f$, $f''$ is degree one. It remains to check that the level sets are still $2$-spheres. We changed $h$ to $h''$ by a small $C^{\infty}$ perturbation, so not only will ${h''}^{-1}(t)$ lie in a small product neighbourhood of $h^{-1}(t)$, but we have also changed the first derivative by only a small amount. This means that the projection map ${h''}^{-1}(t) \to h^{-1}(t)$ coming from the product structure on a small tubular neighbourhood of $h^{-1}(t)$, is a local diffeomorphism, hence a covering map. Therefore ${h''}^{-1}(t)$ is also a $2$-sphere.
\end{proof}

For each critical point, we can choose a small partial sweepout which contains it, which is isotopic to a ``standard model'' for that type of critical point. We will call such a small partial sweepout a {\sl move}. We now give a precise definition of a move, and then list all the critical points that may arise, and describe their ``standard models''.

\begin{definition}{\sl Moves}
\index{moves}

A move is a partial sweepout whose image neighbourhood is a tubular neighbourhood for the orbit of a critical point. We may assume that the image neighbourhood has a compatible product structure $N_I$, where each $N_t$ is the orbit of a $3$-ball which is disjoint from its images. Furthermore, we require the intersection of the sweepout with the image neighbourhood to be ``as simple as possible''. This means that the sweepout surfaces in $N_I$ must be isotopic to the explicit descriptions of the configurations in $N_I$ which we give below.

We call $N_I$ the {\sl move neighbourhood}. We call $I$ the {\sl move interval}.
\end{definition}

The following is a complete list of the moves that may occur, and a description of the images of the sweepout surfaces in the move neighbourhood. In each case we draw pictures of the sweepout surfaces in one of the components of $N_t$. The other components are the disjoint images of this under $G$.

\begin{description}

\item[Critical points of the height function on $M$] \hbox{}$\phantom {ZZ}$

\begin{enumerate}

\item {\bf Appear and vanish moves} (Index 0 and 3)
\index{appear} \index{vanish} \index{move!appear} \index{move!vanish} \index{$A$} \index{$V$}

In an appear move, a sphere orbit appears. In a  vanish move a sphere orbit disappears. At the critical time the singular set consists of the orbit of a single point. Passing through the critical time a $2$-sphere either appears or vanishes.

\begin{figure}[ht!]
\begin{center}
\epsfig{file=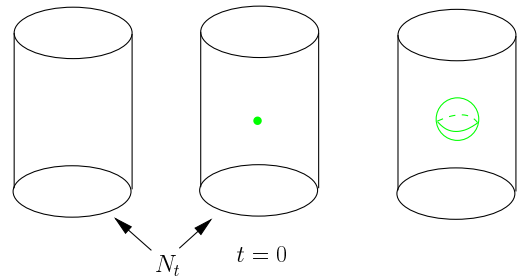, height=105pt}
\end{center}\vspace{-10pt}
\caption{An appear move}
\end{figure}

The sweepout spheres are disjoint from $\d N$. Before the singular time, $N$ is disjoint from the sweepout spheres. After the singular time, each component of $N_t$ contains a $2$-sphere.

This can be modelled by $x^2 + y^2 + z^2 = t$, for $t \in [-1, 1]$. The singular time is $t = 0$. The time reverse of this is a vanish move.

\item {\bf Cut and paste moves} (Index 1 and 2)
\index{cut} \index{paste} \index{move!cut} \index{move!paste} \index{$C$} \index{$P$}

In a cut move, a sphere orbit splits into two. In a paste move, two sphere orbits join together. At the singular time two sphere orbits share a single point in common. Each of the images of this pair of spheres also has a common point.

\begin{figure}[ht!]
\begin{center}
\epsfig{file=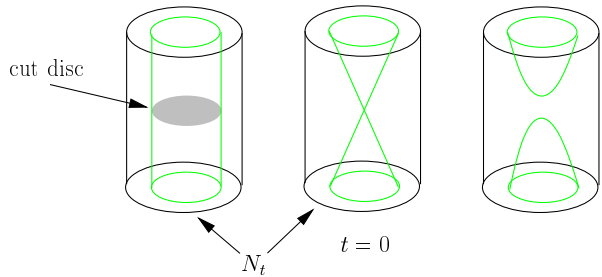, height=105pt}
\end{center}
\caption{A cut move}
\end{figure}

The green sphere intersects $\d N$ in a pair of circles. Before the singular time, the pair of circles bound a properly embedded annulus in $N$. After the singular time the circle bound a pair of disjoint discs.

This can be modelled by $z = x^2 + y^2 + t$, for $t \in [-1, 1]$. The singular time is $t = 0$. The time reverse of this is a paste move.

\end{enumerate}
\end{description}

\begin{definition}{\sl Cut disc}
\index{cut disc}

A cut move can be specified by giving an embedded disc in $S^3$, whose boundary is contained in a single sweepout sphere, whose interior is disjoint from the sweepout spheres, and which is disjoint from its images under $G$. We call this the {\sl cut disc}. 
\end{definition}

\begin{description}
\item[Critical points of the height function on the double set] 

The double set is $2$-dimensional, so the height function has three sorts of critical points, births, deaths and saddles.

\begin{enumerate}

\item {\bf Birth/death of double curves} (index 0 and 2)
\index{double curve!birth} \index{double curve!death} \index{move!double curve birth} \index{move!double curve death} \index{$B$} \index{$D$}

At the critical time, two spheres of different colours intersect in a single point in each of the components of $N_t$.

\begin{figure}[ht!]
\begin{center}
\epsfig{file=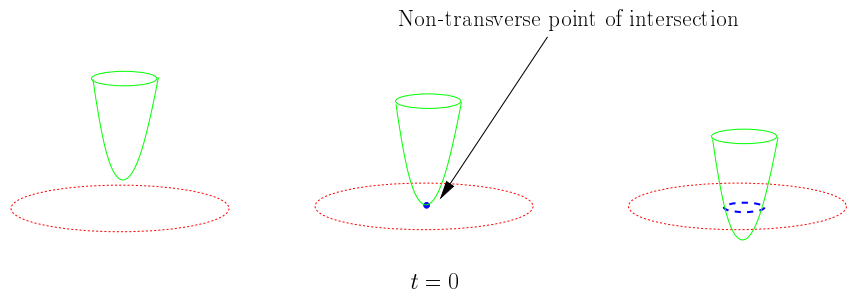, height=110pt}
\end{center}
\caption{Birth of a blue double curve}
\end{figure}

The sweepout surfaces intersect one of the components of $\d N_t$ in a red circle and a green circle, which bound a red and a green disc inside that component of $N_t$, respectively. Before the singular time, the two discs are disjoint. At the singular time the two discs intersect at a single point. After the singular time the two discs intersect in a single (blue) double curve.

This can be modelled by choosing the green surface to be $z = x^2 + y^2$, and the red surface to be $z = t$, for $t \in [-1,1]$. The singular time is $t = 0$. A death move is the time reverse of this.

\item {\bf Saddle moves} (index 1)
\index{saddle} \index{move!saddle} \index{$S$}

During a saddle move, either a single double curve orbit splits into two double curve orbits, or two double curve orbits join together to form a single double curve orbit. At the critical time, each of the components of $N_t$ contains the singular point of a figure eight double curve of intersection between two spheres of different colours.

\begin{figure}[ht!]
\begin{center}
\epsfig{file=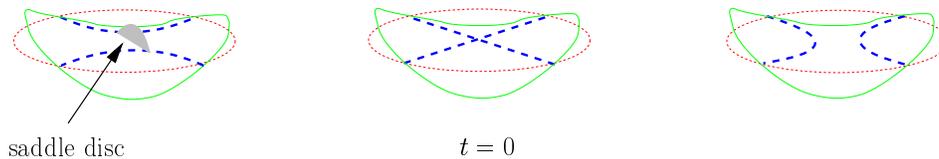, height=60pt}
\end{center}
\caption{A saddle move} \label{figure:saddle}
\end{figure}

The component of $N_t$ shown above contains two properly embedded discs, one of which is red and the other is green. The boundaries of the discs are two circles in $\d N_t$ which intersect at four distinct points. The points of intersection occur in the same order on each circle, so each point is adjacent to the same two points, whether you travel along the red circle or the green circle.

Before the singular time, the green disc intersects the red discs in a pair of double arcs, which connect two pairs of adjacent points in this component of $\d N_t$.

After the singular time, the discs intersect in a pair of double arcs that connect each point to the other adjacent point in $\d N_t$.

This can be modelled by choosing the green surface to be $z = x^2 - y^2$, and the red surface to be $z = t$, for $t \in [-1, 1]$. The singular time is $t = 0$. The time reverse of this is also a saddle move.

\end{enumerate}
\end{description}
\begin{definition}{\sl Saddle disc}
\index{saddle disc}

The saddle move can be specified by giving the orbit of a disc which is disjoint from its images, properly embedded in $(S^3, \orbit S_t)$, whose interior is disjoint from the sweepout surfaces, and whose boundary consists of two arcs, one of which lies in the red sphere, and the other lies in the green sphere, which meet at the blue double arcs. This is illustrated in Figure \ref{figure:saddle}. The saddle move can be thought of as pushing the green sphere through the red sphere along this disc.
\end{definition}

\begin{description}

\item[Critical points of the Morse function on the triple point set] \hbox{}
$\phantom{55}$\medskip

{\bf Triple point births and deaths} (index 0 and 1)\qua
\index{triple point move} \index{move!triple point} \index{$T \pm$}
In a triple point move, a pair of triple points is either created or destroyed. The triple point set is one dimensional, so there are two sorts of singularities, births and deaths of triple points.

\begin{figure}[ht!]
\begin{center}
\epsfig{file=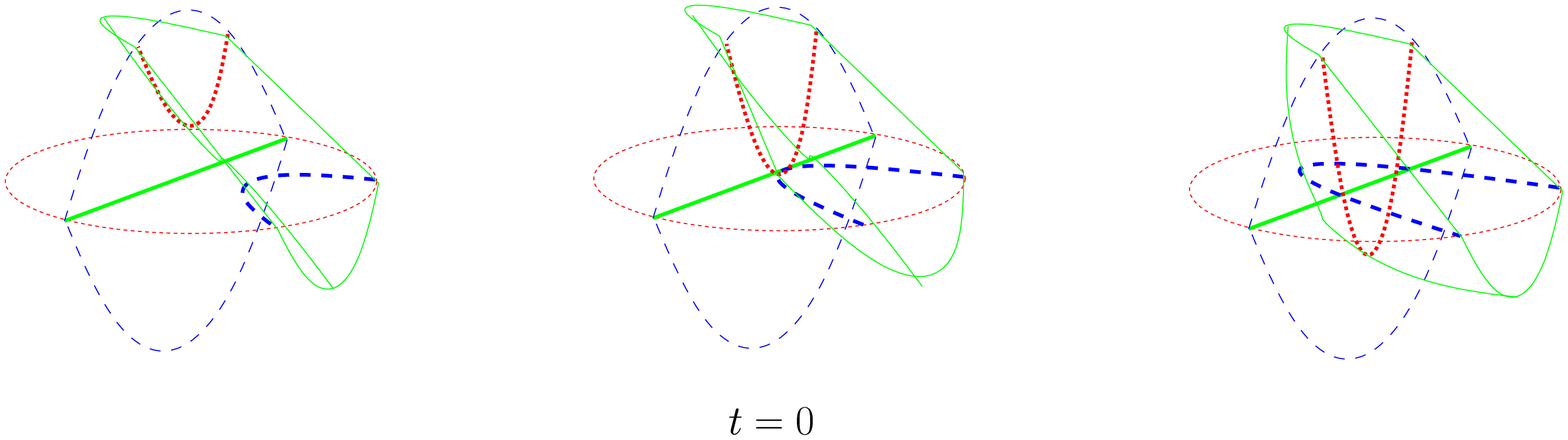, height=100pt}
\end{center}
\caption{A triple point birth}
\end{figure}

The component of $N_t$ shown above contains three properly embedded disc, one of each colour. The boundary of each disc is a circle in $\d N_t$. Each pair of discs intersects in a single arc, and the three arcs are disjoint, so that each disc contains two disjoint double arcs.

Before the triple point birth, each disc intersects the other discs in a pair of disjoint double arcs. At the singular time the three double arcs intersect in a common point. After the singular time each double arc contains two triple points.

This can be modelled by choosing the red surface to be $z = 0$, the blue surface to be $x = 0$, and the green surface to be $z = y^2 - x +t$, where $t \in [-1, 1]$. The singular time is at $t = 0$. A triple point death is the time reverse of this.

\end{description}

\begin{definition}{\sl Football}
\index{football}

A {\sl football} is a region homeomorphic to a $3$-ball in $S^3$, whose boundary consists of one double curve free bigons of each colour. 

The images of a football under $G$ are disjoint, this is an immediate consequence of Corollary \ref{corollary:bigon orbit disjoint boundaries}.
A triple point birth move creates a football orbit, and a triple point death move destroys a football orbit.

\end{definition}

\eject

\begin{definition}{\sl Isotopy layer}
\index{isotopy layer}\index{isotopy!product}
\index{isotopy interval}\index{interval!isotopy}

Let $(M, \phi)$ be a sweepout, let $I$ be a closed subinterval of $\R$, and let $(h^{-1}(I), \phi)$ be the corresponding partial sweepout.

We say that $(h^{-1}(I), \phi)$ is an {\sl isotopy layer}, if there are no moves during the time interval $I$. We say that $I$ is an {\sl isotopy interval}.
\end{definition}

\begin{definition}{\sl Move layer}
\index{move layer}\index{layer!move}

Let $(M, \phi)$ be a sweepout, and let $I$ be a closed subinterval of $\R$. Suppose the partial sweepout $(h^{-1}(I), \phi)$ is a union of two partial sweepouts with disjoint interiors, one of which is a move, and the other is a partial sweepout which is a product. Then we say that $(h^{-1}(I), \phi)$ is a {\sl move layer}.
\end{definition}

\begin{definition}{\sl Move position}
\index{move position}

Let $(M, \phi)$ be a sweepout in general position, together with a choice of move neighbourhood for each move, so that the sweepout is a sequence of isotopy layers and move layers. Then we say that the sweepout is in {\sl move position}.
\end{definition}

\begin{lemma}
Given a sweepout in general position, we can choose move neighbourhoods for every move, so that the sweepout is in move position.
\end{lemma}

\begin{proof}
The critical points of the sweepout occur at distinct times, so we can choose move neighbourhoods which are disjoint in time.
\end{proof}

We will often wish to choose tubular neighbourhoods for subsets of the sweepout surfaces that intersect the sweepout spheres in a way which is ``as simple as possible''. We now make this precise.

\begin{definition}{\sl Thin regular neighbourhood}
\index{thin regular neighbourhood} \index{regular neighbourhood!thin}

Suppose that $\S$ is a subset of the sweepout spheres at some time $t$. 

Let $N \subset S^3$ be a regular neighbourhood of $\S$. If the intersection of $N$ with each of the sweepout spheres forms a regular neighbourhood in the sweepout spheres for $\S$, then we say that $N$ is a thin regular neighbourhood for $\S$.
\end{definition}

\begin{remark}
\index{thin tubular neighbourhood} \index{tubular neighbourhood!thin}
If $\S$ is a manifold, then we may choose a thin regular neighbourhood which is a tubular neighbourhood of $\S$. We call this a {\sl thin tubular neighbourhood} for $\S$.
\end{remark}

\begin{lemma} \label{lemma:isotopy}
Let $(M_I, f \cross h)$ be a partial sweepout with image neighbourhood $N_I$, which contains no moves, and let $K_0$ be an equivariant submanifold of $N_0$. Then we can extend $K_0$ to an equivariant continuously varying family of submanifolds $K_I$ contained in $N_I$, such that the intersection of the sweepout surfaces with $K_t$ only changes up to isotopy for $t \in I$. 
\end{lemma}

The idea is to think of the partial sweepout $(M, f \cross h)$ as an isotopy between the images of $M_0$ and $M_1$, and then use the isotopy extension theorem.

\begin{theorem}
Let $B$ be a manifold, which may have boundary, and let $A$ be a properly embedded compact submanifold of $B$. Let $f\co A \cross I \to B$ be an isotopy for which each $f_t$ is a proper embedding. Then $f$ extends to an ambient isotopy $F\co B \cross I \to B$.
\end{theorem}

We now prove Lemma \ref{lemma:isotopy}.

\begin{proof}
The image neighbourhood $N_I$ has a compatible product structure, so we can identify $N_I$ with $N_0 \cross I$, and think of $f\co M_I \to N_0$ as an isotopy between $f(M_0) = S_0$ and $f(M_1) = S_1$ which are properly embedded. We first show how to extend this isotopy to an equivariant ambient isotopy of $N_0$.

By the isotopy extension theorem, the isotopy $f$ extends to an isotopy $F'\co N_0 \cross I \to N_0$, which need not be equivariant. We can choose a continuously varying family of tubular neighbourhoods $U_t$ of $S_t$. As there are no moves, we can choose this neighbourhood to be sufficiently small, so that $\orbit U_t$ is a tubular neighbourhood for $\orbit S_t$.

Then $\orbit U_t$ is a $3$-manifold on which $G$ acts freely, so $\orbit U_t / G$ is also a $3$-manifold. Furthermore $\orbit U_I / G$ is an isotopy between $\orbit U_0 /G$ and $\orbit U_1 /G$ in $N_0 / G$, so by the isotopy extension theorem this extends to an isotopy $F''$ on $N_0 / G$. Let $F$ be the lift of $F''$ to $N_0$. The $F$ is an equivariant ambient isotopy that extends $f$.

Now let $K_I$ be the image of $K_0$ under the equivariant ambient isotopy $F$. The intersection of $K_t$ with the sweepout surfaces only changes by isotopy, by construction.
\end{proof}

Finally we show that every component of $M$ is either $S^3$ or a connect sum of $S^2 \cross S^1$'s.

\begin{proposition} \label{proposition:connect_sum}

Let $M$ be a closed $3$-manifold, and let $h\co M \to S^3$ be a Morse function, which has the property that at non-singular times the level sets are unions of $2$-spheres.

Then each component of $M$ is either $S^3$ or a connect sum of $S^{1} \cross S^{2}$'s.
\end{proposition}

\begin{proof}
Choose times $t_i$ between each of the critical times of the height function. The pre-images of these times is a collection of embedded $2$-spheres in $M$. They divide $M$ into pieces which contain at most one singular point of a singular level set. So $M - \S$ consists of pieces of the following three types:

\begin{enumerate}

\item If there is no critical point of $M$ in a component of $M - \S$, then the component is a  product $S^2 \cross I$.

\item If a component of $M - \S$ contains a critical point of index $0$ or $3$, in which a sphere either appears or disappears, then it is a $3$-ball.

\item If a component of $M - \S$ contains a critical point of index $1$ or $2$, in which a single sphere splits into two, or two spheres join together, then it is a $3$-ball with two $3$-balls removed.

\end{enumerate}

\begin{figure}[ht!]
\begin{center}
\epsfig{file=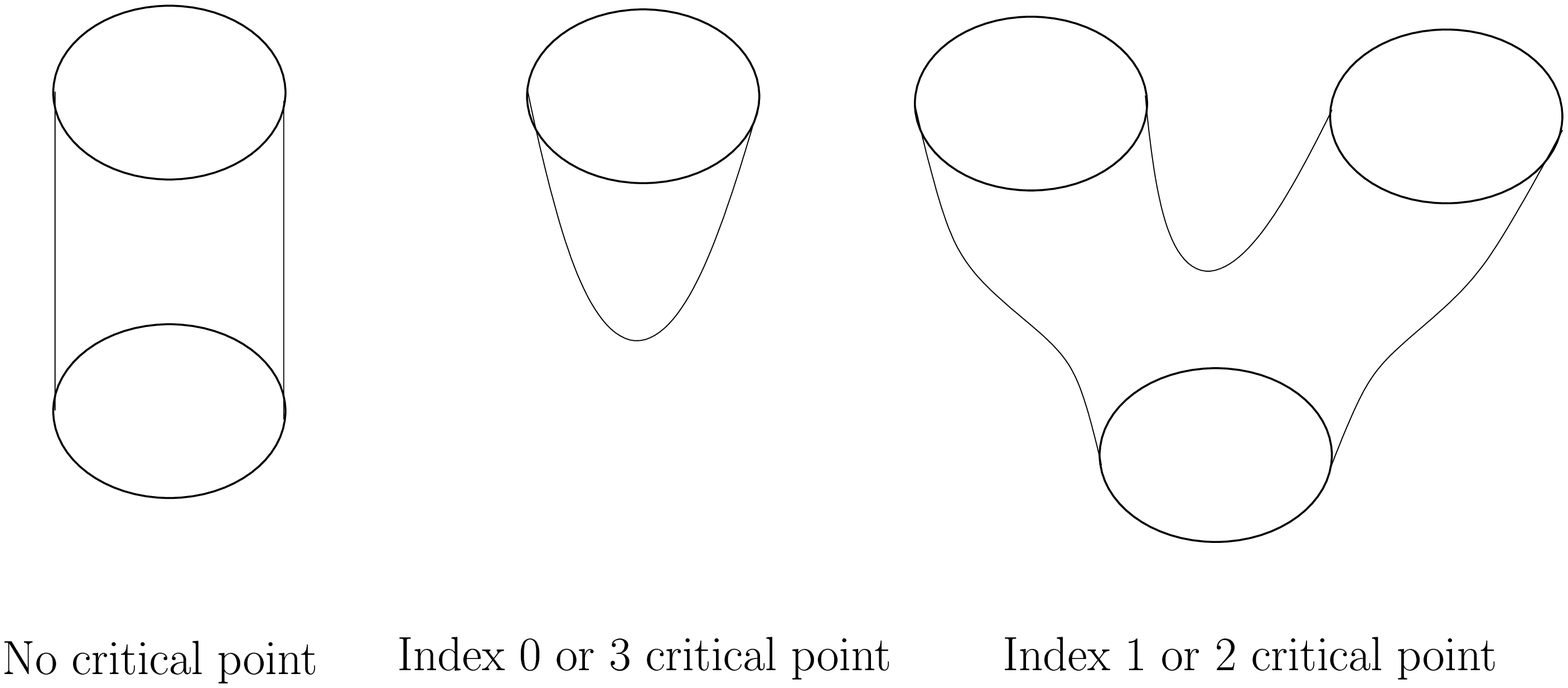, height=110pt}
\end{center}
\caption{Components of $M - \S$}
\end{figure}

In each case, capping off the $2$-sphere boundaries of the piece with $3$-balls produces a copy of $S^3$, so $M$ is either $S^3$, or a connect sum of $S^2 \cross S^1$'s.
\end{proof}

\subsection{Complexity}

\begin{definition}{\sl Complexity}
\index{complexity} \index{sweepout!complexity}

We say that the {\sl complexity} of the sweepout at the generic time $t$ is the ordered pair $(n,d)$, where $n$ is the number of triple points in the sweepout at time $t$, and $d$ is the number of simple closed double curves with no triple points at time $t$. We order the pairs lexicographically.
\end{definition}

Triple points come in multiples of six, double curves come in multiples of three.

\begin{definition}{\sl The graphic}
\index{graphic}

The graphic is the graph of complexity against time, with labels added to show which move occurs. We will mark every move, even if it does not change complexity, so an unmarked interval in the graphic corresponds to a time interval in which no moves occur.
\end{definition}

\begin{figure}[ht!] 
\centering
	\mbox{\subfigure[Key to symbols in the graphic]{%
			\begin{minipage}[b]{0.42\textwidth}\small
				\centering
				\begin{tabular}{cc}
					Move & Symbol \\
					Appear & A \\
					Vanish & V \\
					Cut & C \\
					Paste & P \\
					Birth & B \\
					Death & D \\
					Saddle & S \\
					Triple point birth & T+ \\
					Triple point death & T- \\
				\end{tabular}
				\vspace{0pt}%
			\end{minipage}%
		}
		\subfigure[The graphic]{%
			\begin{minipage}[b]{0.55\textwidth}
				\centering
				\epsfig{file=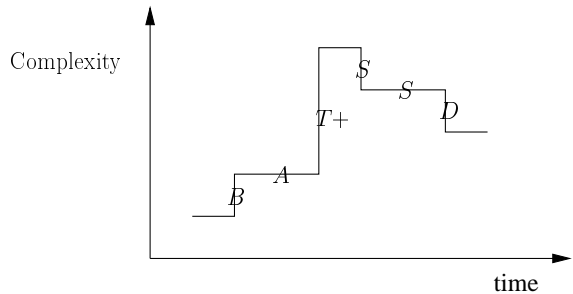, height=110pt}
				\vspace{0pt}%
			\end{minipage}%
		}
	}
\caption{The graphic}
\end{figure}

\subsection{Orientations}

We have chosen a continuously varying normal vector to $S_{t}$, which defines an inside and outside of the red spheres. We think of the inside as coloured red. This also gives us normal vectors for the images $gS_t$ and $g^2 S_t$, and we have coloured their insides green and blue respectively.

Each triple point lies at the intersection of three surfaces, one of each colour, so the triple point lies in the boundary of regions shaded with all possible combinations of colours. In particular, each triple point lies in the boundary of a ``clear'' region, ie, a region on the ``outside'' of all of the spheres.

Choose tangent vectors to the double curves at the triple point, so that the tangent vectors point toward the clear region. The vectors have a cyclic ordering (an orientation) coming from the map $g$ ie (red, green, blue).

\begin{definition}{\sl Positive and negative triple points}
\index{triple point!orientation} \index{orientation!triple points}

If this orientation agrees with the orientation of $S^3$, we call the triple point {\sl positive}, if it disagrees we call it {\sl negative}. 
\end{definition}

Every triple point is therefore labelled either positive or negative, and triple points that are connected by double arcs containing no triple points in their interiors have opposite sign.

\begin{figure}[ht!]
\begin{center} 
\epsfig{file=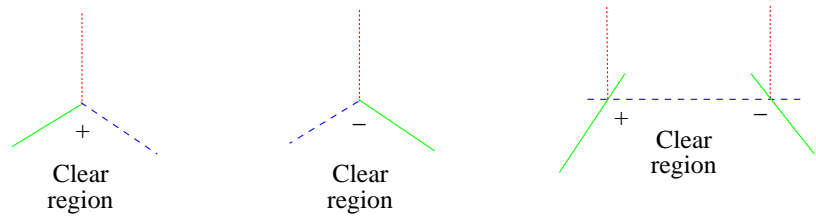, height=81pt}
\end{center}
\caption{Signs of triple points}
\end{figure}

\begin{definition}{\sl Adjacent}
\index{adjacent}

We say that two triple points are {\sl adjacent} if they are connected by a double arc which contains no triple points in its interior.
\end{definition}

Triple points are adjacent if they are connected by a red arc, so triple points may still be adjacent, even if they are not adjacent in the blue-green graph.

\begin{lemma}{} \label{neighbour}
Adjacent triple points cannot lie in the same orbit under $G$.
\end{lemma} 

\begin{proof}
The map $f$ preserves the sign of the triple point. Adjacent triple points have opposite sign, so they can not be images of each other.
\end{proof}

\begin{corollary} \label{corollary:bigon orbit disjoint boundaries}
A bigon orbit consists of a red bigon, a blue bigon and a green bigon, which have disjoint boundaries.
\end{corollary}

\begin{proof}
The vertices of a red bigon are adjacent triple points, so must have distinct orbits under $G$. Therefore the images of the arcs in the boundary of the red bigon have distinct endpoints, so they must be distinct also.
\end{proof}

\subsection{Modifying sweepouts} \label{subsection:modifications}

\begin{definition}{\sl A modification neighbourhood}
\index{modification neighbourhood}

Let $(M, \phi)$ be a sweepout in move position. Let $N_I$ be a compatible product neighbourhood with constant boundary pattern. If $\d (N_I)$ is disjoint from the chosen move neighbourhoods of the sweepout, then we say $N_I$ is a {\sl modification neighbourhood}.
\end{definition}

\begin{notation}
If $P \subset M$, and $\phi$ is a map defined on $M$, then we will abuse notation and write $\phi$ to denote the restriction map $\phi |_P$.
\end{notation}

\begin{definition}{\sl A pre-image sweepout}
\index{sweepout!pre-image} \index{pre-image sweepout}

Let $(M, \phi)$ be a sweepout in move position, and let $N_I$ be a modification neighbourhood. Then $(\phi^{-1}(N_I), \phi)$ is a partial sweepout which we call the {\sl pre-image partial sweepout for $N_I$}.
\end{definition}

Every modification neighbourhood is the image neighbourhood of its pre-image sweepout.

\begin{remark}
We can construct a partial sweepout by drawing a sequence of configurations in $N$ which differ by moves.
\end{remark}

\begin{figure}[ht!]
\begin{center}
\epsfig{file=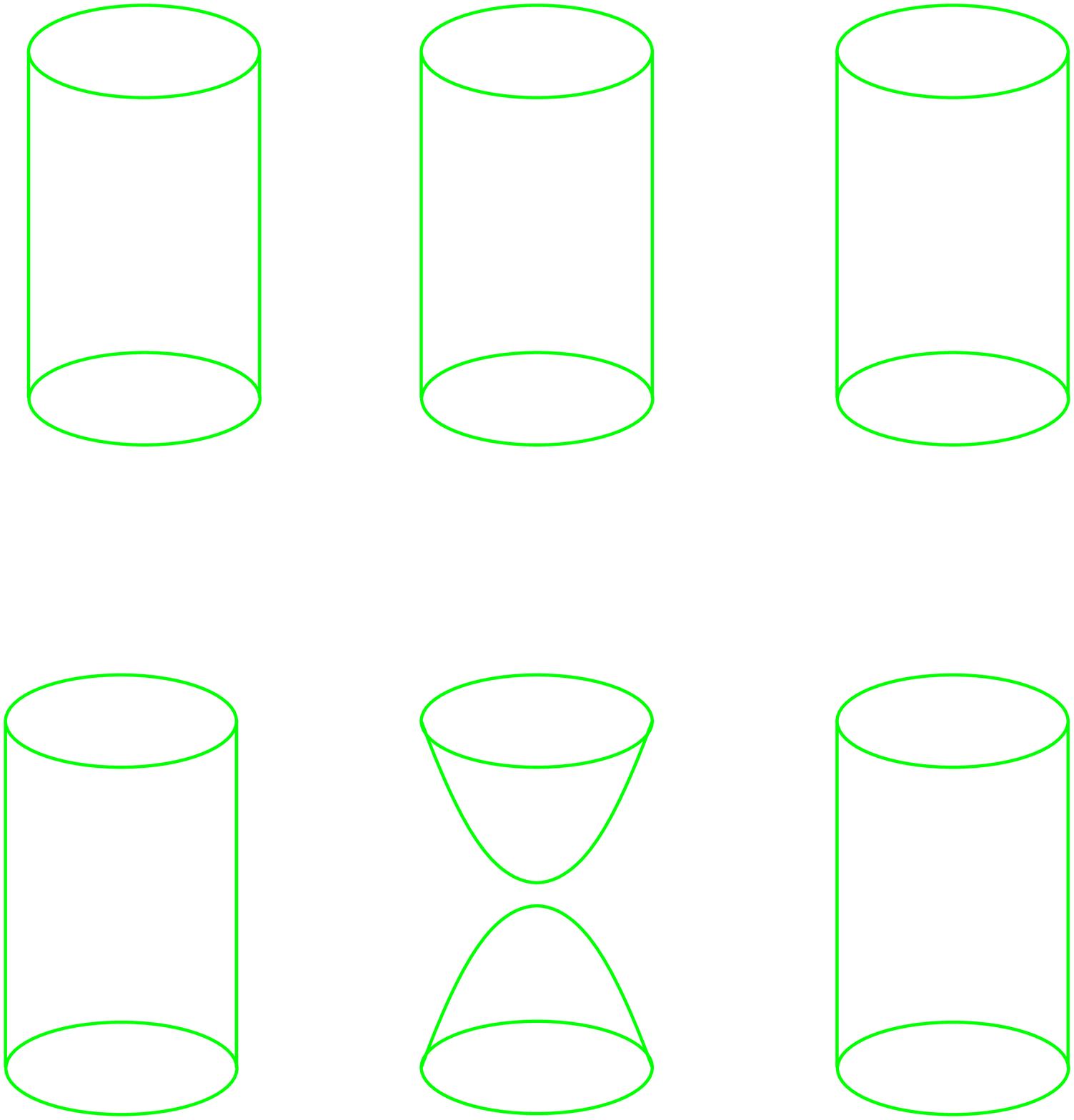, height=200pt}
\end{center}
\caption{Two partial sweepouts with the same boundary}
\end{figure}

\begin{definition}{\sl Partial sweepouts with the same boundary}
\index{partial sweepout!same boundary}

Let $(P, \theta)$ and $(P', \theta')$ be two partial sweepouts with the same image neighbourhoods. If $\theta^{-1} \circ \theta'$ is a diffeomorphism between $\d P$ and $\d P'$ then we say that the two partial sweepouts {\sl have the same boundary}.
\end{definition}

\begin{theorem} \label{theorem:mod}
Let $(M_1, \phi_1)$ be a sweepout in move position, let $N_I$ be a  modification neighbourhood with pre-image sweepout $(P_1, \phi_1)$, and let $(P_2, \phi_2)$ be a partial sweepout with the same boundary as the pre-image sweepout. Assume that $N_0$ is not all of $S^3$. Define $M'$ to be the manifold $(M_1 \complement P_1) \cup P_2$, using the gluing map $\phi_1^{-1} \circ \phi_2$, and define $\phi'$ as follows:

$$ \phi'  = g' \cross h' = \begin{cases} \phi_2\text{ on }P_2 \\ \phi_1\text{ on }M_1 - P_1 \end{cases}$$

If all but finitely many level sets ${h'}^{-1}(t)$ are a union of $2$-spheres, then $(M', \phi')$ is a sweepout.
\end{theorem}

\begin{proof} 
It remains to show that $f' =  \pi_{S^3} \circ \phi'$ is degree one. Let $f = \pi_{S^3} \circ \phi_1$. As $N_t$ is not all of $S^3$ we may assume (perhaps after adjusting by an isotopy) that there is a point $p \in S^3$ which is disjoint from $N_t$ for all $t \in I$. Then $f^{-1}(p)$ is the same as ${f'}^{-1}(p)$, so both maps have the same degree.
\end{proof}

\begin{definition}{\sl Modification}
\index{modification}

Using the notation from Theorem \ref{theorem:mod} above, we say that $(M', \phi')$ is a {\sl modification} of the sweepout $(M, \phi)$.
\end{definition}

The following two modifications will be useful:

\begin{description}

\item[Removing a simple closed double curve which bounds a disc]$\phantom S$ 
\index{modification!removing double curves}

Let $(M, \phi)$ be a sweepout in move position. Suppose for some time interval $I$ we can choose a continuous family of double curves $\{ \g_t \mid t \in I \}$, which bound a continuous family of discs $\{ \D_t \mid t \in I \}$, which contain no double curves in their interiors, and which are disjoint from all the move neighbourhoods in the sweepout. By Lemma \ref{lemma:isotopy}, we can choose a compatible product neighbourhood $N_I$ so that each $N_t$ is a thin tubular neighbourhood for the orbit of the disc $\D_t$, and $N_I$ is disjoint from the move neighbourhoods. This means that the intersection of each $N_t$ with the sweepout surfaces consists of the orbit of a disc and an annulus, which intersect transversely in the orbit of the simple closed double curve $\g_t$. Let $(P, \phi)$ be the pre-image sweepout for $N_I$. This is illustrated by the top row of pictures in Figure \ref{figure:removing a double curve}.

We can remove the double curve by doing a cut move parallel to $\D_t$, and then a double curve death move. If we follow this pair of moves by its time reverse, a birth followed by a paste, we return the configuration to its initial state.  This creates a new partial sweepout with the same boundary as $(P, \phi)$, but which does not contain any double curves for a subinterval of $I$. This is illustrated by the bottom row of pictures in Figure \ref{figure:removing a double curve}.

\begin{figure}[ht!]
\begin{center}
\epsfig{file=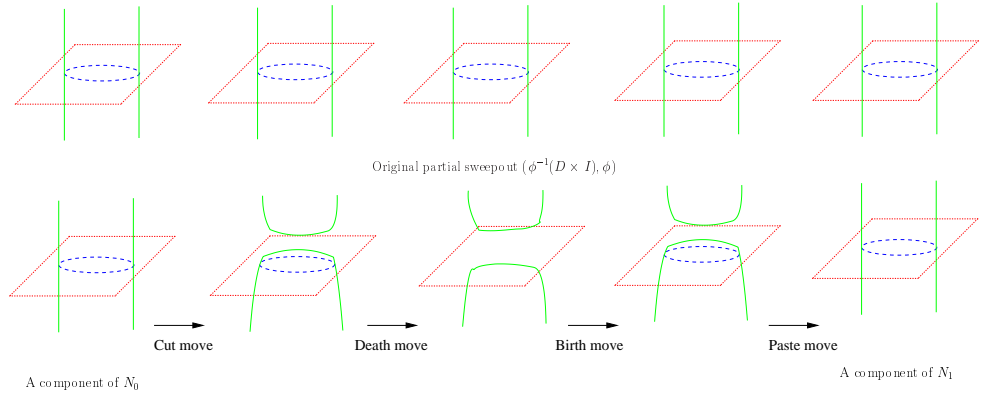, width=.98\hsize}
\end{center}
\caption{Removing a double curve which bounds a double curve free disc} \label{figure:removing a double curve}
\end{figure}

We will call the pair of moves which removes the double curve a {\sl compound double curve death}, and we call the pair of moves which replaces the double curve a {\sl compound double curve birth}. No other moves may occur during a compound double curve move, but moves may occur in between the compound double curve birth and the compound double curve death.

\begin{figure}[ht!]
\begin{center}
\epsfig{file=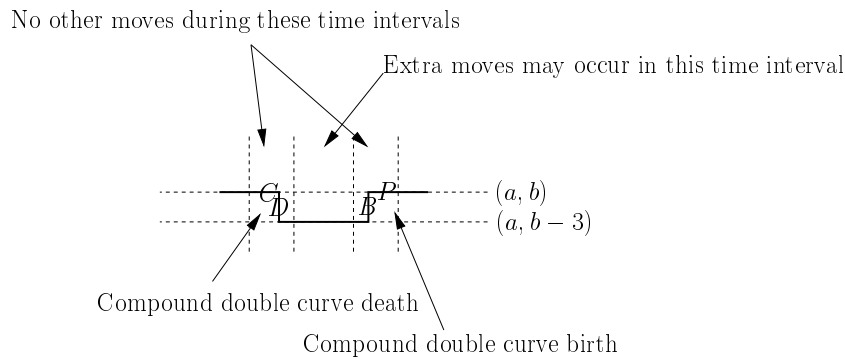, width=.84\hsize}
\end{center}
\caption{Removing a double curve changes the graphic}
\end{figure}

\item[Removing a double curve free bigon] 
\index{modification!removing bigons}

Let $(M, \phi)$ be a sweepout in move position. Suppose there is a time interval $I$ for which we can choose a continuous family of bigons $\{ A_t \mid t \in I \}$, so that $A_t$ is always double curve free, and disjoint from the move neighbourhoods of the sweepout. 
By Corollary \ref{corollary:bigon orbit disjoint boundaries} such a bigon is disjoint from its images under $G$.
By Lemma \ref{lemma:isotopy}, we can choose a compatible product neighbourhood $N_I$ disjoint from the move neighbourhoods, so that each $N_t$ is a thin tubular neighbourhood for the bigon orbit $\a_t$. As the three images of the bigon $A_t$ are disjoint, $N_t$ will be the orbit of a $3$-ball, which is also disjoint from its images. So it will suffice to describe how to change the sweepout in a single component of $N_I$, as the changes can be extended equivariantly to the other components.

Consider the component of $N_t$ which contains the bigon $A_t$.
This intersects the sweepout surfaces in three discs, one of each colour. One of these discs, which we shall think of as being ``horizontal'' is a tubular neighbourhood of the bigon $A_t$ in the sweepout sphere which contains $A_t$. Each of the other discs, which we shall think of as ``vertical'', intersect the horizontal disc transversely in a single double arc. These double arcs contain the boundary of the bigon $A_t$. The vertical discs intersect each other in a pair of vertical double arcs. This is illustrated in the top row of pictures in Figure \ref{figure:removing a bigon}.

We can remove the bigon by pushing the vertical surfaces through each other sideways. In order to do this in general position, we first do a saddle move, using a saddle disc parallel to $A_t$. This saddles the vertical double arcs and creates a football containing $A_t$ in its boundary. We can then do a triple point death move to destroy the football, removing the bigon, and reducing the number of triple points in this component by two, and by six in total. If we then do the time reverse of these two moves, we have created a partial sweepout with the same boundary as the original one, but which has fewer triple points for a subinterval of $I$. This is illustrated in the bottom row of pictures in Figure \ref{figure:removing a bigon}.

\begin{figure}[ht!]
\begin{center}
\epsfig{file=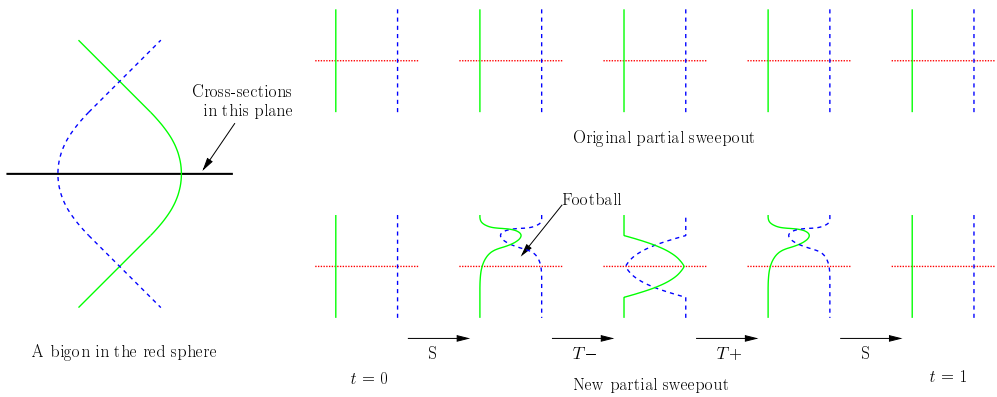, width=\hsize}
\end{center}
\caption{Removing a double curve free bigon}\label{figure:removing a bigon}
\end{figure}

We will call the pair of moves that removes the bigon a {\sl compound triple point death}, and the pair of moves that replaces the bigon a {\sl compound triple point birth}. No other moves may occur during the pair of moves that make up a compound move. We will call triple point moves that are not part of compound triple point moves {\sl non-compound} triple point moves.

\begin{figure}[ht!]
\begin{center}
\epsfig{file=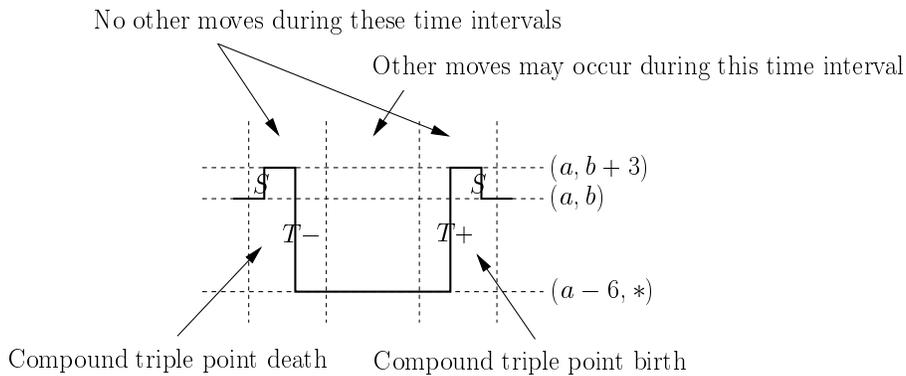, height=140pt}
\end{center}
\caption{Removing a double curve free bigon changes the graphic}
\end{figure}

The saddle move may leave complexity unchanged, or may increase it by $(0,3)$. We have illustrated the latter case above.

\end{description}

We can do one of these modifications whenever we can find a simple closed curve that bounds a disc, or a double curve free bigon, that is disjoint from the move neighbourhoods of the sweepout for some time interval. If this time interval contains a local maximum in the graphic, the effect of the modification will be to reduce the height of the local maximum. We say that we have {\sl undermined} the local maximum.

In the discussion above, we always wrote $\a_t$ to refer to a continuously varying family of subsets of the sweepout spheres defined for some time interval. In future, we will just write $\a$ without the subscript to refer to such a family, if it is clear from context that we mean a family defined for a time interval.

\section{Reduction to special cases} \label{section:undermining}

In this section we will prove the main theorem, assuming the following three lemmas, which we will prove in the remaining sections.

\begin{lemma1}
Every sweepout contains triple points.
\end{lemma1}

\begin{lemma2}
Suppose a non-triple point move occurs while the configuration contains triple points. Then there is a vertex free bigon orbit disjoint from the move neighbourhood of the move.
\end{lemma2}

\begin{definition}{\sl A special case local maximum}
\index{special case local maximum} \index{local maximum!special case}

We say that a local maximum in the graphic is a {\sl special case local maximum} if it consists of a compound triple point birth, followed by a compound triple point death, with no other moves occurring in between.
\end{definition}

\begin{lemma3}
A special case local maximum can either be undermined, or else contains an invariant unknotted circle.
\end{lemma3}

Lemma \ref{lemma2} enables us to use the two modifications described in the previous section to change all the local maxima into special case local maxima, so that we can apply Lemma \ref{lemma3}. If we could undermine all of the local maxima, then we would have a sweepout with no triple points, contradicting Lemma \ref{lemma1}, so there must be a special case local maximum we cannot undermine, which therefore contains an unknotted invariant curve.

So in order to prove the main theorem, it suffices to show that given a sweepout of complexity $\geqslant (6,0)$, we can change it to produce a new sweepout with lower complexity. Therefore the following lemma proves the main theorem, assuming the three lemmas above.

\begin{lemma} \label{lemma:reduction}
Suppose $(M, \phi)$ is a sweepout with maximum complexity $(p,q) \geqslant (6,0)$. Then either we can find a new sweepout with lower maximum complexity, or else there is an invariant unknotted curve.
\end{lemma}

First we deal with the case where the complexity of the sweepout is $(p,q)$, with $q > 0$, and show that if there are no invariant unknotted curves we can reduce the complexity to $(p, q-3)$. We divide the proof into three steps, which we describe below. Finally, we show how to deal with the case $q=0$.

\setcounter{step}{0}
\begin{step}{Undermine non-triple point moves}

Lemma \ref{lemma2} implies that we can repeatedly undermine every non-triple point move, until they all lie below $(p,0)$. This produces a sweepout with complexity at most $(p, q+3)$, but now all moves at complexities at least $(p,0)$ are either compound moves, or triple point births and deaths. 
\end{step}

\begin{step}{Replace all non-compound triple point moves with compound triple point moves}

We can replace a non-compound triple point death with the following sequence of moves. Choose one of the bigons of the football, call it $A$, and do a compound triple point death move using it. We can choose to do the saddle move of the compound triple point death move using a saddle disc which lies inside the football. This splits the football into a smaller football close to $A$, and a $3$-ball region bounded by a pair of discs which meet in a simple closed double curve. The triple point death move of the compound move removes the football close to $A$, and we can now do a double curve death move to remove the simple closed double curve. This is illustrated below.

\begin{figure}[ht!]
\begin{center}
\epsfig{file=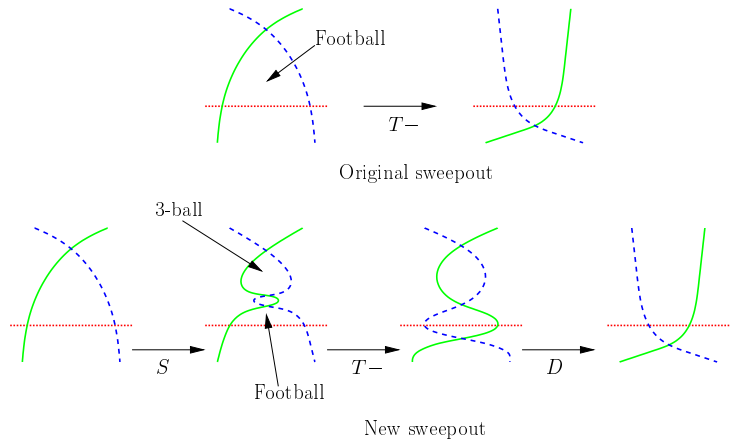, height=165pt}
\end{center}
\caption{Replacing a non-compound triple point move with a compound triple point move}
\end{figure}

This does not increase complexity above $(p,q+3)$, as non-compound triple point moves have heights at most $(p,q)$. This introduces non-triple point moves at heights strictly less than $(p, 0)$, but now all moves above $(p,q-3)$ are compound moves.

\begin{figure}[ht!]
\begin{center}
\epsfig{file=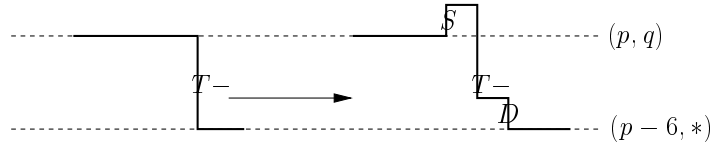, height=54pt}
\end{center}
\caption{Replacing a non-compound triple point move with a compound triple point move}
\end{figure}

\end{step}

\begin{step}{Reducing compound local maxima}

We now consider the local maxima of complexity $(p,q)$ or higher, and show we can change the sweepout to reduce them all to complexity at most $(p,q-3)$. First, in Case 1, we show how to undermine the compound local maxima of complexity $(p,q+3)$, so that they have complexity at most $(p,q)$. In Case 2, we then show how to undermine compound local maxima of height $(p,q)$, with $q>0$. Finally, in Case 3, we show how to undermine compound local maxima of height $(p,0)$.

\setcounter{case}{0}
\begin{case}{A local maximum of complexity $(p,q+3)$}

A compound move always changes complexity, so the local maximum consists of one compound move, immediately followed by another compound move. If a local maximum has complexity $(p, q+3)$, then it must contain at least one compound triple point move, as these are the only compound moves at this height. The other move may be either another compound triple point move, or a compound double curve move. If the local maximum contains a compound double curve move, then there are two cases, depending on whether the compound double curve move occurs first, followed by the triple point move, or whether the compound triple point move occurs first, followed by the compound double curve move. The picture below shows all three possibilities.

\begin{figure}[ht!]
\begin{center}
\epsfig{file=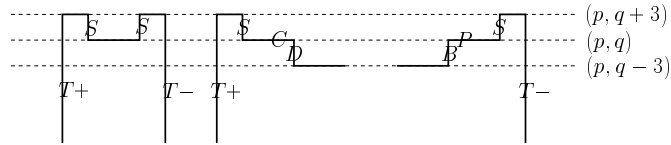, height=54pt}
\end{center}
\caption{Local maxima of complexity $(p, q+3)$}
\end{figure}

The first case shown is a special case local maximum, which by Lemma \ref{lemma3}, can either be undermined, or contains an invariant unknotted curve. It now suffices to show how to undermine the second case, as the third
 local maximum is just the time reverse of the second one.

Let $\a$ be the bigon orbit of the bigon in the compound triple point move, and let $\D$ be the disc orbit of the disc in the compound double curve move. The bigon orbit is double curve free, so $\a$ must be disjoint from the disc orbit $\D$, so we can just swap the order of the compound moves. 

\begin{figure}[ht!]
\begin{center}
\epsfig{file=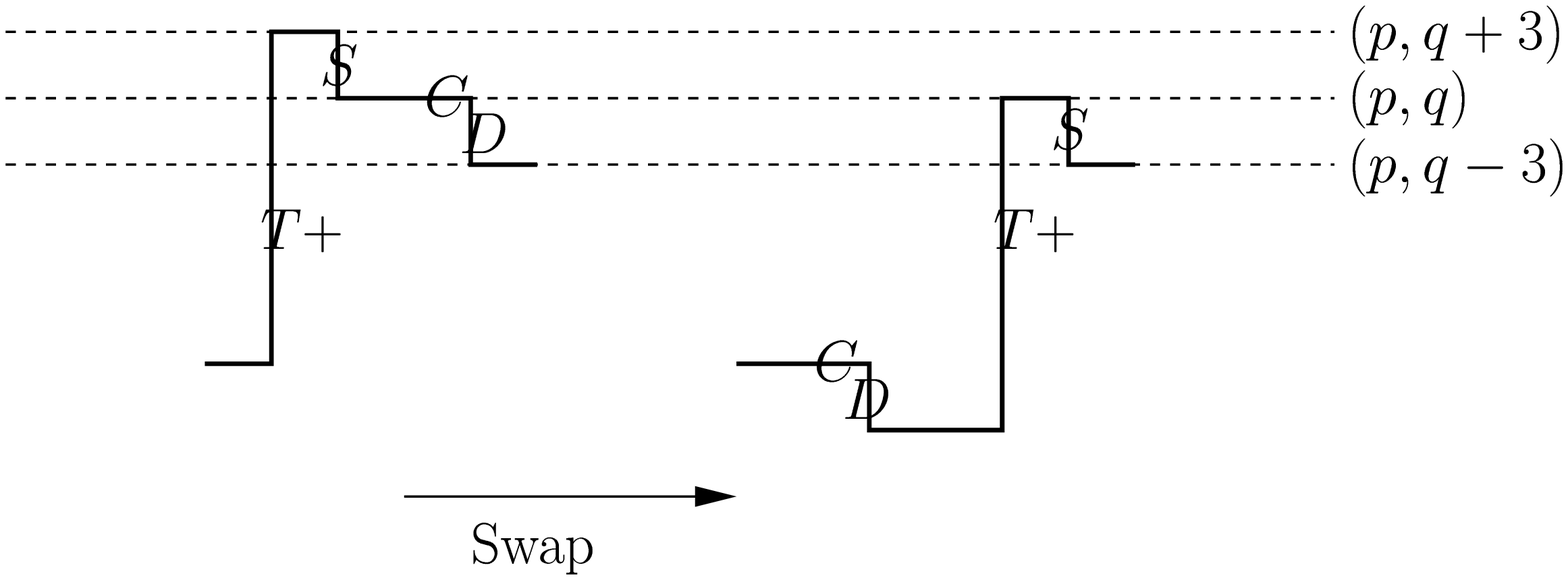, height=100pt}
\end{center}
\caption{A compound triple point birth followed by a compound double curve death}
\end{figure}

We may now assume that we have reduced the complexities of all the local maxima to at most $(p,q)$.
\end{case}

\begin{case}{A local maximum at height $(p,q)$, with $q>0$}

A compound move always changes complexity, so a local maximum of height $(p,q)$ consists of one compound move, followed by another compound move. Each compound move may be either a compound double curve move, or a compound triple point move. The picture below shows all possible local maxima of complexity $(p,q)$, up to time reverse. The time reverse of the second picture may also occur.

\begin{figure}[ht!]
\begin{center}
\epsfig{file=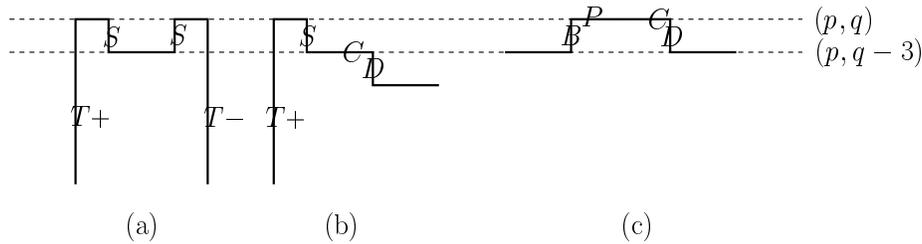, height=90pt}
\end{center}
\caption{Local maxima of complexity $(p,q)$} \label{picture573}
\end{figure}

By Lemma \ref{lemma3}, the first case, Figure \ref{picture573}(a), can either be undermined, or else the configuration contains an unknotted invariant curve. The second case, Figure \ref{picture573}(b), and its time reverse can be undermined by the argument given above in Case 1. For the remaining case, Figure \ref{picture573}(c), we show that it can either be undermined, or reduced to one of the first two cases.

The remaining case consists of a compound double curve birth followed by a compound double curve death. Let $\D_1$ be the orbit of the double curve free disc in the first compound double curve move, and let $\D_2$ be the orbit of the double curve free disc in the second compound double curve move.

\begin{figure}[ht!]
\begin{center}
\epsfig{file=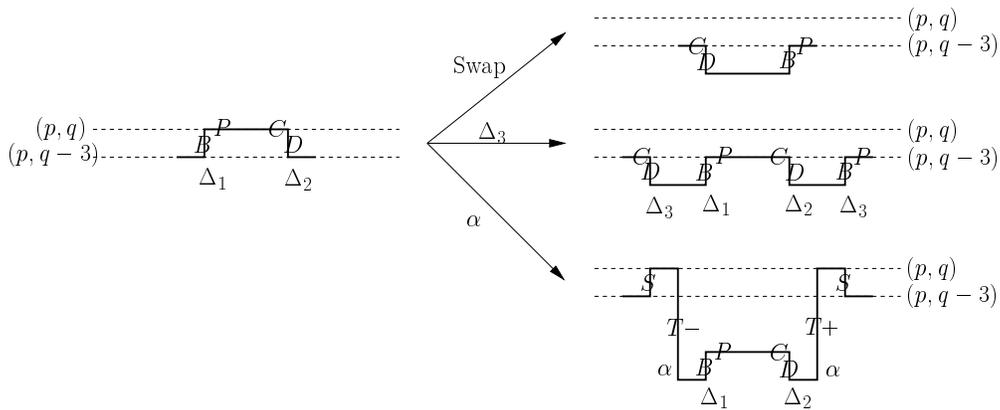, width=\hsize}
\end{center}
\caption{A compound double curve birth followed by a compound double curve death}
\end{figure}

If $\D_1$ and $\D_2$ are disjoint, then we can just swap the order of the compound moves. This undermines the local maximum.

If $\D_1$ and $\D_2$ are not disjoint, then they must share a common boundary, which is the orbit of a simple closed curve $\g$. The green and blue components of $\g$ are contained in the red spheres, while the red component is disjoint from the red spheres. As there are triple points, by Lemma \ref{4lemma} there are at least four vertex free bigons with disjoint interiors, so there is a vertex free bigon orbit $\a$ disjoint from the orbits of $\D_1$ and $\D_2$. So either $\a$ is double curve free, or else it contains an innermost simple closed curve, which bounds a double curve free disc $\D_3$.

If there is a simple closed  double curve which bounds a double curve free disc $\D_3$, which is disjoint from $\D_1 \cup \D_2$, then we can undermine the local maximum by removing the double curve for the duration of the local maximum, by applying a double curve removal modification.

If there is a double curve free bigon $\a$, which is disjoint from $\D_1 \cup \D_2$, then we can modify the sweepout by inserting a compound triple point death, followed by a compound triple point birth, which removes the bigon for the duration of the local maximum. This reduces the complexity of the compound double curve local maximum, but may not reduce the overall complexity of the sweepout, as the extra compound triple point moves may still have maximum complexity $(p,q)$. However, we can remove all local maxima consisting only of compound double curve moves in this way, and then the remaining local maxima of height $(p,q)$ will have at least one compound triple point move. This means that they are all covered by the previous two cases, and so can be undermined.
\end{case}

\begin{case}{A local maximum of height $(p,0)$}

Finally, we explain how to deal with the case when the sweepout has complexity $(p,0)$. In this case, the local maxima of greatest complexity do not contain any simple closed double curves without triple points. We can again use Lemma \ref{lemma1} to undermine the non-triple point moves of greatest complexity, but as there are no simple closed double curves without triple points (because $q=0$), all the compound moves will be compound triple point moves. This means that the resulting sweepout may have local maxima of complexity $(p,0)$ and $(p,3)$, but all the local maxima at these heights will consist either of compound triple point moves or of non-compound triple point moves.

We can now apply Step 2 to replace all the non-compound triple point moves with compound triple point moves, and so now all local maxima are special case local maxima, and we can apply Lemma \ref{lemma3} directly.\end{case}

This completes Step 3, the final step in the proof of Lemma \ref{lemma:reduction}.
\end{step}

We have now proved the main result, assuming Lemmas \ref{lemma1}, \ref{lemma2} and \ref{lemma3} above. The remaining sections are devoted to proving these three lemmas.

\section{Every sweepout contains triple points} \label{section:triple points}

In this section we prove the following lemma:

\begin{lemma} \label{lemma1}
Every sweepout contains triple points.
\end{lemma}

At the beginning of the sweepout there is a clear connected invariant region, which we will call an initial region. We show that if a sweepout has no triple points, then there is always an initial region. However, this gives a contradiction, as at the end of the sweepout there are no sweepout spheres, and $S^3$ is a single region coloured with all three colours, so there are no clear regions at all.

\begin{definition}{\sl An initial region}
\index{initial region} \index{region!initial}

We say that a region is {\sl initial} if it is  clear, connected and invariant under $f$. 
\end{definition}

\begin{remark}
The union of the clear regions gets mapped to itself by $f$, but need not be connected. Individual connected clear regions need not be invariant.
\end{remark}

\begin{lemma} \label{lemma:breakup}
A sweepout without triple points contains initial regions at all times.
\end{lemma}

\begin{proof}
At the beginning of the sweepout, there is an initial region. This is because before the first appear move, there are no sweepout surfaces, and so all of $S^3$ is a clear connected invariant region, and so is initial.

We now complete the proof of the lemma by considering each type of move in turn, and showing that if there is an initial region before the move, then there must also be an initial region after the move, assuming that there are no triple points.

A single move is supported in the orbit of a $3$-ball which is disjoint from its images, so a single move cannot eliminate an initial region by shrinking it down to a point. However a single move might eliminate an initial region by disconnecting it.

\setcounter{case}{0}
\begin{case}
{Appear and vanish moves}

If a double curve free $2$-sphere orbit appears inside an initial region, then the interior of the $2$-sphere orbit is a new region, namely a $3$-ball orbit coloured by the same colour as the sphere. The region on the outside of the $2$-sphere orbit is still connected, and so is still an initial region.
\end{case}

\begin{case}
{Cut and paste moves}

A paste move does not disconnect any region, so we need only consider cut moves. If a cut move disconnects an initial region $C$, then the cut disc $D$ must be contained inside $C$. The images of $D$ must also lie inside $C$, so we can find a path $\gamma$ from $D$ to $gD$ which does not intersect any other sweepout surface. There are two cases, depending on whether or not the curve $\g \cup g\g$ intersects the cut disc transversely.

\begin{figure}[ht!]
\begin{center}
\epsfig{file=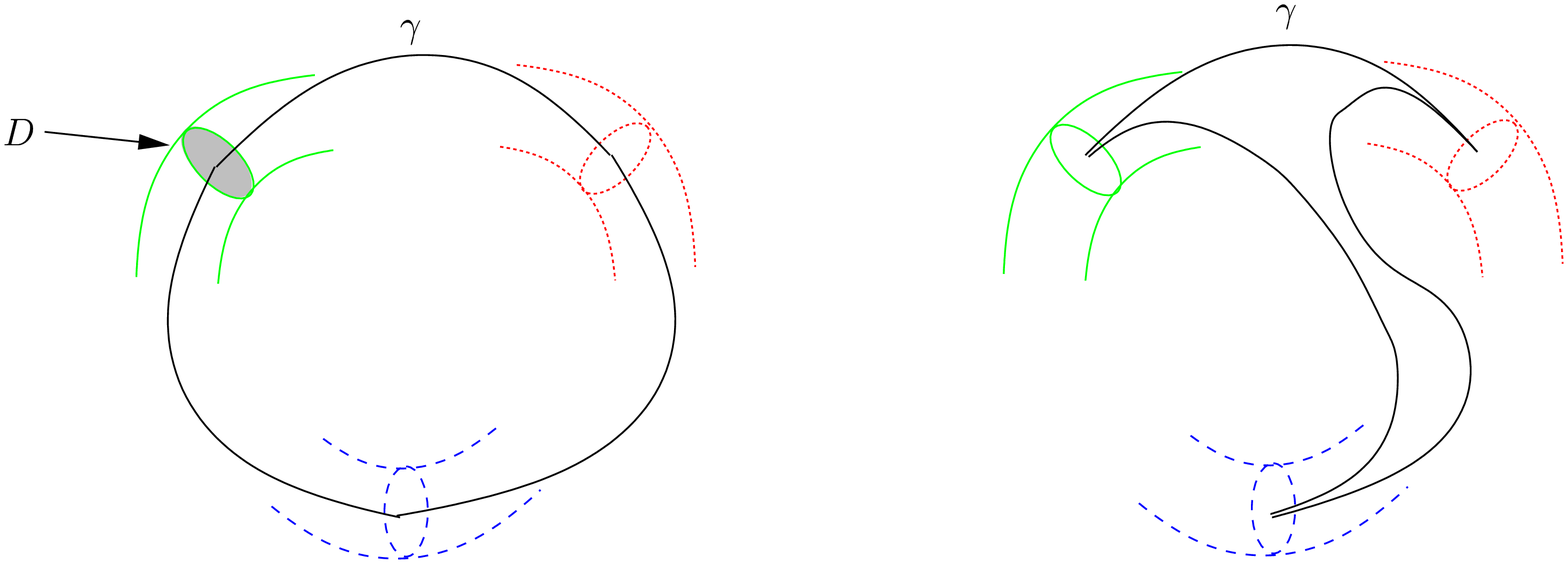, width=.95\hsize}
\end{center}
\caption{A cut move}
\end{figure}

In the right hand case the region containing $\orbit \gamma$ is still an initial region after the cut move.
In the left hand case, the disc $D$ divides the green sphere into two discs. The union of $D$ with one of these discs is a sphere. As the simple closed curve $\orbit \gamma$ lies in the clear region it is disjoint from the green sphere, so it intersects this $2$-sphere precisely once transversely, a contradiction, as $H_1(S^3) = 0$, so the left hand case cannot occur.
\end{case}

\begin{case}
{Double curve births and deaths}

\begin{figure}[ht!]
\begin{center}
\epsfig{file=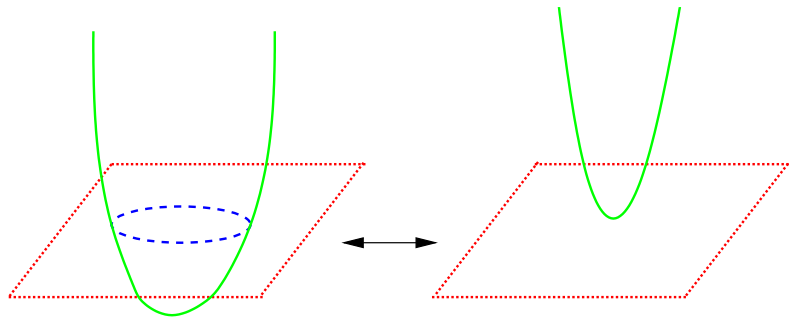, height=120pt}
\end{center}
\caption{Birth/death of a blue double curve}
\end{figure}

The region that is created or destroyed in the move cannot be invariant. The other regions that intersect the move neighbourhood are not disconnected by the move. 
\end{case}

\begin{case}
{Saddle moves}

\begin{figure}[ht!]
\begin{center}
\epsfig{file=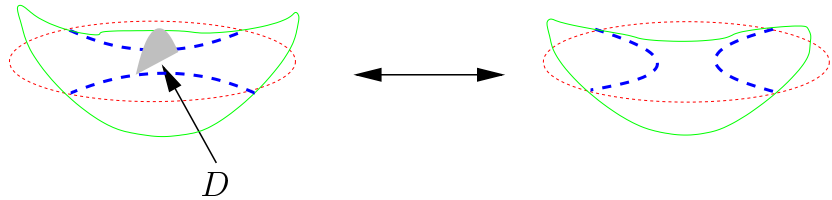, height=78pt}
\end{center}
\caption{A saddle move}
\end{figure}

The only region which could become disconnected as a result of a saddle move is the region containing the saddle disc $D$ in the diagram above. If the disc $D$ intersects two distinct blue circles, then as there are no triple points, we can choose a path from one side of $D$ to the other, in the interior of the initial region, parallel to either one of the blue double curves. So in this case, the saddle move does not disconnect the initial region.

Now suppose that the disc $D$ intersects a single blue double curve. The disc $D$ has three disjoint images under $G$. The region is connected, so we can find a path $\gamma$ from $D$ to $gD$, which only intersects the discs in its endpoints, and is disjoint from all the sweepout surfaces. Then $\orbit \gamma$ is a simple closed curve that connects the three discs in one of the following two ways:

\begin{figure}[ht!]
\begin{center}
\epsfig{file=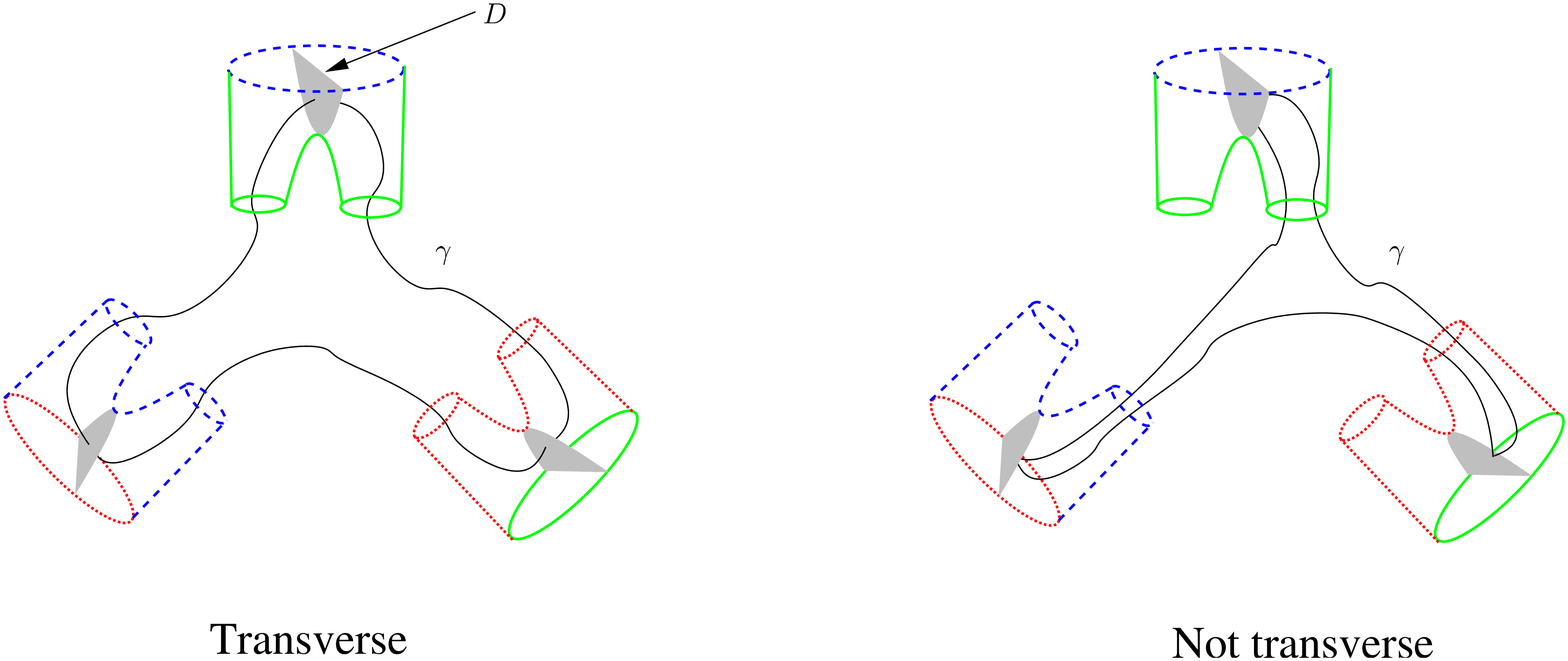, height=150pt}
\end{center}
\caption{The invariant curve $\orbit \g$}
\end{figure}

In the one-sided case, the region containing $\orbit \gamma$ is still an initial region after the saddle move. We now show that the transverse case cannot arise. 

The green double curve, which we will call $\alpha$, divides one of the red spheres into two discs. Consider the disc which contains part of the boundary of the saddle disc $D$. This disc is divided into two subdiscs by the saddle arc. As the orbit of the saddle disc $\orbit D$ intersects the red spheres in precisely two arcs, at most one of the subdiscs of the red disc intersects an image of the saddle disc $D$. Let $D'$ be one of the red subdiscs whose interior is disjoint from $\orbit D$. This is illustrated below. The red sphere is drawn as horizontal, and the blue sphere intersects it transversely in the double curve $\alpha$.

\begin{figure}[ht!]
\begin{center}
\epsfig{file=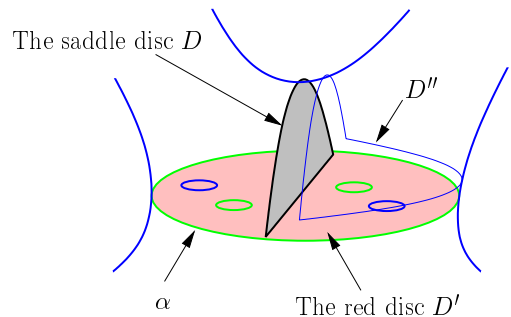, height=120pt}
\end{center}
\caption{The disc $D''$ is parallel to $D \cup D'$}
\end{figure}

The union of $D$ and $D'$ form a disc whose boundary lies in one of the blue spheres. Let $D''$ be a disc parallel to $D \cup D'$. We may choose this disc so that it is disjoint from the double curves, and from $\orbit D$, so it intersects the sweepout surfaces in simple closed curves only. Cut move do not disconnect the initial region, by Case 2, so we may remove the intersections of the sweepout surfaces with $D''$. We may need to change the path $\gamma$ if it passes through a cut disc, however we may always choose a new path from $D$ to $gD$, which we shall also call $\gamma$, as cut moves do not disconnect the initial region. Now the boundary of $D \cup D'$ bounds a blue disc with no double curves in its interior. However, the union of this blue disc with the red disc $D'$, and the saddle disc $D$, is a $2$-sphere, which the curve $\orbit \g$ intersects precisely once transversely, a contradiction. This shows that the transverse case cannot arise.
\end{case}

Therefore there is always an initial region after a move in a sweepout without triple points. This completes the proof of Lemma \ref{lemma:breakup}, which in turn completes the proof of Lemma \ref{lemma1}, that every sweepout must contain triple points.
\end{proof}

\section{Disjoint bigons for non-triple point moves} \label{section:disjoint}

In this section we prove the following lemma:

\begin{lemma} \label{lemma2}
Let $N_I$ be a move neighbourhood for a non-triple point move in a sweepout. Then there is a vertex free red bigon disjoint from $N$ during the time interval $I$.
\end{lemma}

This means that the orbit of this red bigon will also be disjoint from the move neighbourhood of the move. The configuration outside the move neighbourhood only changes by an isotopy during the move interval, so it suffices to find a bigon disjoint from the move neighbourhood at a given time during the move interval. We show that a configuration with triple points has at least four red bigons, and that a move neighbourhood can't intersect all of these at once.

\begin{definition}{\sl A red pseudo-bigon}
\index{pseudo-bigon} \index{bigon!pseudo}

A {\sl red pseudo-bigon} is a closed disc in the red sphere whose boundary consists of one green arc and one blue arc, which may contain triple points in their interiors.
\end{definition}

\begin{figure}[ht!]
\begin{center} 
\epsfig{file=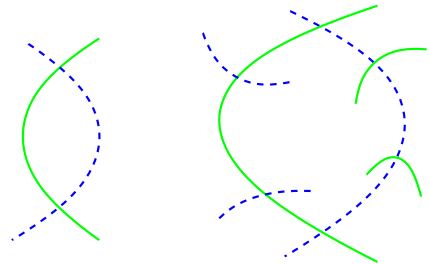, height=100pt}
\end{center}
\caption{Red pseudo-bigons}
\end{figure}

\begin{proposition}
If there are triple points, then the blue-green diagram contains at least four pseudo-bigons with disjoint interiors.
\end{proposition}

\begin{proof}
Consider a green curve which contains triple points. This divides a red sphere into two discs with disjoint interiors. Each disc must contain at least one blue arc, which divides each disc into a pair of pseudo-bigons with disjoint interiors.
\end{proof}

\begin{lemma} \label{innermost}
An innermost pseudo-bigon is a vertex free bigon.
\end{lemma}

\begin{proof}
Choose an innermost pseudo-bigon. If a boundary arc of this pseudo-bigon contains triple points in its interior, then we can find a double curve arc inside the pseudo-bigon with both endpoints on the same boundary arc which forms a smaller pseudo-bigon, a contradiction. So the boundary arcs of an innermost pseudo-bigon contain no triple points in their interiors. 

If an innermost pseudo-bigon has no triple points in the interiors of its boundary arcs, but triple points in its interior, then the interior contains a pair of intersecting circles, which create a pseudo-bigon, again giving a contradiction. Therefore an innermost pseudo-bigon has no triple points in its interior, and has boundary arcs with no triple points in their interiors, so is a vertex free bigon.
\end{proof}

\begin{corollary}{} \label{4lemma}
If a red sphere contains triple points, then the red sphere contains at least four vertex free red bigons with disjoint interiors.
\end{corollary}

\begin{remark}

Two vertex free disjoint red bigons may have translates which intersect in simple closed double curves.

\end{remark}

We now prove Lemma \ref{lemma2}. \\

\begin{proof}[Proof of Lemma \ref{lemma2}]
We deal with the possible non-triple point moves one by one.

\setcounter{case}{0}

\begin{case}{There is always a bigon orbit disjoint from a vanish or appear move}

At the singular time $t$ the move neighbourhood for a vanish or appear move contains only the orbit of a singular sphere, consisting of the orbit of a single point. The bigon orbits must all be contained in different sweepout spheres, so they are all disjoint from the move neighbourhood.
\end{case}

\begin{case}{There is always a bigon orbit disjoint from a cut or paste move}

At the singular time $t$ the move neighbourhood for a cut or paste move contains two discs of the same colour, which meet at a single point in their interiors. These discs are contained in the interiors of at most two vertex free bigon orbits, so there is a disjoint bigon orbit, as there are at least four vertex free bigon orbits by Corollary \ref{4lemma}.
\end{case}

\begin{case}{There is always a bigon orbit disjoint from a birth or death double curve move} \label{birthbigon}

At the singular time $t$ the move neighbourhood for a birth or death move contains two discs of different colours, which intersect at a single point. The orbit of these discs contains two red discs, which can be contained in at most two different vertex free red bigons. Therefore there is a red bigon disjoint from $N$, and so the orbit of this red bigon is a bigon orbit which is also disjoint from $N$.
\end{case}

\begin{case}{There is always a bigon orbit disjoint from a saddle move} \label{double} \label{disjoint_saddle}
 
At the move time $t$, the red spheres contain blue and green double curve circles, and precisely one figure eight double curve of each colour, which are in the same orbit under $G$. We will call the green figure eight curve $\g$. The blue figure eight curve is then $g\g$. All of these double curves may contain triple points.

The orbit of the move neighbourhood for the saddle move intersects the red sphere in a pair of discs, one of which contains the singular point of $\g$, and one of which contains the singular point of $g\g$. Furthermore, the intersection of the blue and green double curves with the discs consists only of a pair of green arcs which cross at the green singular point, and a pair of blue arcs which cross at the blue singular point, so it suffices to find a red bigon which is disjoint from the singular points of $\g$ and $g \g$.

The green figure eight curve $\g$ divides the red sphere in which it lies into three surfaces. Precisely two of these surfaces have closures which are discs, call these discs $D$ and $D'$. Let $D''$ be the closure of the complement of $D \cup D'$.

\begin{figure}[ht!] 
\begin{center}
\epsfig{file=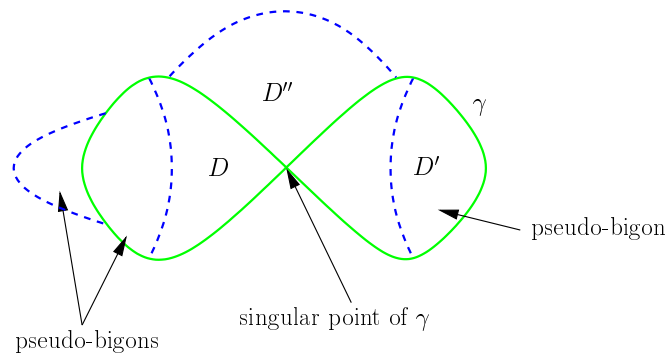, height=132pt,angle=0}
\end{center}
\caption{The figure-eight curve} 
\end{figure}

We consider the following three cases in turn: there is a blue arc in $D''$ with both endpoints on the same component of $\g - \{ \text{singular point} \}$, there is a blue arc in $D''$ with endpoints on different components of $\g - \{ \text{singular point} \}$, or else there are no blue arcs in $D''$.

If there is an arc in $D''$ with both endpoints on the same component of $\g$, then this arc creates a pseudo-bigon in $D''$. However, there must also be a blue arc in one of the discs $D$ or $D'$, creating a pseudo-bigon whose interior is disjoint from the one in $D''$. At most one of these pseudo-bigons can contain the blue singular point of $g\g$, so there is a pseudo-bigon disjoint from the move, and hence a vertex free bigon disjoint from the move.

If there is an arc in $D''$ with endpoints on both $D$ and $D'$, then there must be blue double arcs in both $D$ and $D'$, so each disc contains a pseudo-bigon. At most one of these pseudo-bigons can contain the blue singular point of $g\g$, so there is a pseudo-bigon disjoint from the move, and hence a vertex free bigon disjoint from the move.

If there are no blue arcs in $D''$, then $\g$ is triple point free, so the green and blue figure eight curves do not intersect any other double curves, and so may be contained in at most two vertex free red bigons. As there are at least four vertex free red bigons, there is a vertex free red bigon disjoint from the move neighbourhood.
\end{case}

This completes the proof of Lemma \ref{lemma2}, that there is always a vertex free red bigon disjoint from a non-triple point move.
\end{proof}

\section{Undermining the special cases}
\label{section:special cases}

Recall that a special case local maximum consists of a compound triple point birth followed by a compound triple point death. In this final section we complete the proof of the main result by proving the following theorem:

\begin{theorem} \label{lemma3}
Given a special case local maximum, we can either undermine the local maximum, or find an invariant unknotted curve.
\end{theorem}

This means that if there are no invariant unknotted curves, we can remove every local maximum, so the minimax sweepout contains no triple points, a contradiction, by Lemma \ref{lemma1}.

Consider a special case local maximum. Let $\a$ be the bigon orbit which is created by the compound triple point birth move, and let $\b$ be the bigon destroyed in the compound triple point death move. At a time during the local maximum in between the compound moves, both bigons are present in the sweepout, and are double curve free. There are only finitely many ways in which such a pair of bigons may intersect, and we deal with each possibility in turn.

There is the trivial case when $\a$ and $\b$ are the same bigon orbit. In this case the configurations before and after the local maximum are isotopic, so we can replace the partial sweepout containing the moves that occur during the local maximum with a partial sweepout that is an isotopy and contains no moves, thus undermining the local maximum.
There are three remaining cases to consider depending on how many vertices the bigon orbits $\a$ and $\b$ have in common. If $\a$ and $\b$ have no vertices in common, then they must be disjoint, so we can just postpone the first compound move till after the second compound move, thus undermining the local maximum. The remaining two cases to consider are when $\a$ and $\b$ have either three or six vertices in common.

Section \ref{section:three_vertices} deals with the case in which $\a$ and $\b$ have three vertices in common. In this case we give explicit constructions of modifications which undermine the local maxima.

Section \ref{section:six_vertices} deals with the case in which $\a$ and $\b$ have six vertices in common. If $\a \cup \b$ is connected, we produce an invariant unknotted curve. If $\a \cup \b$ is not connected, then we show how to undermine the local maximum.

We now fix some notation which we will use for the rest of this section.

\begin{definition}{\sl Special case modification neighbourhood}
\index{modification neighbourhood!special case}
\index{special case modification neighbourhood}

Suppose $(M, \phi)$ is a sweepout containing a special case local maximum. Suppose that $N_I$ is a modification neighbourhood with the following properties:

\begin{itemize}

\item $N_I$ contains the move neighbourhoods for the compound triple point moves in its interior, and is disjoint from all the other move neighbourhoods of the sweepout. 

\item $N_t$ is a tubular neighbourhood for $\a_t \cup \b_t$ for the times $t$ in between the compound triple point moves. 

\end{itemize}

Then we say that $N_I$ is a {\sl special case modification neighbourhood} for the special case local maximum.
\end{definition}

The special case modification neighbourhood may be connected or disconnected, and may have components that are not $3$-balls. The sweepout only changes by an isotopy outside the special case modification neighbourhood during the local maximum.

\begin{figure}[ht!]
\begin{center}
\epsfig{file=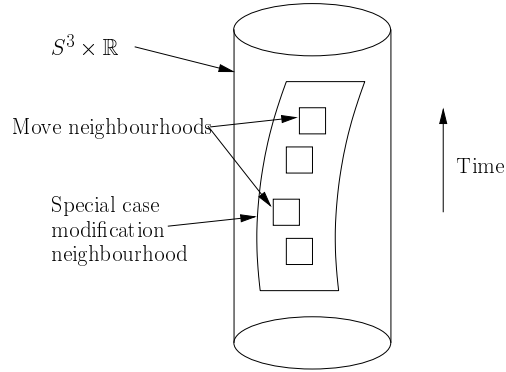, height=140pt}
\end{center}
\caption{Schematic picture of the special case modification neighbourhood}
\end{figure}

\begin{lemma}
Let $(M, \phi)$ be a sweepout containing a special case local maximum. Then we can choose a special case modification neighbourhood for the given special case local maximum.
\end{lemma}

\begin{proof}
Choose move neighbourhoods for the local maximum, and label the corresponding time intervals as shown below.

\begin{figure}[ht!]
\begin{center}
\epsfig{file=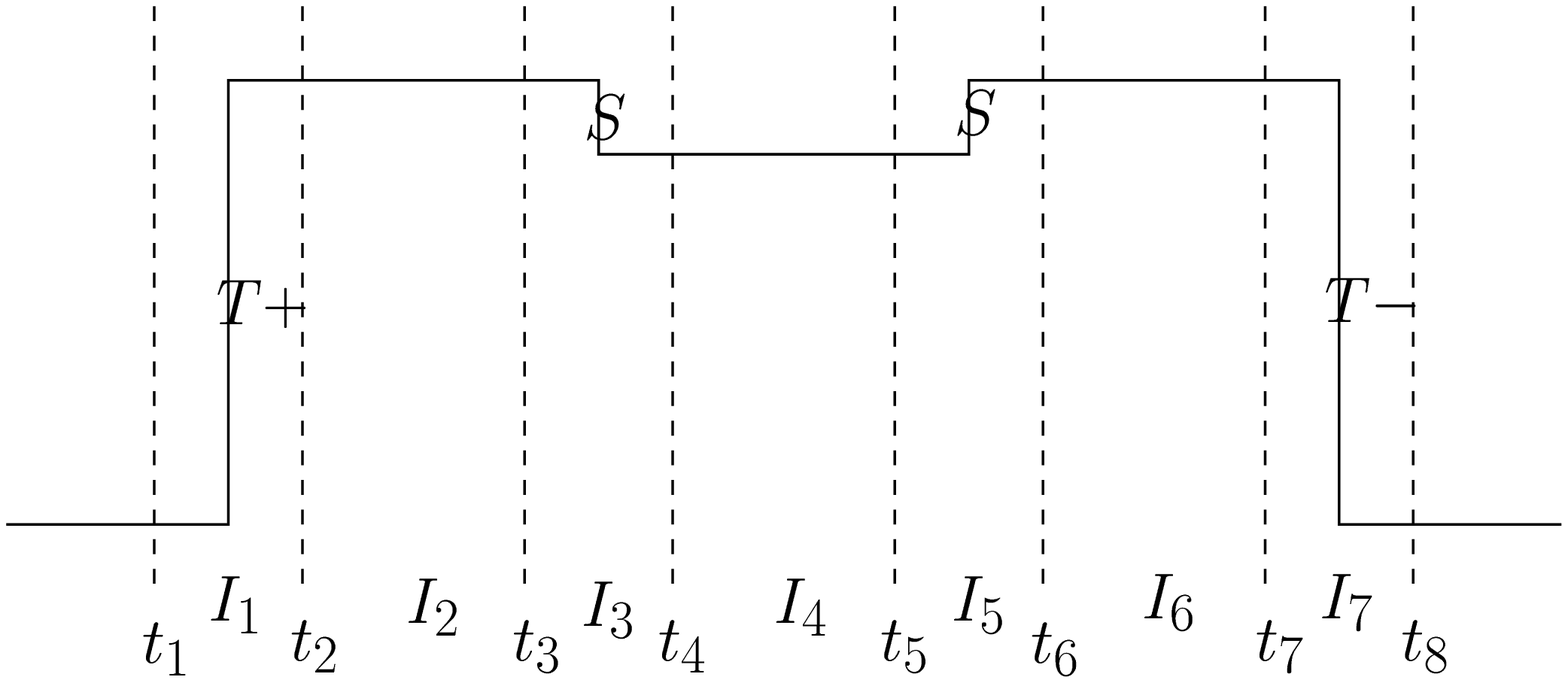, height=120pt}
\end{center}
\caption{The graphic for a special case local maximum}
\end{figure}

We will show how to construct the special case modification neighbourhood going forward in time. A similar construction will work going back in time. Let $J_{I_5}$ be the move neighbourhood for the second saddle move. Let $K_{I_7}$ be the move neighbourhood for the triple point death move.

Choose a thin tubular neighbourhood for $\a \cup \b$ at time $t_4$. The pre-image sweepout for $S^3_{I_4}$ is a product, so by Lemma \ref{lemma:isotopy} we can choose a compatible product neighbourhood $N_{I_4}$ so that $N_t$ is a thin tubular neighbourhood for $\a \cup \b$ for each $t \in I_4$.

The saddle move occurring during the time interval $I_5$ is determined by a choice of saddle disc $D$ at time $t_5$. The saddle move neighbourhood $J_{t_5}$ is a tubular neighbourhood of this disc $D$ at time $t_5$. The saddle disc $D$ is parallel to the bigon $\b$, which is contained in $N_{t_5}$, so at time $t_5$ we can choose the thin tubular neighbourhood $N_{t_5}$ of $\a \cup \b$ to be large enough to contain $J_{t_5}$. 

Let $N'_{t_5} = N_{t_5} - J_{t_5}$. The partial sweepout $(S^3 \cross I_5) - J_{I_5}$ is a product, so by Lemma \ref{lemma:isotopy} we can extend $N'_{t_5}$ to a continuously varying family $N'_{I_5}$, so that the intersection of each $N_t$ with the sweepout surfaces is isotopic to the intersection of the sweepout surfaces with $N_t$. In particular, if we define $N_{I_5}$ to be $N'_{I_5} \cup J_{I_5}$, then as $\d N_t$ is contained in $S^3 - J_t$, the boundary pattern only changes up to isotopy in $\d N_t$ for the time interval $I_5$.

There are no moves during the time interval $I_6$, so by Lemma \ref{lemma:isotopy}, we can extend $N_I$ across this interval, so that the sweepout in $N_t$ only changes up to isotopy.

As $N_{t_6}$ contains the football created by the saddle move, $N_{t_7}$ will also contain this football. At time $t_7$, the move neighbourhood for $K_{t_7}$ is a tubular neighbourhood of this football, so we can choose $N_{t_7}$ to be large enough to contain $K_{t_7}$. Now by using Lemma \ref{lemma:isotopy} in the same way as for the saddle move, we can extend $N_I$ across the time interval $I_7$, so that $N_I$ has constant boundary pattern.
\end{proof}

It will be convenient to fix notation for the different components of the orbit of each bigon.

\begin{notation}
Let $A$ be the red bigon component of $\a$, and $B$ the red bigon component of $\b$. We will write $N^A$ for the component of $N$ that contains the red bigon $A$.
\end{notation}

\subsection{Three vertices in common} \label{section:three_vertices}

In this section we show: 

\begin{lemma}

A special case local maximum in which $\a$ and $\b$ have exactly three vertices in common can be undermined.

\end{lemma}

We shall prove this by giving explicit constructions of partial sweepouts, with the same boundary as the pre-image sweepout for the special case modification neighbourhood, but which contain fewer triple points. 

The red bigon $A$ will be disjoint from its images under $g$, though it will intersect images of the red bigon $B$. There are two different ways in which the bigon-orbits can share three vertices. Case 1 is if the red bigons $A$ and $B$ have a vertex in common, we will call this the {\sl coplanar} case. Case 2 is if the red bigon $A$ shares a vertex with the green bigon $gB$. We will call this the {\sl orthogonal} case. If the red bigon $A$ shares a vertex with the blue bigon $g^2B$, then the red bigon $B$ shares a vertex with the green bigon $gA$, so this is just the time reverse of Case 2. In both of these cases we will show that we can modify the sweepout to reduce the local maximum.

\setcounter{case}{0}

\begin{case}
{Coplanar}

\begin{figure}[ht!]
\begin{center}
\epsfig{file=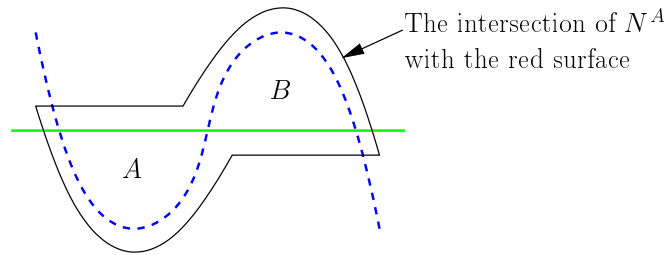, height=95pt}
\end{center}
\caption{Coplanar red bigons in the red sphere}
\end{figure}

In this case $N_I$ is the orbit of a $3$-ball which is disjoint from its images. Each component of $N_I$ contains an image of $A \cup B$ under $G$. It suffices to describe the modification in the component $N^A_I$ of $N_I$ that contains the red bigons $A \cup B$, we do the corresponding equivariant modifications in the other components of $N_I$. During the local maximum, we may assume that the red and blue surfaces remain fixed, and only the green surface moves. The red and blue spheres intersects $N^A$ in a pair of discs, which intersect in a single double arc. It is convenient to give $N^A$ a $\{ \text{disc} \} \cross I$ product structure, so that the red disc is horizontal, and the blue disc is vertical. This is shown below:

\begin{figure}[ht!]
\begin{center}
\epsfig{file=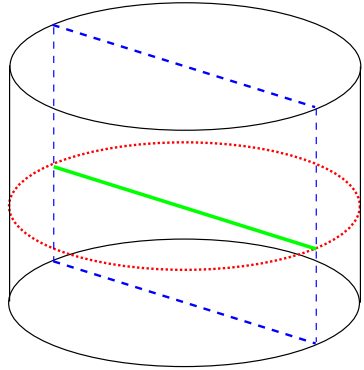, height=140pt}
\end{center}
\caption{Case 1: The red and blue surfaces in $N^A$}
\end{figure}

We can choose the product structure so that during the central product interval, the green surface is also vertical. This is shown below:

\begin{figure}[ht!]
\begin{center}
\epsfig{file=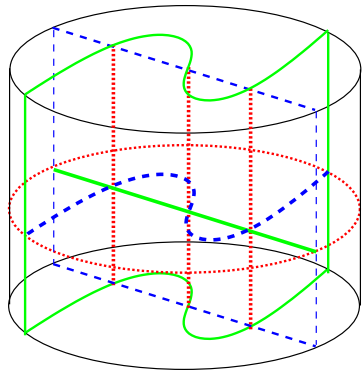, height=140pt}
\end{center}
\caption{Case 1: The green surface in $N^A$ during the local maximum}
\end{figure}

In order to describe the partial sweepout for the modification, it will be convenient to draw pictures of $N^A$ split into four parts, along the red and blue surfaces. We will call these pieces {\sl quadrants}.

\begin{figure}[ht!]
\begin{center}
\epsfig{file=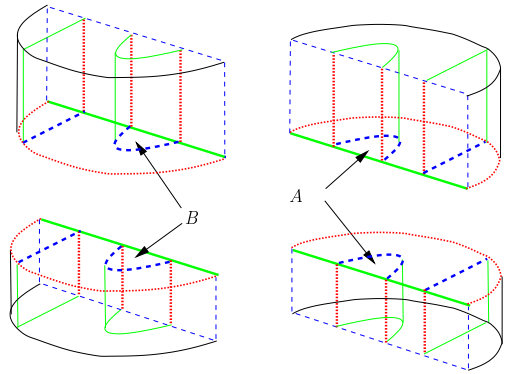, height=140pt}
\end{center}
\caption{The configuration between the compound moves} \label{picture360}
\end{figure}

We can now construct the configurations before and after the local maximum. Figure \ref{picture361}(a) shows the configuration before the local maximum. It is produced from the configuration in Figure \ref{picture360} by removing the bigon $A$. Figure \ref{picture361}(b) shows the configuration after the local maximum. It is produced from the configuration in Figure \ref{picture360} by removing the bigon $B$.

\begin{figure}[ht!]
\begin{center}

\mbox{

\subfigure[Before the local maximum]{\epsfig{file=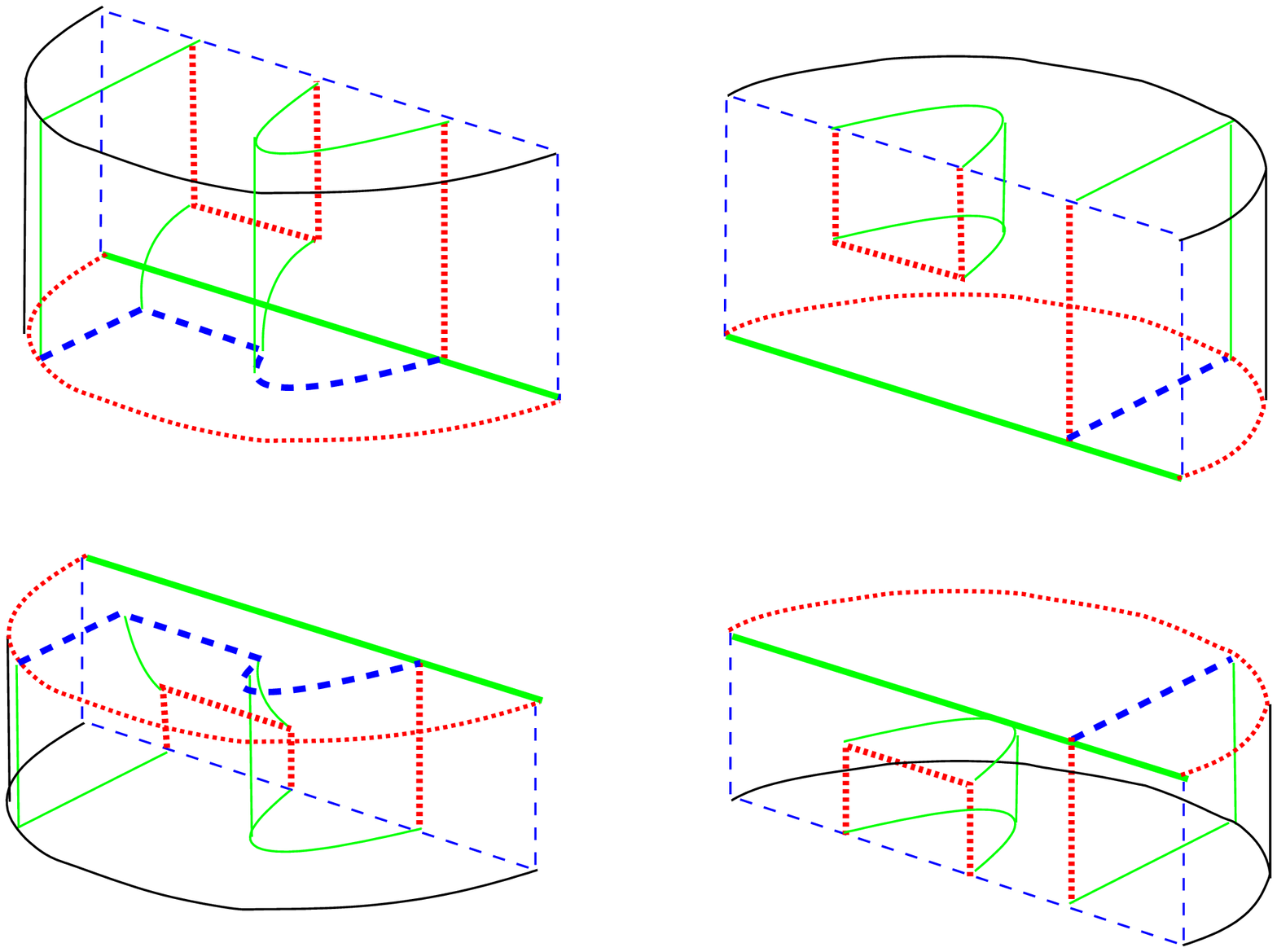, height=120pt}}\hspace{0.25in} 

\subfigure[After the local maximum]{\epsfig{file=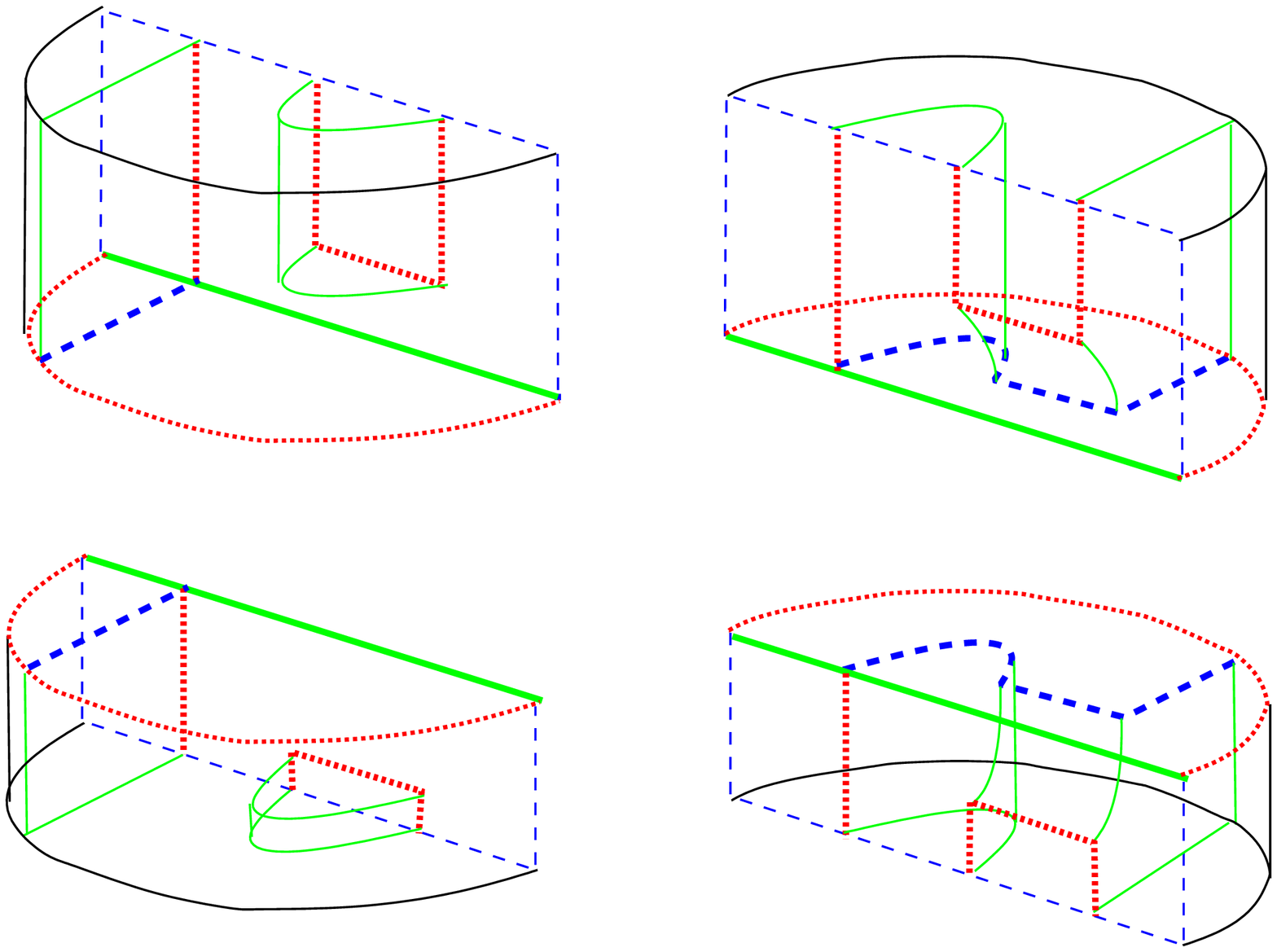, height=120pt}}

}

\end{center}
\caption{The original partial sweepout in $N^A$} \label{picture361}
\end{figure}

We now describe a new partial sweepout with the same boundary as the pre-image sweepout for the special case modification neighbourhood, which we will use to modify the sweepout. Notice that the configurations in Figure \ref{picture361}(a) and Figure \ref{picture361}(b) are not isotopic, as for example, the top right quadrant contains different numbers of discs.

During the original partial sweepout, the components of the green surface in each quadrant are discs, and so the configuration is determined (up to isotopy) by the boundaries of these discs. We can make the boundary curves isotopic by applying appropriate saddle moves to the configuration.  Figure \ref{picture361}(a) shows the result of applying a saddle move parallel to $A$ to Figure \ref{picture360}. Figure \ref{picture361}(b) shows the result of applying a saddle move parallel to $B$ to Figure \ref{picture360}. This preserves the fact that the green surfaces in each quadrant are discs, so it suffices to check that the configurations are isotopic on the boundaries of the quadrants. The configurations on the boundaries of each quadrant are clearly isotopic, and furthermore, we can do these isotopies so that they agree on the surfaces which are identified. Also, the curves consist of the same combination of coloured arcs, so we can do an isotopy without changing the number of coloured arcs, so there will always be exactly one triple point.

\begin{figure}[ht!]
\begin{center}

\mbox{
\subfigure[After a saddle move parallel to~$B$]{\epsfig{file=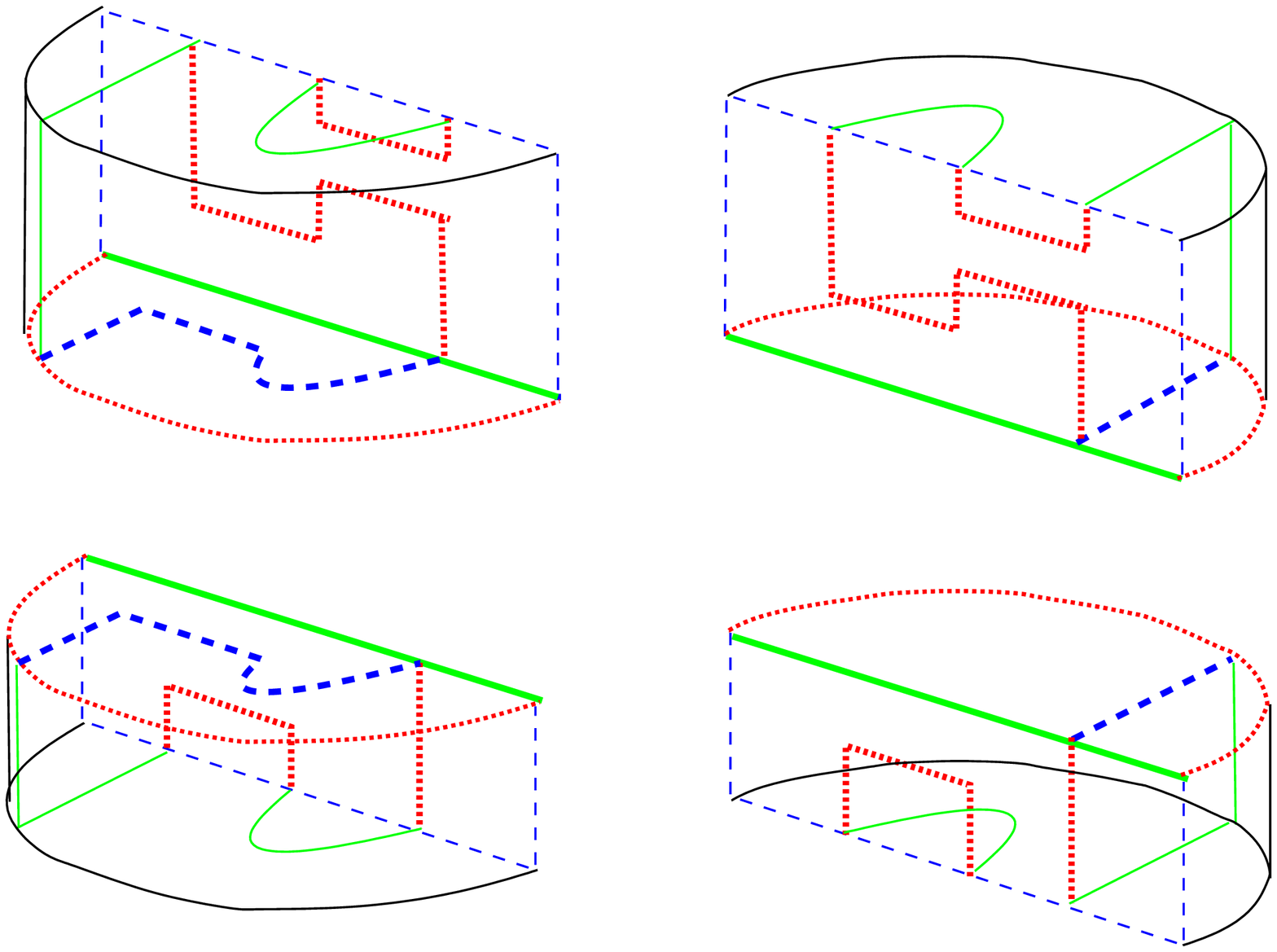, 
height=130pt}}\hspace{0.22in}
\subfigure[Before a saddle move parallel to~$A$]{\epsfig{file=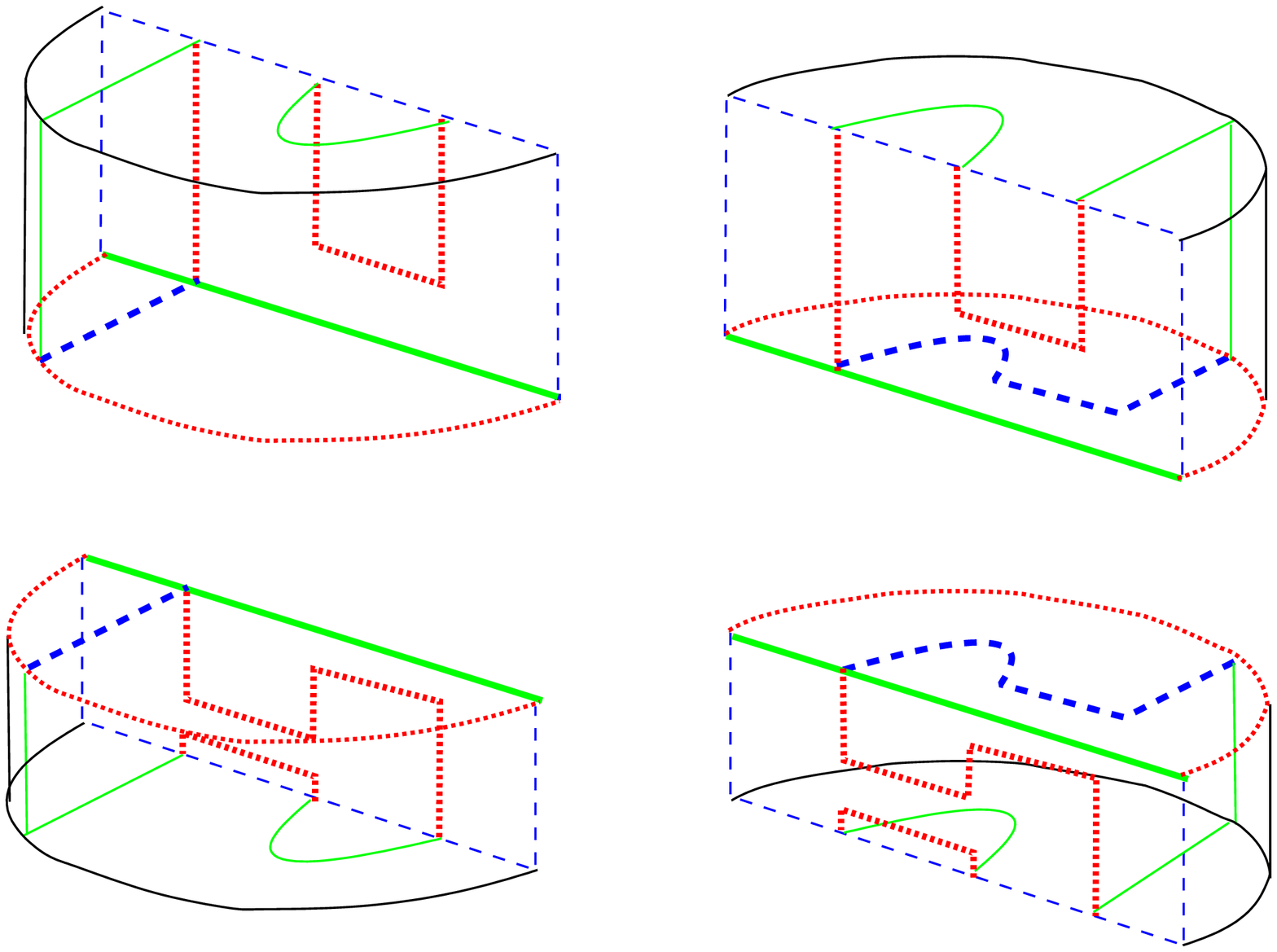, 
height=130pt}}
}

\end{center}
\caption{The new partial sweepout in $N^A$} \label{picture363}
\end{figure}

So we can construct a partial sweepout starting at the configuration in \ref{picture361}(a), and then applying a saddle move to create the configuration in Figure \ref{picture363}(a). This is then isotopic to Figure \ref{picture363}(b), and we can then apply another saddle move to create the final configuration Figure \ref{picture361}(b). This partial sweepout has the same boundary as the initial partial sweepout in $N^A$, but has only one triple point, so we can use it to modify the sweepout to undermine the local maximum.
\end{case}

\begin{case}{Orthogonal}

We now consider the case in which the images of $A$ and $B$ which share a vertex have different colours. We may assume that $A$ shares a vertex in common with $gB$. If $A$ shares a vertex with $g^2 B$, then $gA$ shares a vertex with $B$, so we could just swap the labels on $A$ and $B$. The neighbourhood $N$ is then the orbit of a $3$-ball which is disjoint from its images, each component of which contains an image of $A \cup gB$.

Let $N^A$ be the component of $N$ which contains $A \cup gB$. It suffices to describe how to modify the sweepout inside $N^A$ only. Figure \ref{picture354} shows the configuration in $N^A$ at a time in between the compound moves. We have drawn the red disc as horizontal, and the green disc as vertical. The red and green discs intersect in a single double arc.

\begin{figure}[ht!]
\begin{center}
\epsfig{file=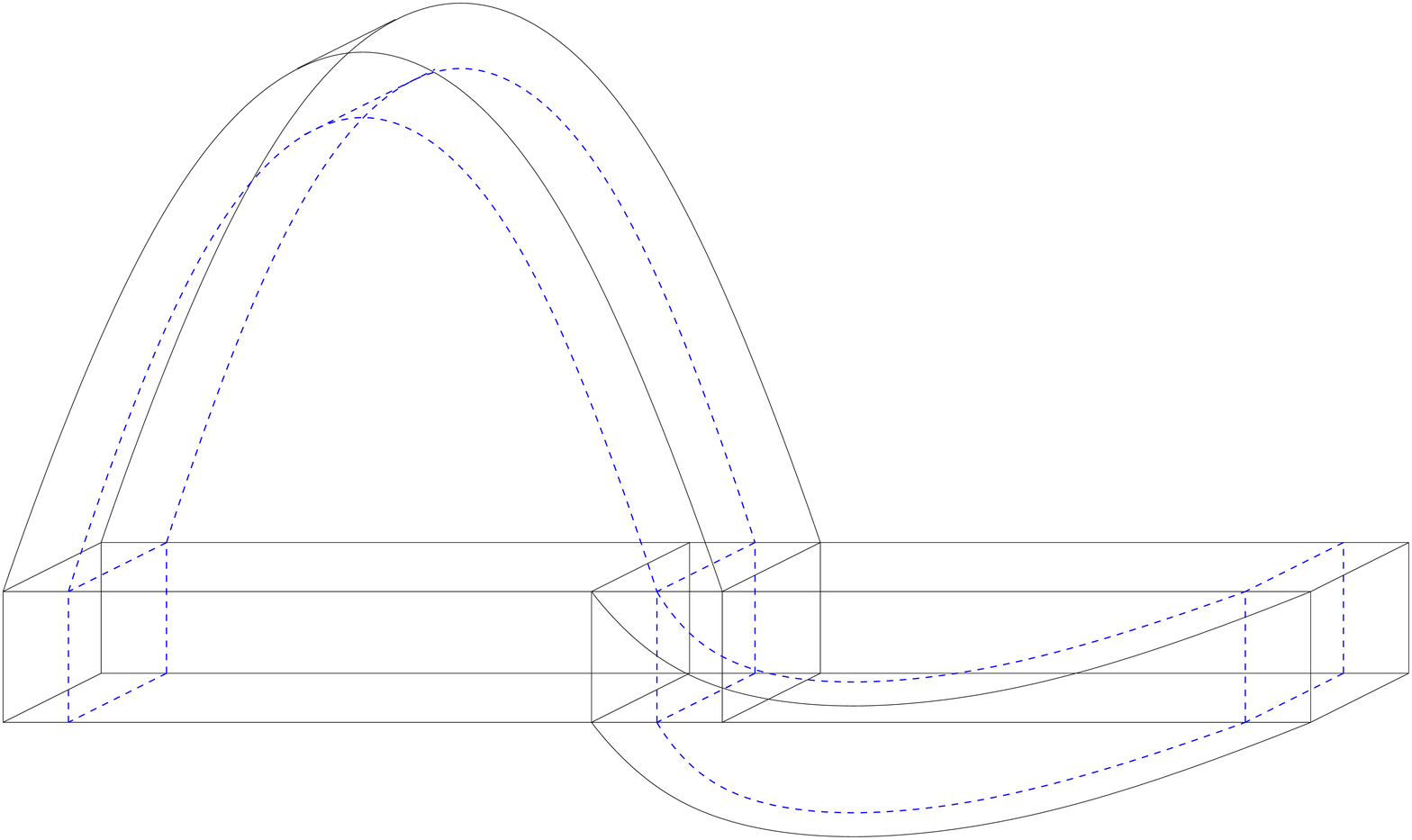, height=160pt}
\end{center}
\caption{Case 2: The blue surface inside $N^A$} \label{picture473}
\end{figure}

In the picture above we have $N^A$ as a thin tubular neighbourhood of a horizontal bigon and a vertical bigon. The picture below is isotopic to the one above, except we have drawn $N^A$ as $\{ \text{disc} \} \cross I$, with the red surface as horizontal, and the green surface as vertical.

\begin{figure}[ht!]
\begin{center}
\epsfig{file=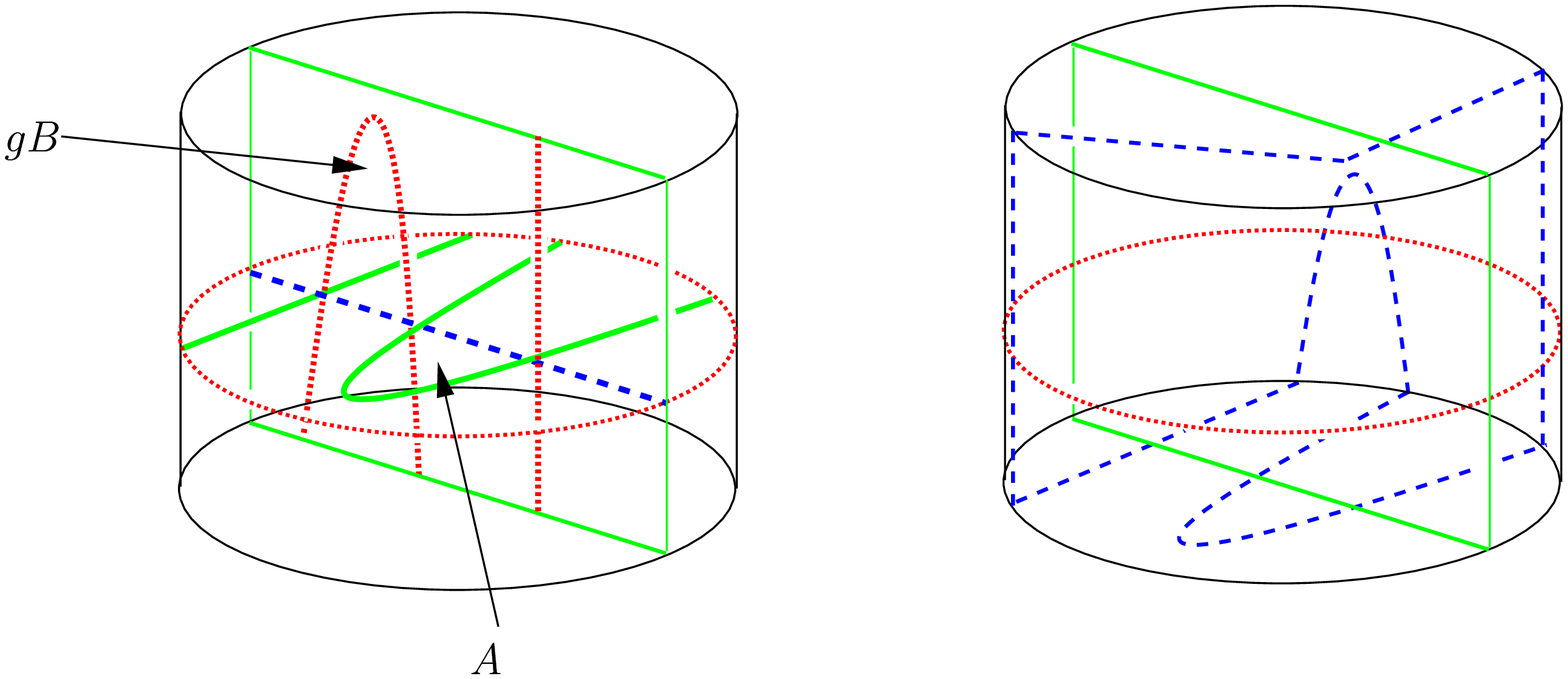, height=160pt}
\end{center}
\caption{Case 2: $N^A$. This picture is isotopic to the one above.} \label{picture354}
\end{figure}

As in the previous case, it will be convenient to draw pictures of $N^A$ divided into four quadrants by cutting it along the red and green surfaces, as shown in Figure \ref{picture355}.

\begin{figure}[ht!]
\begin{center}
\epsfig{file=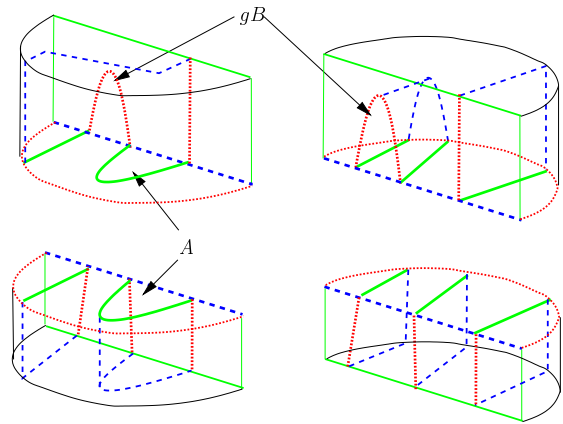, height=160pt}
\end{center}
\caption{Case 2: $N^A$} \label{picture355}
\end{figure}

We may assume that the red and green surfaces in $N^A$ remain fixed, and only the blue sphere moves, so going forward in time, we remove the green bigon $gB$ by pushing the blue surface down, as shown in Figure \ref{picture357}(b), and going back in time, we remove the red bigon $A$ by pushing the blue surface sideways, as shown in Figure \ref{picture357}(a).

\begin{figure}[ht!]
\begin{center}

\mbox{

\subfigure[Before the local maximum]{\epsfig{file=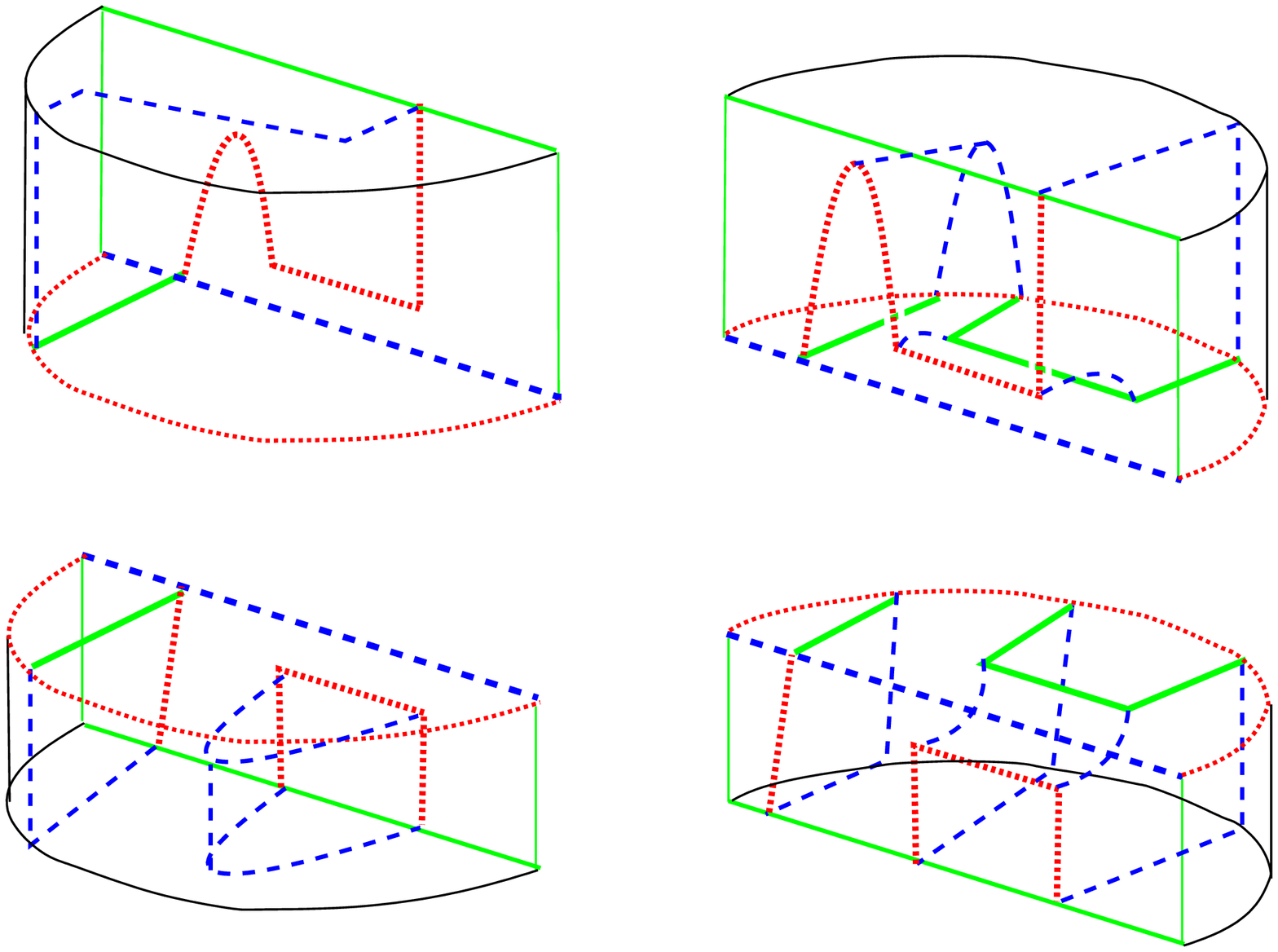, height=130pt}}\hspace{0.25in}

\subfigure[After the local maximum]{\epsfig{file=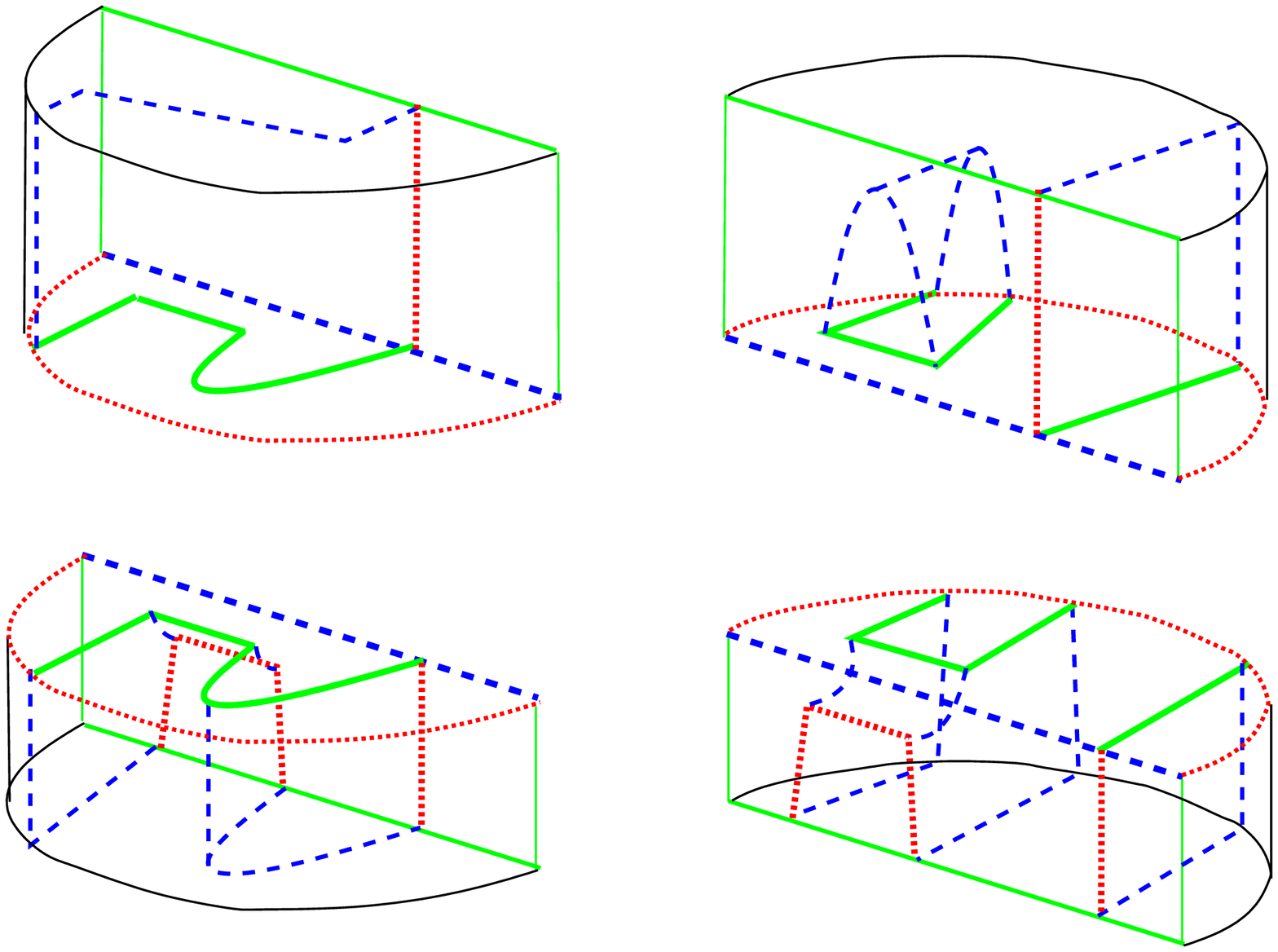, height=130pt}}

}

\end{center}
\caption{The original sweepout, after an isotopy} \label{picture357}
\end{figure}

Figure \ref{picture357}(a) and Figure \ref{picture357}(b) are not isotopic, as for example, the top right quadrant contains different numbers of blue discs. However, by applying an appropriate saddle move to Figure \ref{picture357}(a) and Figure \ref{picture357}(b), we can create a pair of isotopic configurations.

We now describe the new partial sweepout. Start with the initial configuration shown in Figure \ref{picture357}(a). A saddle move is then applied to the green double curves, using a saddle disc parallel to the green disc $gB$. The resulting configuration is shown in Figure \ref{picture365}(a). As all the components of the blue surface in each quadrant are still discs, it suffices to draw their boundaries. Figure \ref{picture365}(a) is isotopic to Figure \ref{picture365}(b), as the boundaries of the blue discs in each quadrant are isotopic. Furthermore, we can do these isotopies so that they agree on the surfaces in the boundary of each quadrant which are identified, and without changing the number of arcs of each colour. In particular, this means that there is only one triple point in the new partial sweepout in $N^A$. Figure \ref{picture365}(b) is produced from Figure \ref{picture357}(b) by doing a saddle move in the bottom quadrants, using a saddle disc parallel to red disc $A$, so we end at the correct final configuration.

\begin{figure}[ht!]
\begin{center}

\mbox{

\subfigure[After the first saddle]{\epsfig{file=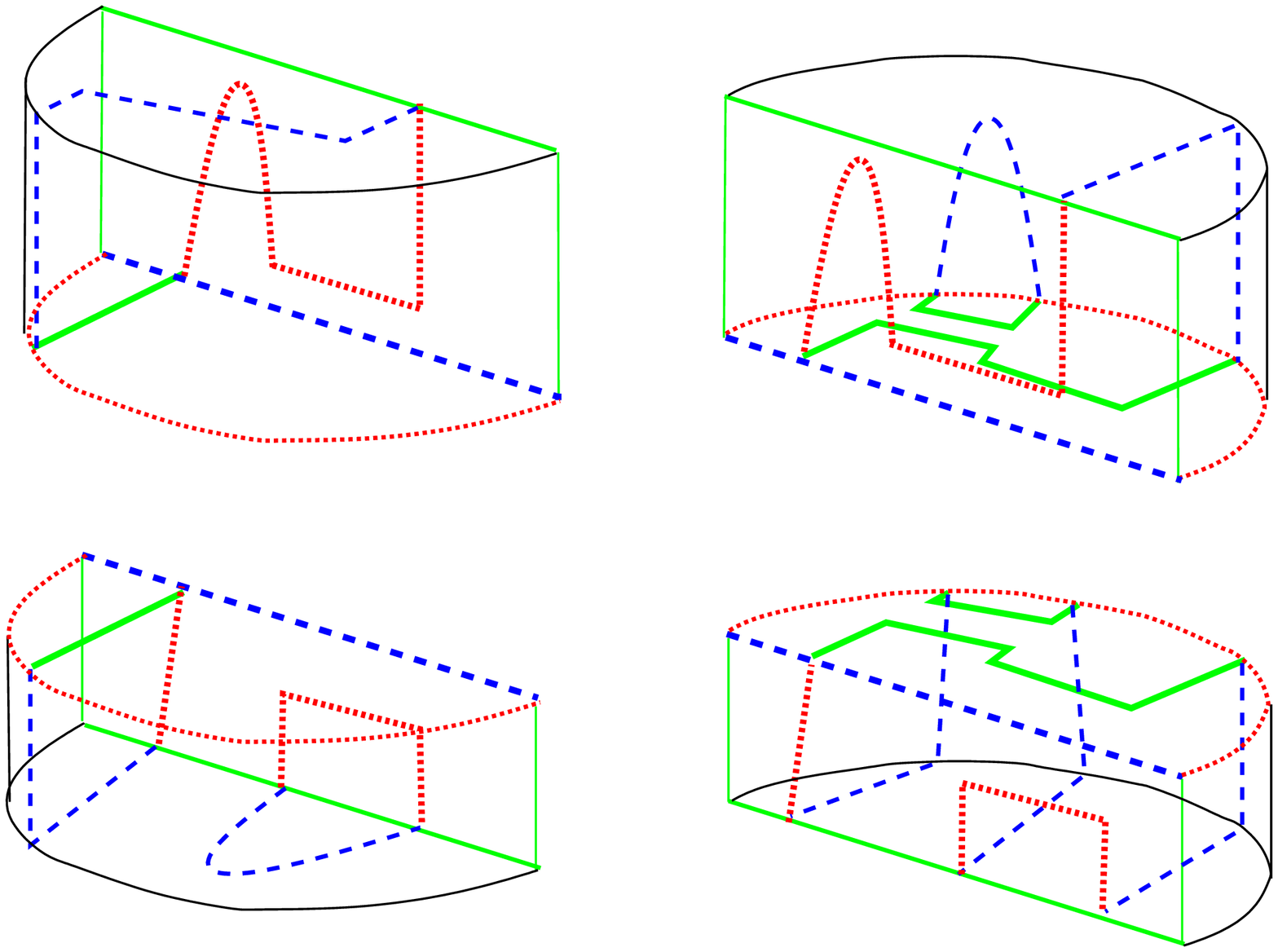, height=130pt}}\hspace{0.25in}

\subfigure[Before the second saddle]{\epsfig{file=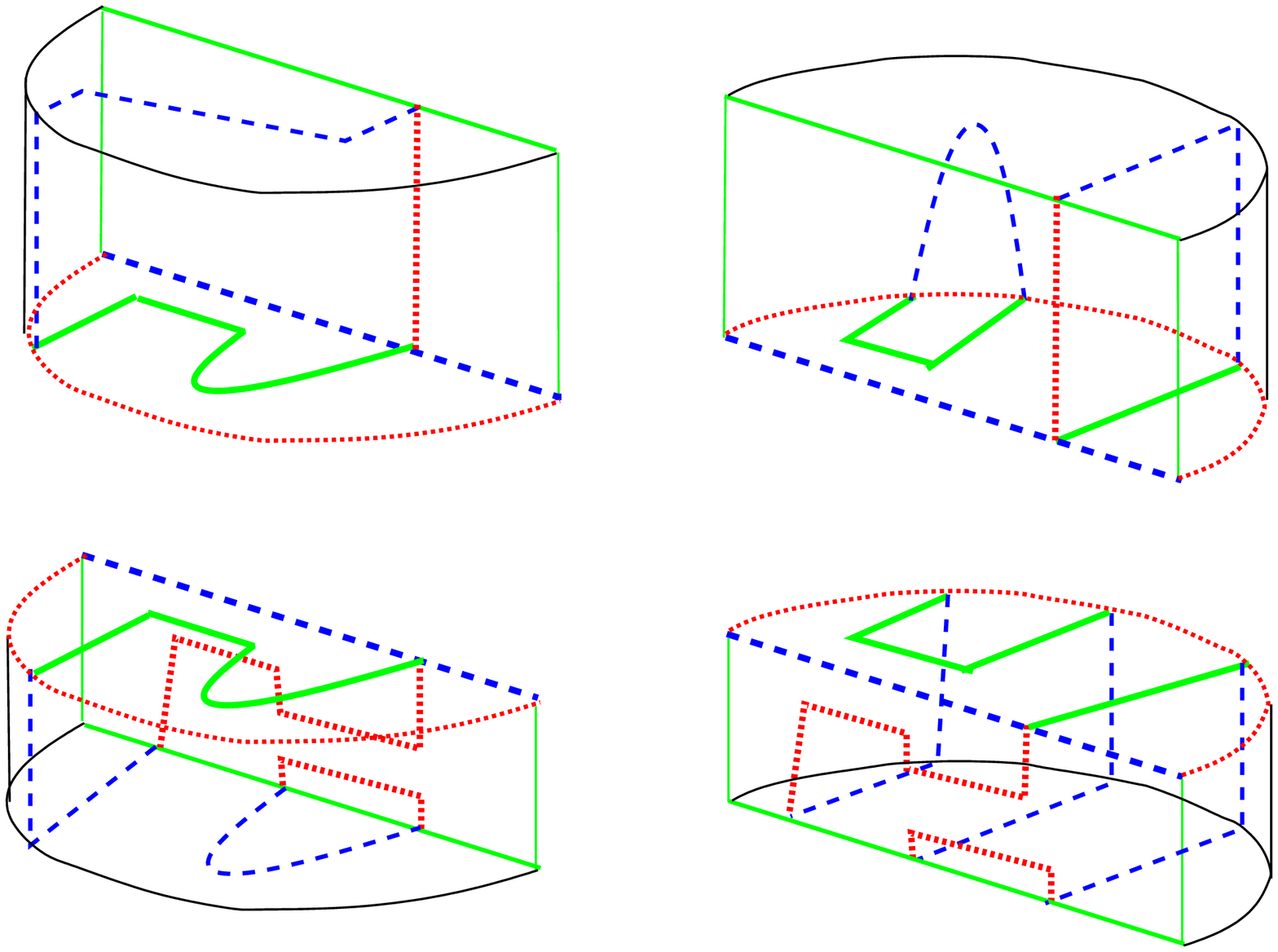, height=130pt}}

}

\end{center}
\caption{The new partial sweepout} \label{picture365}
\end{figure}
\end{case}

\subsection{Six vertices in common} \label{section:six_vertices}

In this section we deal with the final case, when the bigon orbits $\a$ and $\b$ in the special case local maximum have all six vertices in common.

There are three different cases to consider depending on how many edge orbits the bigon orbits $\a$ and $\b$ have in common. If they have both edge orbits in common, then $\a$ and $\b$ are in fact the same bigon orbit, and the configurations before and after the local maximum are isotopic, so the local maximum can be undermined by replacing the sweepout in the special case modification neighbourhood by an isotopy. The remaining two cases are when $\a$ and $\b$ share one edge orbit in common, which we discuss in Section \ref{section:one_edge}, or no edge orbits in common, which we discuss in Section \ref{section:no_edge}. 

In order to undermine these local maximum, we shall prove a general undermining lemma, which we now describe. Let $N_I$ be a modification neighbourhood in which $N$ is the orbit of a $3$-ball which is disjoint from its images. We call an intersection of a double curve with $\d N$ a {\sl double point} of the boundary pattern. A {\sl boundary bigon} is a disc in $\d N$ whose boundary consists of exactly two arcs of intersection with the sweepout surfaces, which contain no double points in their interiors. Given a boundary bigon, we can use this disc as a saddle disc for a saddle move, which simplifies the boundary pattern by removing a pair of double points. We say a boundary pattern is {\sl saddle reducible} if all double points can be removed by a sequence of saddle moves using boundary bigons.

Given a modification neighbourhood $N_I$ which has a saddle reducible boundary pattern, and no triple points in $N_0$ and $N_1$ we show that we can construct a new partial sweepout with the same boundary as the original one, but which contains no triple points for the whole time interval.

In some cases, the special case modification neighbourhood will have saddle reducible boundary, and we can apply the lemma directly. In other cases, $N$ will not be the orbit of a $3$-ball, but we will show how to extend the neighbourhood to a larger one which is a $3$-ball with a saddle reducible boundary pattern.

We now give a precise definition of what it means for a boundary pattern to be saddle reducible, and we show that if a modification neighbourhood has a saddle reducible boundary pattern then we can replace the partial sweepout in this neighbourhood by one which contains no triple points.

\begin{definition}{\sl Double points}
\index{double points}

We say a transverse intersection of a double curve with $\d N$ is a {\sl double point} in $\d N$.
\end{definition}

\begin{definition}{\sl A boundary-bigon}
\index{boundary bigon} \index{bigon!boundary}

A {\sl boundary-bigon} is a disc in $\d N$ whose boundary consists of a pair of arcs of different colours, with no double points in their interiors.
\end{definition}

\begin{remark}

Given a simple closed curve which bounds an innermost disc in the boundary pattern, we can use this as a cut disc for a cut move, which removes the curve of intersection. Similarly, given a boundary-bigon in the boundary pattern, we can use it as a saddle disc for a saddle move, producing a boundary pattern with two fewer double points.

\end{remark}

\begin{definition}{\sl Saddle reducible}
\index{saddle-reducible}

Let $N_I$ be a modification neighbourhood, where $N$ is the orbit of a
$3$-ball which is disjoint from its images. We say that the boundary
pattern of $N_I$ is {\sl saddle reducible} if there is a time $t$ at
which there is sequence of cut moves and saddle moves, using cut discs
and saddle discs contained in $\d N_t$, which remove all the curves in
the boundary pattern of $N_t$.
\end{definition}

\begin{remark}
The boundary pattern of a modification neighbourhood $N_I$ only changes up to isotopy, so if it is saddle reducible at some time $t \in I$, it will be saddle reducible for all times $t \in I$. 
\end{remark}

Not all boundary patterns are saddle reducible, for example:

\begin{figure}[ht!]
\begin{center}
\epsfig{file=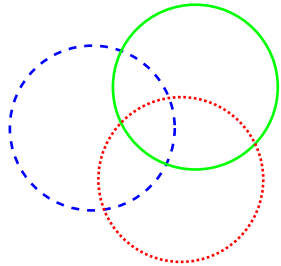, height=100pt}
\end{center}
\caption{This boundary pattern is not saddle reducible}
\end{figure}

\begin{lemma} \label{lemma:undermining}
Let $(M, \phi)$ be a sweepout. Let $N_I$ be a modification
neighbourhood, where $N_t$ is the orbit of a $3$-ball which is
disjoint from its images. Suppose that $N_I$ contains no triple points
at the beginning and the end of the time interval, and has a saddle reducible boundary pattern.

Then there is a new sweepout $(M, \phi)$, which is the same as $(M, \phi)$ outside of $N_I$, and which has no triple points in $N_I$.
\end{lemma}

\begin{proof}
We construct a partial sweepout with the same boundary as the pre-image sweepout for $N_I$, but which contains no triple points. 
Let $C_I$ be a subset of $N_I$, such that $C_I$ has a compatible product structure, and $C_t \cong \orbit S^2 \cross [0,1]$ is a collar neighbourhood for $\d N_t$, so that the sweepout surfaces intersect $C_t$ in a product $\{ \text{boundary pattern} \} \cross [0,1]$. We may assume that $\d N_t$ is $\orbit S^2_t \cross \{ 0 \}$. Let $C'_I$ be the ``outer half'' of $C_I$, so that $C'_t \cong \orbit S^2 \cross [0,\half]$.

As the intersection of $C_t \cross \{s\}$ with the sweepout surfaces is the same for all $s \in [0,1]$, the closure of $N_t - C'_I$ also has saddle reducible boundary pattern, so there is a sequence of cut moves and saddle moves which make the sweepout surfaces disjoint from $\orbit S^2 \cross \{ \half\}$. As these moves use cut discs and saddle discs contained in $\orbit S^2 \cross \{ \half\}$, we may assume that the move neighbourhoods of the moves are contained in $C_I$. The sphere orbit $\orbit S^2 \cross \{ \half \}$ is now the boundary of a $3$-ball orbit containing no triple points. We can remove all the sweepout spheres inside the $3$-ball orbit by first removing double curves using cut and death moves, and then removing the remaining disjoint spheres by vanish moves.

We can construct a similar partial sweepout starting with the final configuration instead. Use the same sequence of cut moves and saddle moves as before, to make the sweepout surfaces disjoint from $\orbit S^2 \cross \{ \half \}$. The resulting collection of closed $2$-spheres inside the three ball orbit bounded by $\orbit S^2 \cross \{ \half\}$ may be different, but remove them all using cut, death and vanish moves. The resulting configuration will be the same as the one above, so we can use the first partial sweepout, followed by the time reverse of the second, to create a partial sweepout with the same boundary as $(\phi^{-1}(N_I), \phi)$, but which has no triple points.
\end{proof}

The lemma below shows that if there are not many double points, we do not need to explicitly construct the boundary pattern.

\begin{lemma} \label{lemma:four boundary points}
If the boundary pattern has at most four double points, then it is saddle reducible.
\end{lemma}

\begin{proof}
If the boundary pattern has only two double points, then all the curves in the boundary are simple, except for a single pair of curves that intersect each other twice, creating four bigons. We can remove any simple closed curves of intersection from the bigons by using cut moves. Then using any of these bigons as a saddle disc removes the double points.

Suppose there are four double points. If all the double points lie on a single pair of curves, as shown in Figure \ref{picture581}(a), then there is a bigon we can use to reduce the number of double points. This reduces us to the previous case, so this boundary pattern is saddle reducible. If the double points lie on two disjoint pairs of curves, Figure \ref{picture581}(b), then the boundary pattern is also saddle reducible. Otherwise, two of the double points lie on a single pair of curves, Figure \ref{picture581}(c), at least one of which contains another double point. However, curves must contain even numbers of double points, so this curve must in fact contain four double points, coming from two distinct intersecting curves, so there is a bigon we can use to reduce the number of double points. 
\end{proof}

\begin{figure}[ht!]
\begin{center}
\epsfig{file=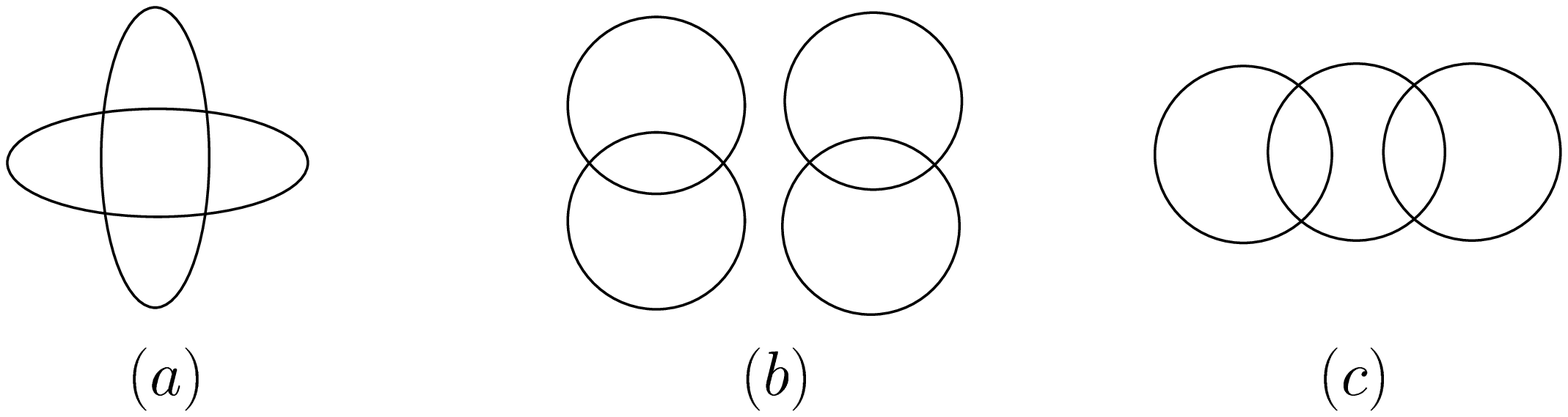, height=72pt}
\end{center}
\caption{All possible boundary patterns with four double points}
\end{figure} \label{picture581}

In some cases $N_t$ will not be an orbit of a $3$-ball, and in these cases we will show that we can extend the neighbourhood to produce a larger neighbourhood which is the orbit of a $3$-ball.  We will then examine the boundary pattern, and show that in each case it is saddle reducible. 

The following lemmas will allow us to make some simplifying assumptions about the configuration outside $N_t$ during the local maximum.

\begin{lemma} \label{lemma:three components of N in S}
Suppose $(M, \phi)$ is a sweepout containing a special case local maximum, with bigon-orbits $\a$ and $\b$, and special case modification neighbourhood $N_I$. Then $N_t \cap S_t$ has at most three components.
\end{lemma}

\begin{proof}
The neighbourhood $N_t$ is the union of two regular neighbourhoods of bigon-orbits. Each component of a regular neighbourhood of a bigon-orbit intersects the sweepout surfaces in three discs, one of each colour. So a regular neighbourhood of $\a$ intersects the red spheres $S_t$ in three discs, as does a regular neighbourhood of $\b$. However, the bigon-orbits $\a$ and $\b$ have a common vertex-orbit, so each red disc must intersect at least one other red disc, so there are at most three components of $N_t \cap S_t$.
\end{proof}

\begin{lemma} \label{disjoint}
Suppose $(M, \phi)$ is a sweepout containing a special case local maximum, with bigon-orbits $\a$ and $\b$, and special case modification neighbourhood $N_I$. If there is a component of the blue-green graph disjoint from $N$, then there is a vertex-free bigon-orbit disjoint from $N_I$.
\end{lemma}

\begin{proof}
Every component of the blue-green graph contains at least four red pseudo-bigons with disjoint interiors. By Lemma \ref{lemma:three components of N in S}, $N_t \cap S_t$ has at most three components, and $N_t$ may intersect at most three of these red pseudo-bigons, so there is a pseudo-bigon disjoint from $N_t$. This pseudo-bigon contains a vertex-free red bigon disjoint from $N_t$. As the configuration in $S^3 - N_t$ only changes by an isotopy, there is therefore a continuously varying family of disjoint bigons for the time interval $I$.
\end{proof}

\begin{lemma} \label{lemma:double_curve_free_bigon}
Suppose $(M, \phi)$ is a sweepout containing a special case local maximum, with bigon-orbits $\a$ and $\b$. Suppose there is a vertex-free red bigon $C$ at some time during the local maximum, which contains some simple closed double curves.
Then we can modify the sweepout to remove the simple closed double curves during the corresponding local maximum in the new sweepout. 
\end{lemma}

\begin{proof}
The double curves in the interior of $C$ are disjoint from the move neighbourhoods of the moves occurring during the special case local maximum, so we can remove them for the duration of the local maximum by doing cut and death moves before the compound triple point birth, and then birth and paste moves after the compound triple point death. 
\end{proof}

Finally, we prove a lemma which will be useful in showing that certain sets are disjoint from their images under $G$.

\begin{lemma} \label{lemma:bigons with same vertices}
Let $A_1$, \ldots, $A_n$ be double curve free bigons of any colours, such that each bigon $A_i$ has the same vertices as $A_1$. Let $A$ be the union of these bigons. Then $A$ is disjoint from its images under $G$. 
\end{lemma}

\begin{proof}
If any image of one of the $A_i$ intersects $A$, then as the bigons are double curve free, it must intersect $A$ in at least one vertex. But this means that there is a bigon which shares a vertex with one of its images, a contradiction, by Corollary \ref{corollary:bigon orbit disjoint boundaries}.
\end{proof}

This lemma guarantees that we can choose a regular neighbourhood of such a set of bigons which is disjoint from its images.

\subsubsection{One edge in common} \label{section:one_edge}

There are two cases to consider, coplanar and orthogonal. In the coplanar case the red bigons $A$ and $B$ share a common edge. In the orthogonal case the red bigon $A$ shares a common edge with either the green or blue image of $B$.

\setcounter{case}{0}

\begin{case}{Coplanar}

In this case $A \cup B$ is a red disc which is disjoint from its images under $G$. In the picture below, $A$ and $B$ share a green edge in common. If they share a blue edge in common, then the picture is the same, but with the colours blue and green reversed.

\begin{figure}[ht!]
\begin{center}
\epsfig{file=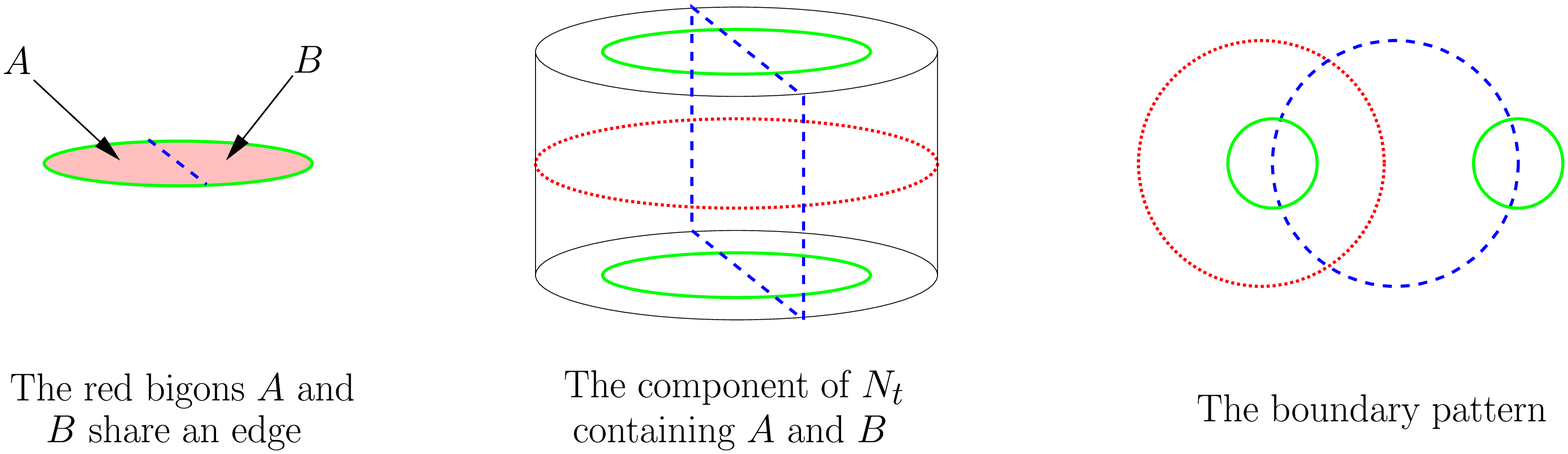, height=100pt}
\end{center}
\caption{$N^A$} \label{picture488}
\end{figure}

The orbit of $N^A$ is a union of three disjoint $3$-balls, with
boundary pattern shown in Figure \ref{picture488}, which is saddle
reducible, so the local maximum can be undermined by Lemma
\ref{lemma:undermining}.
\end{case}

\begin{case}{Orthogonal}

We may assume that the red bigon $A$ shares an edge in common with the green bigon $gB$. If $A$ shares an edge with the blue bigon $g^2B$, then we can reduce to the previous case by swapping the labels on $A$ and $B$.

\begin{figure}[ht!]
\begin{center}
\epsfig{file=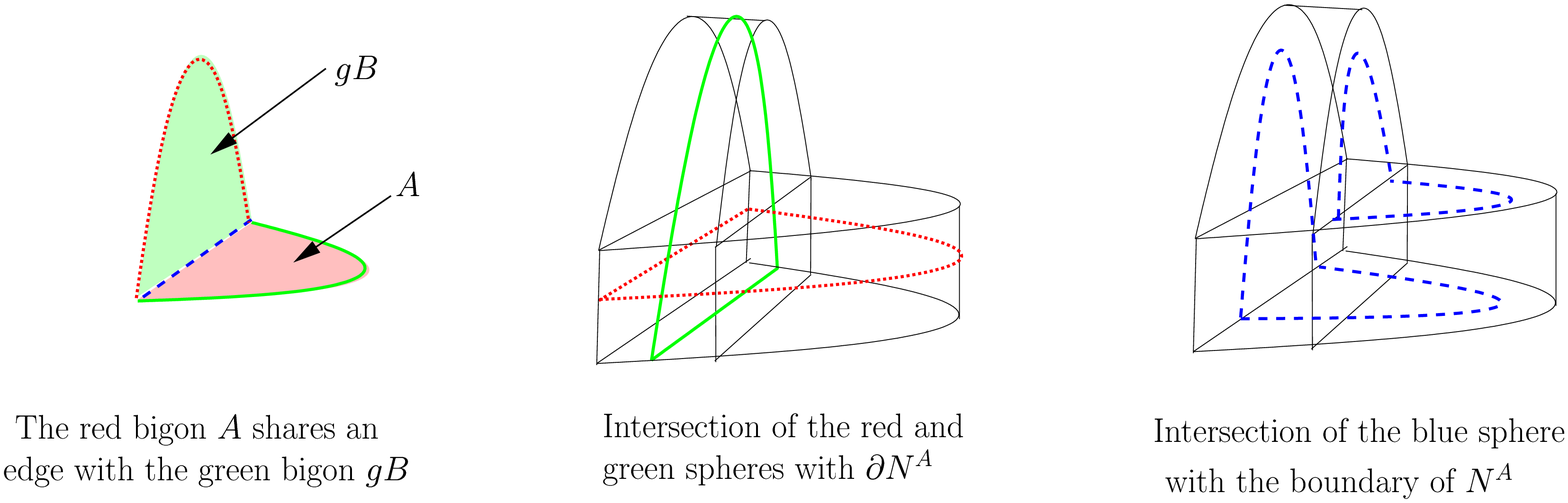, width=\hsize}
\end{center}
\caption{$N^A$} \label{disc}
\end{figure}

\begin{figure}[ht!]
\begin{center}
\epsfig{file=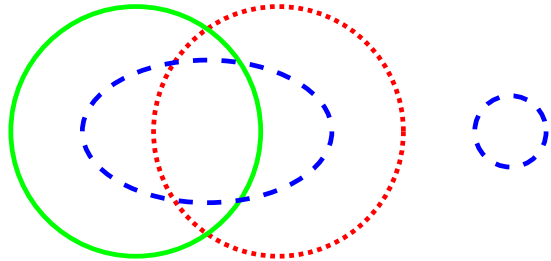, height=100pt}
\end{center}
\caption{The boundary pattern} \label{picture489}
\end{figure}

The orbit of $N^A$ is three disjoint $3$-balls, and $N^A$ has a saddle
reducible boundary pattern, so we can undermine the local maximum by
Lemma \ref{lemma:undermining}.
\end{case}

\subsubsection{No edges in common} \label{section:no_edge}

There are two different ways in which the bigon-orbits $\a$ and $\b$ may have six vertices in common, but no edges in common. This is illustrated in the picture below. We have coloured tubular neighbourhoods of one bigon black, and tubular neighbourhoods of the other bigon pink.

\begin{figure}[ht!]
\begin{center}
\epsfig{file=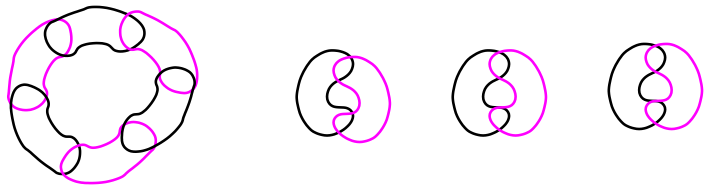, height=70pt}
\end{center}
\caption{The bigon-orbits $\a$ and $\b$ may share six vertices in two different ways}
\end{figure}

If the union of the bigon-orbits is connected, we can show that the action of $G$ is standard.

\begin{lemma}
Let $\a$ and $\b$ be a pair of bigon-orbits which have all six vertices in common, and with $\a \cup \b$ connected. Then there is an invariant unknotted curve. 
\end{lemma}

\begin{proof}
Let $\g$ be the orbit of the green edges of $A$ and $B$. Then $\g$ is an invariant simple closed curve consisting of a pair of arcs of each colour. The blue and green arcs lie in the red sphere. The red arcs lie in the boundary of bigons whose other edge lies in the red sphere, so we can isotope them into the red sphere across these bigons, so $\g$ is unknotted, as it is isotopic to a curve which lies in a sphere.
\end{proof}

\begin{figure}[ht!]
\begin{center}
\epsfig{file=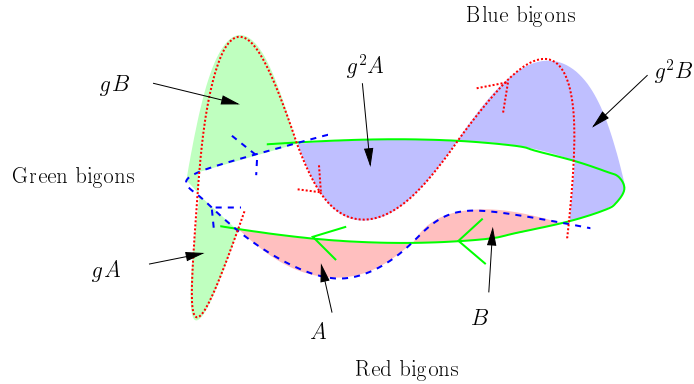, height=145pt}
\end{center}
\caption{An example of an unknotted invariant curve}
\end{figure}

So the remaining situations we need to analyse are when $N_I$ has three connected components. There are two main cases to consider, the coplanar case, when the red bigons $A$ and $B$ have  a pair of vertices in common, and the orthogonal case, when the red bigon $A$ has a pair of vertices in common with either the green bigon $gB$ or the blue bigon $g^2B$.

\setcounter{case}{0}
\begin{case}
{Coplanar}

In this case the red bigons $A$ and $B$ have two vertices in common, but no edges in common. Then $N^A$ contains a simple closed green curve, which we will call $a$, and a simple closed blue curve, which we shall call $gb$.

\begin{figure}[ht!]
\begin{center}
\epsfig{file=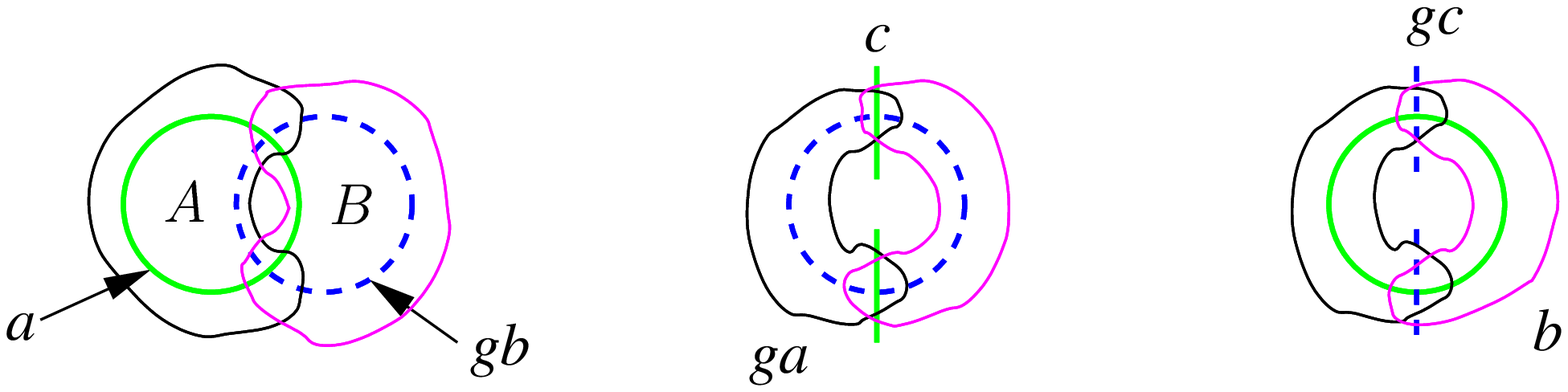, height=70pt}
\end{center}
\caption{The intersection of $N$ with the red sphere}
\end{figure}

In the picture above, we have drawn a thin tubular neighbourhood of $\a$ in black, and a thin tubular neighbourhood of $\b$ in pink.

The double curves $b \cup ga$ divides the red spheres into components, at least two of which are discs. Only one of these discs can intersect $N^A$, so one of the curves bounds a disc $\D$ in a red sphere which is disjoint from $N^A$. If there are vertex-free bigon-orbits disjoint from $N^A$, then we can use them to undermine the local maximum, so we may assume that this disc cannot contain triple points in its interior. This means that $\D$ contains a single simple double arc. By pinching off any simple closed double curves in $\D$ for the duration of the local maximum, we may change the sweepout so that $\D$ contains no simple closed double curves in its interior. 

We may assume that $b$ bounds the disc $\D$, the case where $ga$ bounds the disc $\D$ is the same, but with the colours blue and green swapped round.

\begin{figure}[ht!]
\begin{center}
\epsfig{file=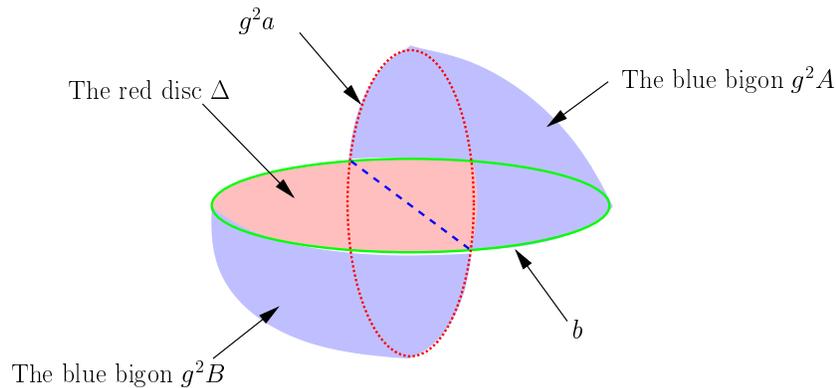, height=147pt}
\end{center}
\caption{The red disc $\D$, together with the blue bigons $g^2A$ and $g^2B$.}
\end{figure}

The surface $\D \cup gA \cup gB$ shown above is disjoint from its images under $G$, by Lemma \ref{lemma:bigons with same vertices}, as it is a union of double curve free bigons which share the same vertices. Therefore we can choose a thin tubular neighbourhood of the bigons which is a $3$-ball containing $N^A$, which is also disjoint from its images under $G$. The boundary pattern for this $3$-ball has exactly two double points, so will be saddle reducible, by Lemma \ref{lemma:four boundary points}. So the special case local maximum can be undermined by Lemma \ref{lemma:undermining}.
\end{case}

\begin{case}
{Orthogonal}

In this case the red bigon $A$ shares both vertices in common with either the green or blue image of $B$. We may assume that $A$ shares vertices with the green bigon $gB$. If $A$ shares vertices with the blue bigon $g^2B$, then we can reduce to the previous case by swapping the labels on $A$ and $B$.

The neighbourhood $N^A$ contains a simple closed blue curve, call this $gb$. It also contains a green double arc, which we shall call $a$ and a red double arc, which we shall call $g^2c$.

\begin{figure}[ht!]
\begin{center}
\epsfig{file=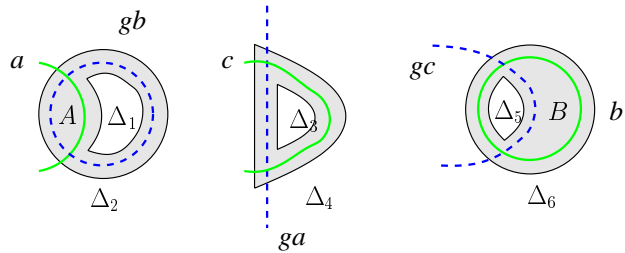, height=95pt}
\end{center}
\caption{The intersection of $N$ with the red sphere} \label{picture584}
\end{figure}

The tubular neighbourhood of $\a \cup \b$ intersects the red spheres in three annuli. These annuli divide the spheres into pieces, which we shall label $\D_1, \ldots, \D_6$, as shown in Figure \ref{picture584}. To be precise, note that each annulus in $N_t \cap S_t$ has a boundary component which is disjoint from the double curves, call this the inner boundary component, and call the other one the outer boundary component. Give the component of $S_t - N_t$ which borders the inner boundary component of $N^A$ the label $\D_1$. Similarly give the components of $S_t - N_t$ which border the inner components of $gN^A$ and $g^2 N^A$ the labels $\D_2$ and $\D_3$. Now give the component of $S_t - N_t$ which borders the outer component of $N^A$ the label $\D_4$. Similarly, give the components of $S_t - N_t$ which border the outer components of $gN^A$ and $g^2 N^A$ the labels $\D_5$ and $\D_6$.
Some of the components of $S_t - N_t$ may have more than one label, but at least two of them must be innermost discs which have only one label.

At least one of the innermost discs is \emph{not} labelled $\D_4$. We will show how to extend $N^A$ to a neighbourhood which is a $3$-ball with saddle reducible boundary by adding an innermost disc labelled $\D_i$, where $i \not = 4$. We now deal with these cases one by one. Note that the case that $\D_1$ is innermost is the same as the case that $\D_5$ is innermost, by swapping colours. Similarly the case for $\D_2$ is the same as the case for $\D_6$.

We may assume that there are no bigon-orbits disjoint from $N_I$, as if there were, we could use them to undermine the local maximum. Therefore the red bigons $\D_1$, $\D_2$ and $\D_3$ may not contain components of the graph of double curves with triple points which are disjoint from the orbit of $N_I$.

\begin{subcase}{One of $\D_1$ or $\D_5$ is innermost}
\end{subcase}

We may assume we are in the case $\D_1$, as the other case is the same as this one, but with the colours green and blue reversed. We will assume that $\D_1$ is vertex-free, otherwise there is a bigon disjoint from $N_I$ which we can use to undermine the local maximum. Furthermore, we may assume that $\D_1$ is double curve free, by removing the double curves.

The three bigons $A$, $gB$ and $\D_1$ share a common pair of vertices, so $A \cup gB \cup \D_1$ is disjoint from its images under $G$, by Lemma \ref{lemma:bigons with same vertices}. We will use a thin tubular neighbourhood of $A \cup gB \cup \D_1$ as our modification neighbourhood, which we may assume contains $N_I$. This is a $3$-ball which is disjoint from its images under $G$.

\begin{figure}[ht!]
\begin{center}
\epsfig{file=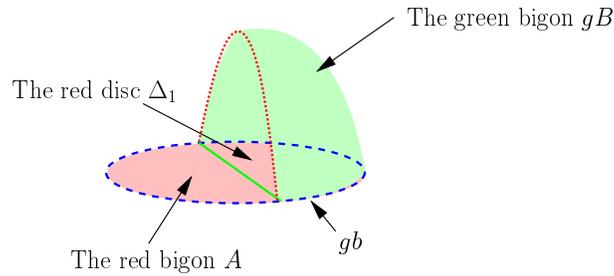, height=105pt}
\end{center}
\caption{The red bigons $A$ and $\D_1$, together with the green bigon $gB$}
\end{figure}

The boundary pattern has four double points, and so is saddle reducible by Lemma \ref{lemma:four boundary points}.

\begin{subcase}{$\D_3$ is innermost}
\end{subcase}

Again, we may assume that $\D_3$ is double curve free. Note that the bigons $gA \cup g^2B \cup \D_3$ have the same vertices, so the union of these bigons is disjoint from its images under $G$. We can choose a tubular neighbourhood of $gA \cup g^2B \cup \D_3$, which is disjoint from its images under $G$, and which contains $N_I$. This is a $3$-ball, which has four double points, as there are two double arcs which intersect the boundary. Therefore it has a saddle reducible boundary pattern, by Lemma \ref{lemma:four boundary points}.

\begin{figure}[ht!]
\begin{center}
\epsfig{file=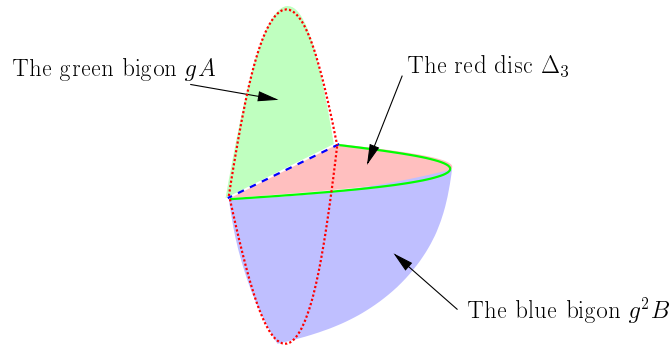, height=130pt}
\end{center}
\caption{The red bigon $\D_3$, together with the green bigon $gA$ and the blue bigon $g^2B$}
\end{figure}

\begin{subcase}{Either $\D_2$ or $\D_6$ is innermost}
\end{subcase}

We will assume $\D_2$ is innermost, as the other case is the same, except with the colours blue and green reversed.

\begin{figure}[ht!]
\begin{center}
\epsfig{file=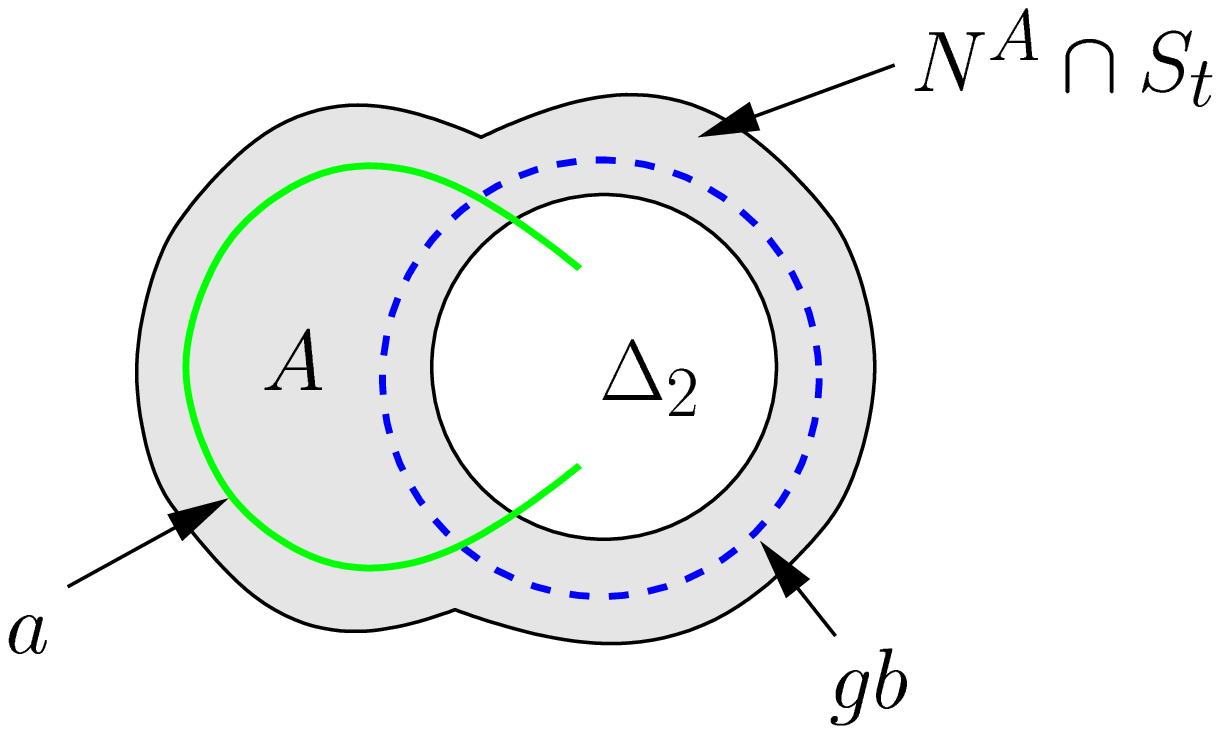, height=90pt}
\end{center}
\caption{$\D_2$ is innermost}
\end{figure}

Any extra triple points inside the innermost disc create disjoint bigons which we could use to undermine the local maximum, so we may assume that $a$ has exactly two triple points, and divides the disc $\D_2$ into two bigons. We can now choose a thin tubular neighbourhood of $A \cup gB \cup \D_2$ which contains $N^A$, and which is disjoint from its images.

\begin{figure}[ht!]
\begin{center}
\epsfig{file=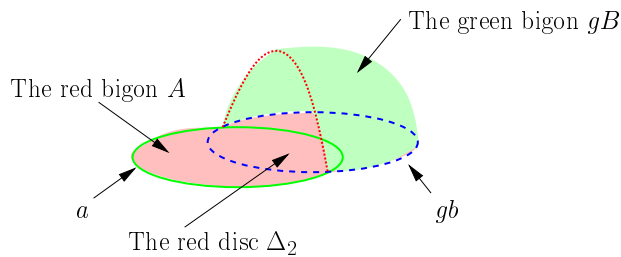, height=95pt}
\end{center}
\caption{The red bigon $A$, togther with the red disc $\D_2$ and the green bigon $gB$}
\end{figure}

A tubular neighbourhood of the bigons shown above is a $3$-ball, whose boundary contains precisely two double points, coming from the red double arc. So the boundary pattern is saddle reducible, by Lemma \ref{lemma:four boundary points}.
\end{case}

We have shown that if there are no unknotted invariant circles, then every special case local maximum can be undermined, so there is a sweepout with no triple points, a contradiction.

This completes the proof that any $\Z_{3}$ action on $S^{3}$ is standard.




\end{document}

%% file: gtoutput.tex
\def\ifplaintex{\expandafter\ifx\csname documentclass\endcsname\relax}

\ifplaintex 
\else
\fi

\expandafter\ifx\csname beginpicture\endcsname\relax
\expandafter\ifx\csname documentclass\endcsname\relax
\input pictex \else
\input prepictex \input pictex \input postpictex \fi\fi

\def\gt{{\mathsurround=0pt\it $\cal G\mskip-2mu$eometry \&\ 
$\cal T\!\!$opology}}        

\def\gtp{{\mathsurround=0pt\it $\cal G\mskip-2mu$eometry \&\ 
$\cal T\!\!$opology $\cal P\!$ublications}}  

\def\lognumber#1{\def\thelognumber{#1}}
\def\volumenumber#1{\def\thevolumenumber{#1}}
\def\papernumber#1{\def\thepapernumber{#1}}
\def\volumeyear#1{\def\thevolumeyear{#1}}

\def\pagenumbers#1#2{\def\startpage{#1}\def\finishpage{#2}}
\def\published#1{\def\publishdate{#1}}
\def\proposed#1{\def\theproposer{#1}}
\def\seconded#1{\def\theseconders{#1}}
\def\received#1{\def\receiveddate{#1}}
\def\revised#1{\def\reviseddate{#1}}
\def\accepted#1{\def\accepteddate{#1}}

\def\asciiaddress#1{\def\theasciiaddress{#1}}

\long\def\asciiabstract#1{\long\def\theasciiabstract{#1}}

\let\\\par\let\thelognumber\relax
\let\thevolumenumber\relax\let\thepapernumber\relax
\let\thevolumeyear\relax\let\thesamplenumber\relax\let\startpage\relax
\let\finishpage\relax\let\publishdate\relax\let\receiveddate\relax
\let\reviseddate\relax\let\accepteddate\relax\let\theasciititle\relax
\let\theasciiauthors\relax\let\theasciiaddress\relax
\let\theasciiabstract\relax
\let\theasciiemail\relax\let\theshortauthors\relax\let\theshorttitle\relax

\long\def\maketitlep{   

\count0=\startpage

\gt\hfill      
\beginpicture
\setcoordinatesystem units <0.33truein, 0.33truein> point at 2.2 0.9
\setplotsymbol ({$\cal G$})
\plotsymbolspacing=9truept
\circulararc 315 degrees from 0 1 center at 0 0
\setplotsymbol ({$\cal T$})
\circulararc 315 degrees from 1 -1 center at 1 0
\endpicture
\break
{\small\ifx\thesamplenumber\relax 
Volume \else Sample
\fi\thevolumenumber\ (\thevolumeyear)
\startpage--\finishpage\nl
Published: \publishdate}
\vglue 0.5truein plus 0.4fil minus 0.1truein

{\parskip=0pt\leftskip 0pt plus 1fil\def\\{\par\smallskip}{\ifplaintex\large
\else\Large\fi\bf\thetitle}\par\medskip}   

\vglue 0pt plus 0.1fil 

{\parskip=0pt\leftskip 0pt plus 1fil\def\\{\par}{\sc\theauthors}
\par\medskip}

\vglue 0pt plus 0.1fil 

{\small\parskip=0pt\let\newline\\
{\leftskip 0pt plus 1fil\def\\{\par}{\sl\theaddress}\par}
\expandafter\ifx\theemail\relax    
\relax\else\vglue 5pt plus 0.02fil minus 2pt\def\\{\stdspace{\rm 
and}\stdspace} 
\cl{Email:\stdspace\tt\theemail}\fi
\ifx\theurl\relax                  
\relax\else\vglue 5pt plus 0.02fil minus 2pt\def\\{\stdspace{\rm 
and}\stdspace}
\cl{URL:\stdspace\tt\theurl}\fi\par}

\vglue 7pt plus 0.3fil minus 3pt

{\bf Abstract}
\vglue 5pt plus 0.1fil minus 2pt

\theabstract

\vglue 7pt plus 0.3fil minus 3pt

{\bf AMS Classification numbers}\quad Primary:\quad \theprimaryclass

Secondary:\quad \thesecondaryclass

\vglue 5pt plus 0.3fil minus 2pt

{\bf Keywords}\quad \thekeywords

\vglue 10pt plus 0.5fil minus 5pt

{\small  Proposed: \theproposer\hfill Received: \receiveddate\nl
Seconded: \theseconders\hfill 
\ifx\reviseddate\relax                         
Accepted: \accepteddate                        
\else
Revised: \reviseddate                          
\fi}
\eject
}       

\let\maketitlepage\maketitlep
\let\maketitle\maketitlepage

\font\phead=cmsl9 scaled 950
\font=cmsl9 scaled 1050
\font\pnum=cmbx10 scaled 913
\font\lnum=cmbx10 
\font\pfoot=cmsl9 scaled 950
\font=cmsl9 scaled 1050
\ifplaintex
\headline{\vbox to 0pt{\vskip -4.5mm\line{\small\phead\ifnum
\count0=\startpage ISSN 1364-0380 (on line)
1465-3060 (printed) \hfill {\pnum\folio}\else\ifodd\count0\def\\{ }%
\ifx\theshorttitle\relax\thetitle\else\theshorttitle\fi\hfill{\pnum\folio}
\else\def\\{ and }{\pnum\folio}\hfill\ifx\theshortauthors\relax\theauthors
\else\theshortauthors\fi\fi\fi}\vss}}
\footline{\vbox to 0pt{\vglue 0mm\line{\small\pfoot\ifnum\count0=\startpage
\copyright\ \gtp\hfill\else
\gt, Volume \thevolumenumber\ (\thevolumeyear)\hfill\fi}\vss
}}
\else
\makeatletter
\def\@oddhead{{\small\lhead\ifnum\count0=\startpage ISSN 1364-0380 (on line)
1465-3060 (printed) \hfill {\lnum\number\count0}\else\ifodd\count0
\def\\{ }\ifx\theshorttitle\relax \thetitle \else\theshorttitle\fi\hfill
{\lnum\number\count0}\else\def\\{ and }{\lnum\number\count0}
\hfill\ifx\theshortauthors\relax 
\theauthors\else\theshortauthors\fi\fi\fi}}\def\@evenhead{\@oddhead}
\def\@oddfoot{\small\lfoot\ifnum\count0=\startpage\copyright\ \gtp\hfill\else
\gt, Volume \thevolumenumber\ (\thevolumeyear)\hfill\fi}
\def\@evenfoot{\@oddfoot}
\makeatother
\fi

\newwrite\gtoutfile
\long\gdef\makeheadfile{  
{\def\\{, }\def\s{ }
\immediate\openout\gtoutfile head.xxx
\immediate\write\gtoutfile{To: math@arxiv.org}
\immediate\write\gtoutfile{Subject: put or rep NNNNN:pppp}
\immediate\write\gtoutfile{--text follows this line--}
\immediate\write\gtoutfile{Proxy-for: \ifx\theasciiauthors\relax
\theauthors\else\theasciiauthors\fi\s<\ifx\theasciiemail\relax\theemail\else\theasciiemail\fi>}
\immediate\write\gtoutfile{\noexpand\\}
\immediate\write\gtoutfile{Authors: \ifx\theasciiauthors\relax
\theauthors\else\theasciiauthors\fi}
{\def\\{ }\immediate\write\gtoutfile{Title: \ifx\theasciititle\relax
\thetitle\else\theasciititle\fi}}
\immediate\write\gtoutfile{Subj-class: GT or SG or MG etc}
\immediate\write\gtoutfile{MSC-class: \theprimaryclass\ifx\thesecondaryclass\relax\else, \thesecondaryclass\fi}
\immediate\write\gtoutfile{Journal-ref: Geom. Topol. \thevolumenumber
(\thevolumeyear) \startpage-\finishpage}
\immediate\write\gtoutfile{Comments: Published by Geometry and Topology at}
\immediate\write\gtoutfile{\s\s http://www.maths.warwick.ac.uk/gt/GTVol\thevolumenumber/paper\thepapernumber.abs.html}
\immediate\write\gtoutfile{\noexpand\\}
\immediate\write\gtoutfile{}
\ifx\theasciiabstract\relax
\immediate\write\gtoutfile{\theabstract}\else
\immediate\write\gtoutfile{\theasciiabstract}\fi
\immediate\write\gtoutfile{}
\immediate\write\gtoutfile{\noexpand\\}
\immediate\write\gtoutfile{}
\immediate\closeout\gtoutfile}}  

\def\maketitlepage{\maketitlep\makeheadfile}
\let\maketitle\maketitlepage